\documentclass{mhkr}
\usepackage[style = numeric, sorting = nyt]{biblatex}
\usepackage{bussproofs}
\usepackage{longtable}
\usepackage{changepage}
\usetikzlibrary{cd}

\newenvironment{mathprooftree}
  {\varwidth{.9\textwidth}\centering\leavevmode}
  {\DisplayProof\endvarwidth}

\title{Proof-theoretic dilator and intermediate pointclasses}
\author{Hanul Jeon}
\email{ \href{mailto:hj344@cornell.edu}{hj344@cornell.edu}}
\urladdr{ \href{https://hanuljeon95.github.io}{https://hanuljeon95.github.io} }
\address{Department of Mathematics, Cornell University, Ithaca, NY 14853} 
\thanks{The research presented in this paper is supported in part by NSF grant DMS–2153975. The author would like to thank James Walsh for helpful comments and Patrick Lutz for a helpful answer to my question about iterated hyperjumps. \added{The author's dissertation \cite{JeonPhD} contains an older version of the results of this paper, with an incorrect proof of $\sfbeta$-completeness.}}

\newcommand{\lag}{\langle}
\newcommand{\rag}{\rangle}
\newcommand{\lr}{\leftrightarrow}

\newcommand{\rank}{\operatorname{rank}}
\newcommand{\field}{\operatorname{field}}

\newcommand{\supp}{\operatorname{supp}}
\newcommand{\Occ}{\operatorname{Occ}}
\newcommand{\Den}{\operatorname{Den}}

\newcommand{\Ord}{\mathrm{Ord}}

\newcommand{\WO}{\mathsf{WO}}
\newcommand{\LO}{\mathsf{LO}}
\newcommand{\Clim}{\mathsf{Clim}}

\newcommand{\Dil}{\mathsf{Dil}}

\newcommand{\CK}{\mathsf{CK}}
\newcommand{\KB}{\mathsf{KB}}

\newcommand{\Id}{\mathsf{Id}}
\newcommand{\Emb}{\mathsf{Emb}}
\newcommand{\RFN}[1][]{\ifthenelse{\equal{#1}{}}{}{#1\mhyphen}\mathsf{RFN}}

\newcommand{\AC}{\mathsf{AC}}

\newcommand{\CA}{\mathsf{CA}}
\newcommand{\RCA}{\mathsf{RCA}}
\newcommand{\ACA}{\mathsf{ACA}}
\newcommand{\ATR}{\mathsf{ATR}}

\newcommand{\TI}{\mathsf{TI}}
\newcommand{\KP}{\mathsf{KP}}

\newcommand{\ZFC}{\mathsf{ZFC}}

\newcommand{\ID}{\mathsf{ID}}

\newcommand{\HJ}{\mathsf{HJ}}
\newcommand{\Lin}{\operatorname{Lin}}

\newcommand{\proj}{\operatorname{proj}}
\newcommand{\Spec}{\operatorname{Spec}}
\newcommand{\SP}{\mathsf{SP}}
\newcommand{\PA}{\mathsf{PA}}

\newcommand{\Th}{\operatorname{Th}}

\newcommand{\LK}{\mathsf{LK}}
\newcommand{\ulc}{\ulcorner}
\newcommand{\urc}{\urcorner}
\newcommand{\Alt}{\mathsf{Alt}}
\newcommand{\Imp}{\mathsf{Imp}}
\newcommand{\RecLO}{\mathsf{RecLO}}

\usepackage[T1]{fontenc}
\DeclareRobustCommand{\msf}[1]{%
  \ifcat\noexpand#1\relax\msfgreek{#1}\else\mathsf{#1}\fi
}
\makeatletter
\newcommand{\msfgreek}[1]{\csname sf\expandafter\@gobble\string#1\endcsname}
\makeatother
\DeclareFontFamily{T1}{CBG-sf}{}
\DeclareFontShape{T1}{CBG-sf}{m}{n}{<-> gsmn1000}{}
\def\GreekFontCBGsf{\usefont{T1}{CBG-sf}{m}{n}}

\DeclareFontFamily{T1}{CBG-sfu}{}
\DeclareFontShape{T1}{CBG-sfu}{m}{n}{<-> gsmu1000}{}
\def\GreekFontCBGsf{\usefont{T1}{CBG-sf}{m}{n}}

\newcommand{\sfbeta}{{\text{\GreekFontCBGsf \char"62}}}

\newcommand{\Seq}{\operatorname{Seq}}
\newcommand{\sfsup}{\mathsf{sup}}
\newcommand{\id}{\mathsf{id}}
\makeatletter
\newcommand{\make@circled}[2]{%
  \ooalign{$\m@th#1\smallbigcirc{#1}$\cr\hidewidth$\m@th#1#2$\hidewidth\cr}%
}
\newcommand{\smallbigcirc}[1]{%
  \vcenter{\hbox{\scalebox{0.77778}{$\m@th#1\bigcirc$}}}%
}
\newcommand{\owedge}{\mathbin{\mathpalette\make@circled{\raisebox{0.02em}{\scalebox{0.89}{$\wedge$}}}}}
\newcommand{\ovee}{\mathbin{\mathpalette\make@circled{\raisebox{-0.02em}{\scalebox{0.89}{$\vee$}}}}}
\makeatother

\newcommand{\added}[1]{#1}
\newcommand{\modified}[1]{#1}

\begin{document}

\begin{abstract}
    There are two major generalizations of the standard ordinal analysis: One is Girard's $\Pi^1_2$-proof theory in which dilators are assigned to theories instead of ordinals. The other is Pohlers' generalized ordinal analysis with Spector classes, where ordinals greater than $\omega_1^\CK$ are assigned to theories. 
    In this paper, we show that these two are systematically entangled, and $\Sigma^1_2$-proof theoretic analysis has a critical role in connecting these two.
\end{abstract}


\maketitle

\section{Introduction}

Proof theory begins from Hilbert's program, an attempt to secure the consistency of all of mathematics by formalizing mathematics and proving the consistency of the formalized mathematics by finitary means. Although G\"odel's incompleteness theorem showed Hilbert's goal is unachievable in its form, Gentzen's proof-theoretic analysis of Peano Arithmetic $\PA$ also attested we can reduce the consistency of $\PA$ to the well-foundedness of $\varepsilon_0 = \sup \{\omega,\omega^\omega,\cdots\}$.
Logicians divide the concepts when they confront paradoxes by circularity. For example, Russell's paradox tells us we cannot form the set of all sets, and we have to separate sets and proper classes.
Likewise, G\"odel's incompleteness theorem also shows there is no single theory proving its own consistency, and we need to separate theories by their \emph{consistency strength}, resulting in the hierarchy of theories by their consistency strength.
We can understand $\varepsilon_0$ as an ordinal representing the strength of Peano Arithmetic, and more formally, \emph{$\varepsilon_0$ is the proof-theoretic ordinal of $\PA$.}

Gentzen's proof-theoretic analysis is the beginning of \emph{Ordinal analysis}. The main goal of ordinal analysis is to calculate the \emph{proof-theoretic ordinal}
\begin{equation*}
    |T|_{\Pi^1_1} = \sup \{\operatorname{otp}(\alpha)\mid \text{$\alpha$ is a recursive well-order such that }T\vdash \text{``$\alpha$ is well-founded''}\}
\end{equation*}
of a given theory $T$. The assertion ``A recursive linear order is well-founded'' is $\Pi^1_1$-complete, that is, every $\Pi^1_1$-statement is equivalent to an assertion of this form. Hence, the proof-theoretic ordinal of a given theory characterizes $\Pi^1_1$-consequences of the theory. 
Calculating the proof-theoretic ordinal of a theory becomes extremely hard as the target theory becomes stronger, and defining the corresponding proof-theoretic ordinal requires a good insight into transfinite objects. One of the reasons is that sentences of higher complexity interact with the proof-theoretic ordinal in a non-trivial manner. There are various ways to overcome this issue, and Girard's \emph{$\Pi^1_2$-proof theory} provides a way.

Girard developed $\Pi^1_2$-proof theory to take a finistic control on infinitary objects \cite{Girard1985IntroductionPi12} by separating the `well-founded part' and the `indiscernible part' \cite{GirardRessayre1985}. The `well-founded part' is represented by \emph{dilators}.
One way to understand a dilator is by viewing it as a \emph{denotation system}: Consider the case expressing the class well-order $\Ord+\Ord$. There is no transitive class isomorphic to $\Ord+\Ord$, but we can still represent it as a collection $\{(i,\xi)\mid i=0,1\lor \xi\in\Ord\}$ with the lexicographic order.
Let us observe that the same construction still works even when we replace $\Ord$ with an arbitrary well-order $\alpha$ or even a linear order.
That is, the map
\begin{equation*}
    X \mapsto \{(i,\xi)\mid i=0,1\lor \xi\in\Ord\} \text{ with the lexicographic order}
\end{equation*}
for a linear order $X$ gives a way to form a new linear order $X+X$.
We can also think of the previous example as we `separated' $\Ord+\Ord$ into the `indiscernible part' $\Ord$ (which can be any linear order, in fact) and the `well-founded part' that instructs how to construct $X+X$ from $X$.
However, not every denotation system is a dilator, as if not every linear order is a well-order: We want to ensure a dilator preserves well-foundedness, so if $D$ is a dilator and $X$ is a well-order, then $D(X)$ is also a well-order.
Thus, we introduce intermediate notions called a \emph{semidilator} and a \emph{predilator}. 
Girard proved that every $\Pi^1_2$-statement is equivalent to `a recursive predilator $D$ is a dilator' for some recursive predilator $D$ (cf., \autoref{Theorem: Girard completeness theorem for dilators}), so dilators are the right object representing $\Pi^1_2$-statements. 
 
Like ordinal analysis provides proof-theoretic ordinal, dilator-based proof-theoretic analysis yields \emph{proof-theoretic dilator} $|T|_{\Pi^1_2}$ of $T$ characterizing the $\Pi^1_2$-consequences of $T$.
Although Girard's $\Pi^1_2$-proof theory provides a tool to analyze theories at the level of $\Pi^1_1\mhyphen\CA_0$ and iterated inductive definitions, it has been underrepresented in proof theory since the mid-1990s.  However, this has changed recently, and $\Pi^1_2$-proof theory and related topics are actively investigated by various people.

$|T|_{\Pi^1_2}$ can be thought as a function from $\Ord$ to $\Ord$, and we can ask the following question:
Does $|T|_{\Pi^1_1}(\alpha)$ encode any proof-theoretic information about $T$? If $\alpha$ is recursive, the following is known:
\begin{theorem}[Aguilera and Pakhomov, {\cite[Theorem 9]{AguileraPakhomov2023Pi12}}] \pushQED{\qed}
    Let $T$ be a $\Pi^1_2$-sound theory extending $\ACA_0$ and $\alpha$ be a recursive well-order. Then $|T|_{\Pi^1_2}(\alpha) = \lvert T + \WO(\alpha)\rvert_{\Pi^1_1}$. \qedhere 
\end{theorem}
Then what is $|T|_{\Pi^1_2}(\alpha)$ for a \emph{non-recursive} $\alpha$? The main goal of this paper is to answer the question for $\alpha<\delta^1_2$, where $\delta^1_2$ is the supremum of ordertypes of $\Delta^1_2$-wellorders (equivalently, the least ordinal satisfying $L_{\delta^1_2}\prec_{\Sigma_1} L$.)
Before presenting the answer, let us illustrate some recent results in proof theory that will motivate us.
Not many proof-theoretic dilators of canonical theories have been computed, but the proof-theoretic dilator of $\ACA_0$ is known: 
\begin{theorem}[Aguilera and Pakhomov \cite{AguileraPakhomov??Pi12PTA}] \pushQED{\qed}
    $\lvert\ACA_0\rvert_{\Pi^1_2} = \varepsilon^+$, where $\varepsilon^+$ is a dilator such that $\varepsilon^+(\alpha)$ is the least epsilon number greater than $\alpha$. \qedhere 
\end{theorem}

\modified{
Meanwhile, Pohlers \cite{Pohlers2015Semiformal, Pohlers2022Performance} considered ordinal analysis over structures other than the natural numbers, assigning non-recursive ordinals to represent the strength of a theory relative to the given structure.
To address the motivation behind Pohlers' framework, let us examine Pohlers' structural conception of ordinal analysis. Polhers stated \cite[p.59]{Pohlers1992Shortcourse} that determining the proof-theoretic ordinal of $T$ is closely related to determining the least model for $\Sigma_1^{L_{\omega_1^\CK}}$-sentences provable from $T$. One related fact about it is the Spector-Gandy theorem, which states that for a $\Pi^1_1$-sentence $\phi$, there is a $\Sigma_1$-sentence $\psi$ in the language of set theory such that
\begin{equation*}
    \bbN \vDash \phi \iff \exists \alpha<\omega_1^\CK \bigl(L_\alpha\vDash \psi\bigr).
\end{equation*}
For sufficiently weak theories, predicative cut elimination suffices to catch the proof-theoretic ordinal.
However, for a sufficiently strong theory $T$ in the sense that $\omega$-models of $T$ include $L_{\omega_1^\CK}$,%
\footnote{$\KP$ is an example of such a theory: Suppose that $M\vDash \KP$ is an $\omega$-model (so $\omega=\omega^M$). If $M$ is well-founded, then clearly the transitive collapse of $M$ contains $L_{\omega_1^\CK}$. If $M$ is ill-founded, then Ville's lemma implies $L_{\omega_1^\CK}$ is contained in a well-founded part of $M$.}
$T$ can see how outer structures of $L_{\omega_1^\CK}$ affects its $\Sigma_1^{L_{\omega_1^\CK}}$-consequences.
This appears in the form of an ordinal collapsing function in impredicative cut elimination. In particular, ordinal analysis for `impredicative' theories requires analyzing the behavior of theories over `outer structures' of $L_{\omega_1^\CK}$.
We want to analyze how the outer structures affect $L_{\omega_1^\CK}$, and we should stratify the outer structures for the analysis.
The natural way to stratify the outer structures is employing admissible sets, like $L_\alpha$ for admissible ordinals $\alpha>\omega_1^\CK$.
Admissible sets are set-theoretic objects, and we want to have an arithmetical counterpart to fit them into the arithmetical framework.
It turns out that the notion of a \emph{Spector class} is such a counterpart:
In a nutshell, a Spector class is a collection of subsets of $\bbN$ having a recursion-theoretic nature. For a Spector class $\Gamma$, its \emph{companion}\footnote{Roughly speaking, the companion of $\Gamma$ for a Spector class $\Gamma$ is a collection of sets that can be coded by $\Gamma\cap\check{\Gamma}$-relations.}
is an admissible set. (See \cite[Ch. 9]{Moschovakis1974EIAS} or \cite[Ch. 6--7]{Moellerfeld2002Thesis} for the details.) Also, the collection of lightface $\Pi^1_1$-sets forms a Spector class.
By viewing the ordinal analysis of $T$ as a procedure computing how many true $\Pi^1_1$-statements $T$ can prove, we can generalize ordinal analysis by replacing $\Pi^1_1$ with a Spector class $\Gamma$, so the generalized ordinal analysis computes how many facts about $\Gamma$ are captured by $T$.

Now, let us illustrate a general framework that Pohlers formulated.
Pohlers suggested gauging the performance of a theory $T$ by gauging the gap between the supremum $\delta^\frakM$ of $\frakM$-definable well-orders and the supremum $\delta^\frakM(T)$ of $T$-provably well-founded $\frakM$-definable well-orders, like the gap between $\omega_1^\CK$ and the proof-theoretic ordinal $|T|_{\Pi^1_1}$ describes the strength of $T$. Indeed, $\delta^\frakM(T)$ is the usual proof-theoretic ordinal $|T|_{\Pi^1_1}$ if we take $\frakM$ to be the standard structure $\bbN$ of natural numbers.
Informally, larger $\frakM$ can be seen as a `measuring stick' with a coarser scale: On the one hand, calculating $\delta^\frakM(T)$ for a larger $\frakM$ is no harder than for a smaller $\frakM$. On the other hand, however, $\delta^\frakM(T)$ for a large $\frakM$ `ignores' more consequences of lower complexity as $\frakM$ grows; Pohlers' No Enhancement theorem \cite[Theorem 3.21]{Pohlers2022Performance} tells $\delta^\frakM(T) = \delta^\frakM(T + \sigma)$ for a $\Sigma^1_1$-statement $\sigma$ with parameters from $\frakM$. Facts about Spector classes (like \cite[\S.9F]{Moschovakis1974EIAS}) suggest $\delta^\frakM(T)$ gauges the degree of $T$ capturing facts about the `next Spector class' above $\frakM$.}

The structure $\frakM$ can be arbitrary in principle; Most structures we consider are an expansion of $\bbN$ under certain conditions. In particular, Pohlers focused on the \emph{iterated Spector classes} $\SP^\xi_\bbN$ to gauge the performance of iterated inductive definitions $\ID_{<\alpha}$.
Notably, we have the following:
\begin{theorem}[Pohlers {\cite[Remark 3.35]{Pohlers2022Performance}}]
    For $\mu\ge 1$ less than the least recursively inaccessible ordinal, we have
    \begin{equation*}
        \delta^{\SP^{\mu}_\bbN}(\Th(\SP^{\mu}_\bbN) + \text{Mathematical Induction}) = \varepsilon_{\omega^{\CK}_\mu+1}.
    \end{equation*}
\end{theorem}

In this paper, we will see that Spector's ordinal spectrum is a special case of a \emph{$\Pi^1_1[R]$ proof-theoretic ordinal} for a $\Sigma^1_2$-singleton $R$.%
\footnote{A real $R$ is a $\Sigma^1_2$-singleton if the set $\{R\}$ is $\Sigma^1_2$-definable. We also call such $R$ an \emph{$\Sigma^1_2$-singleton real}.}
The ($\Pi^1_1$-)proof-theoretic ordinal
\begin{equation*}
    |T|_{\Pi^1_1} = \sup \{\operatorname{otp}(\alpha) \mid \text{$\alpha$ is a recursive linear order such that }T\vdash\WO(\alpha)\}
\end{equation*}
uses the ordertype of $T$-provably recursive well-orders to gauge the strength of $T$. 
Likewise, $\Pi^1_1[R]$ proof-theoretic ordinal $|T|_{\Pi^1_1[R]}$ uses the ordertype of $T$-provably \emph{$R$-recursive} well-orders to gauge the strength of $T$. (See \autoref{Subsection: Pi 1 1 R proof theory} for the precise definition of $|T|_{\Pi^1_1[R]}$.)
We will see in \autoref{Theorem: Successor Spector class as a Pi 1 1 R pointclass} that the iterated Spector class $\SP^{\xi+1}_\bbN$ for $\xi$ less than the least recursively inaccessible ordinal is equivalent to $(\bbN;R)$ for some $\Sigma^1_2$-singleton $R$, corresponding to the $\xi$th iterate of the hyperjump of $\emptyset$.
Under the perspective of \autoref{Proposition: Connecting Pohlers' ordinal and Pi 1 1 R PTO}, Pohlers' result can be rephrased as follows:
\begin{theorem} \pushQED{\qed} \label{Theorem: Pohlers' higher ordinal analysis}
    Let $\xi$ be a successor ordinal less than the least recursively inaccessible. Then 
    \begin{equation*}
        |\ACA_0 + \text{$\HJ^\xi(\emptyset)$ exists} + \Th(\bbN;\HJ^\xi(\emptyset))|_{\Pi^1_1[\HJ^\xi(\emptyset)]} = \varepsilon_{\omega_\xi^\CK+1}. \qedhere 
    \end{equation*}
    Here $\Th(\bbN;\HJ^\xi(\emptyset))$ is the collection of true \emph{first-order} sentences over $(\bbN;\HJ^\xi(\emptyset))$.
\end{theorem}

Pohlers' result and the $\Pi^1_2$-proof theoretic analysis by Aguilera and Pakhomov hint at the connection between these two. 
The next theorem, which is the main result of this paper, illustrates how they are systematically entangled:
\begin{theorem*}[\autoref{Theorem: Main theorem}]
    Let $T$ be a $\Pi^1_2$-sound theory extending $\ACA_0$ and $(D,\varrho)$ be a recursive \emph{locally well-founded genedendron} generating $R$. Furthermore, assume that $T$ proves $(D,\varrho)$ is a locally well-founded genedendron.
    If $\alpha$ is an $R$-recursive well-order such that $D_\alpha$ is ill-founded, then 
    \begin{equation*}
        |T|_{\Pi^1_2}(\alpha) = |T[R] + \WO(\alpha)|_{\Pi^1_1[R]}.
    \end{equation*}
    Here $T[R]$ is the theory $T$ + `$R$ exists.'
\end{theorem*}
The statement of the main theorem uses a new concept named \emph{genedendron}. Less formally, a genedendron is a funtorial family of trees $\lag D_\alpha\mid \alpha\in\Ord\rag$ with a constant partial function $\varrho_\alpha$ taking an infinite branch of $D_\alpha$ and returning a real. The value of $R$ does not depend on the choice of $\alpha$ and an infinite branch. We say a genedendron \emph{generates $R$} if the value of $\varrho_\alpha$ is $R$, and the least ordinal $\alpha$ making $D_\alpha$ ill-founded represents the complexity of the real $R$. We will provide an extended discussion on a genedendron in \autoref{Subsection: What is a genedendron} and \autoref{Section: Genedendrons}.
We will calculate a genedendron $(D,\varrho)$ generating $\HJ(\emptyset)$ such that $D_\alpha$ is ill-founded for a $\HJ(\emptyset)$-recursive well-order of ordertype $\omega_1^\CK$. Hence we will have
\begin{equation*}
    |\ACA_0 + \text{$\HJ(\emptyset)$ exists}|_{\Pi^1_1[\HJ(\emptyset)]} = \varepsilon_{\omega_1^\CK+1}.
\end{equation*}
Combining it with \autoref{Proposition: Adding true Sigma 1 1 R sentence does not change the Pi 1 1 R PTO}, we have
\begin{equation*}
    |\ACA_0 + \text{$\HJ(\emptyset)$ exists} + \Th(\bbN;\HJ(\emptyset))|_{\Pi^1_1[\HJ(\emptyset)]} = \varepsilon_{\omega_1^\CK+1}.
\end{equation*}
which partially reproduces Pohlers' result. It shows the connection between the proof-theoretic dilator and the proof-theoretic-ordinal for \emph{intermediate pointclass} of the form $\Pi^1_1[R]$ for a $\Sigma^1_2$-singleton $R$.%
\footnote{If we line up pointclasses with the prewellordering property under the supremum of the ordertype the pointclass possesses, then $\Pi^1_1[R]$ lies between $\Pi^1_1$ and $\Sigma^1_2$. It is why we call $\Pi^1_1[R]$ an intermediate pointclass.}

\subsection{What is a genedendron?} \label{Subsection: What is a genedendron}
The reader may notice that the statement of the main theorem used a new terminology \emph{genedendron}. We will present its formal definition in \autoref{Section: Genedendrons}, but let us introduce its informal explanation. 

Let us start with the following question: \emph{What is the relationship between $\HJ^\xi(\emptyset)$ and $\omega_\xi^\CK$?}
One possible answer is the following: $\omega_\xi^\CK$ describes the \emph{complexity} of a real $\HJ^\xi(\emptyset)$. Then what is the precise meaning of complexity?
If we view $\HJ^\xi(\emptyset)$ as a $\sfbeta$-model of arithmetic, then $\HJ^\xi(\emptyset)$ has height $\omega_\xi^\CK$.
This viewpoint works for the iterated hyperjumps of the empty set but does not work in general since not every real takes the form of a $\sfbeta$-model of arithmetic.
Despite that, the keyword $\sfbeta$-model will shed light on the correct way to define complexity for $\Sigma^1_2$-singleton reals since we will use $\sfbeta$-logic to define the correct notion of complexity.

Turning our viewpoint to the recursion-theoretic side, Suzuki \cite{Suzuki1964Delta12} (see also \cite[\S V.6]{Hinman1978}) defined a notion of rank for $\Sigma^1_2$-singleton reals. Although Suzuki's notion is sufficient to establish properties of $\Delta^1_2$ sets (e.g., \cite[Theorem 6.8, 6.9]{Hinman1978}), this rank is difficult to calculate by hand. Despite that, the way Suzuki defined the rank is worth examining: 
Following terminologies in \cite[\S V.6]{Hinman1978}, for each $\Sigma^1_2$-formula $\exists^1 Y\forall^1 Z \phi(X,Y,Z)$ \added{witnessing $R$ being a $\Sigma^1_2$-singleton}, we associate a \emph{matrix} $F_\phi$ for the real $R$. 
The matrix $F_\phi$ is a well-order, so we can define the rank of $R$ by the least ordertype of $F_\phi$ for a $\Sigma^1_2$-formula $\phi$ \added{witnessing $R$ being a $\Sigma^1_2$-singleton}.

Suzuki's rank is hard to compute for a specific real since a matrix is hard to calculate. However, this definition still illustrates the idea that defining the rank follows from extracting an ordinal from an implicit $\Sigma^1_2$-definition of a real.
Suzuki's definition shows a similarity with a framework of $\Sigma^1_2$-ordinal analysis presented in \cite{Jeon2024HigherProofTheoryI}: $\Sigma^1_2$-sentences are represented by a recursive \emph{pseudodilator}, and a pseudodilator is associated with its \emph{climax}, an ordinal less than $\delta^1_2$.
However, we need more from what \cite{Jeon2024HigherProofTheoryI} presented for the following reasons: First, we have not provided a way to calculate a pseudodilator from a given $\Sigma^1_2$-formula.
Second, we do not have a way to extract a real from a pseudodilator whose corresponding $\Sigma^1_2$-formula implicitly defines the real.

We overcome these two issues \added{via} proof-theoretic means
To illustrate how it works, let us examine a proof-theoretic proof for the Kleene normal form theorem, stating that for every $\Pi^1_1$-formula $\phi$, we can compute a recursive linear order $\alpha$ such that $\phi\lr \WO(\alpha)$.
Let $\phi\equiv \forall^1X \psi(X)$ for an arithmetical formula $\psi(X)$. Then consider the  $\omega$-logic presented in \cite[\S6]{Girard1987Logical1} with an extra unary symbol $\overline{X}$.
The $\omega$-logic enjoys the following property we will call the \emph{preproof property}:
\begin{theorem}[Girard {\cite[Theorem 6.1.13]{Girard1987Logical1}}] \pushQED{\qed}
    Let $\calL$ be a language of a recursive $\omega$-theory $T$ and let $\Gamma\vdash\Delta$ be a sequent in $\calL$.
    Then we can  construct a recursive $\omega$-preproof $\pi$ for $\Gamma\vdash\Delta$ such that $\bigwedge\Gamma \to\bigvee \Delta$ holds over every $\omega$-model of $T$ iff $\pi$ is well-founded. 
    Furthermore, if $\pi$ is ill-founded with an infinite branch $B$, then there is a $B$-recursive $\omega$-model of $T+\lnot\left(\bigwedge\Gamma \to\bigvee \Delta\right)$.\footnote{\added{Girard used two-sided sequent calculus to formulate $\sfbeta$-logic. We will use Tait-style calculus to formulate $\sfbeta$-logic.}}

    In particular, if $\phi$ is a formula in the language of first-order arithmetic, then $\phi$ is true iff the $\omega$-preproof $\pi$ for $\vdash\phi$ is well-founded. \qedhere 
\end{theorem}

Let $\pi$ be an $\omega$-preproof for the sequent $\vdash\psi(\overline{X})$ given by the preproof property.
Then $\forall^1 X \phi(X)$ holds iff $\pi$ is well-founded. Furthermore, if $\pi$ has an infinite branch, then we can construct an $\omega$-model $M=(\bbN,\overline{X}^M)$ from the infinite branch such that $M\vDash \lnot \phi(\overline{X})$.
In particular, if $\forall^1 X\phi(X)$ fails so $\pi$ has an infinite branch, then we can find a real $R$ satisfying $\lnot\phi(R)$ by examining sequents of the form $\vdash \overline{X}(\overline{n})$ and $\overline{X}(\overline{n})\vdash$.

Observe from the previous explanation that the construction of the $\omega$-preproof $\pi$ is an important step for the proof of the $\omega$-completeness theorem (cf. \cite[Theorem 6.1.12]{Girard1987Logical1}.) Similarly, we may use a proof of $\sfbeta$-completeness theorem (See \autoref{Subsection: completeness of beta logic}) to construct a predilator $D$ from a $\Pi^1_2$-sentence $\phi$ such that $\phi\lr \Dil(D)$.
Furthermore, the resulting predilator $D$ comes from a functorial $\sfbeta$-preproof $\lag \pi_\alpha\mid\alpha\in\Ord\rag$.
If $\phi$ fails, then we have the least $\alpha$ such that $\pi_\alpha$ has an infinite branch, and we should be able to extract a real $R$ such that $\lnot\psi(R)$, where $\psi(X)$ is a $\Sigma^1_1$-formula satisfying $\phi\equiv \forall^1 X \psi(X)$.
The ordinal $\alpha$ should be read as a complexity of the real $R$.\footnote{One may wonder that the previous argument only gives a way to define an ordinal complexity for $\Pi^1_1$-singleton reals. However, a slight modification gives a way to compute an ordinal complexity for $\Sigma^1_2$-singleton reals. In computability-theoretic view, every $\Sigma^1_2$-singleton real is Turing reducible to a $\Pi^1_1$-singleton real (See \autoref{Proposition: Implicit Sigma 1 2 reduces to Implicit Pi 1 1}), so their computability-theoretic properties are not too different.}

In particular, if $\lag \pi_\alpha\mid\alpha\in\Ord\rag$ is a $\sfbeta$-preproof for a $\Pi^1_2$-formula $\forall^1 X\psi(X)$ such that $\lnot \psi(X)$ implicitly defines a real, then there is the least $\alpha$ making $\pi_\alpha$ ill-founded, and we can understand $\alpha$ as a complexity of a real implicitly defined by $\lnot\psi(X)$.
It captures the idea of the definition of a \emph{genedendron} we will introduce in \autoref{Section: Genedendrons}.
Since $\sfbeta$-preproofs are \emph{prequasidendroids} that we will introduce in \autoref{Section: Quasidendroids}, we will define genedendrons in terms of prequasidendroids. 
Genedendrons give the right notion of a complexity for $\Sigma^1_2$-singleton reals called the \emph{$\Sigma^1_2$-altitude} (See \autoref{Definition: Sigma 1 2 altitude}), and we will see in \autoref{Subsection: Hyperjump 0 - Sigma 1 2 alt} that the $\Sigma^1_2$-altitude of $\HJ(\emptyset)$ is indeed $\omega_1^\CK$ by proving a form of $\sfbeta$-completeness theorem for a certain $\sfbeta$-system for Peano arithmetic with additional predicates for Kleene's $\calO$.
Just as we use cut-elimination of infinitary systems to obtain $\Pi^1_1$- or $\Pi^1_2$-consequences of a theory, we use the preproof property for a $\sfbeta$-system to extract information about a $\Sigma^1_2$-singleton real. In this sense, we may view the computation in \autoref{Subsection: Hyperjump 0 - Sigma 1 2 alt} as a form of $\Sigma^1_2$-proof theoretic analysis.

\tableofcontents

\section{Preliminaries} \label{Section: Prelims}
In this section, we review preliminary knowledge necessary for this paper except for quasidendroids and $\sfbeta$-logic, which will appear in later sections separately. We also include a summary of Pohlers' work \cite{Pohlers2022Performance} to support arguments in the Introduction, although its technical details are unnecessary for our work.

In this paper, we usually denote natural numbers in lowercase Latin letters and real numbers in uppercase Latin letters. However, lowercase Greek letters do not necessarily mean they are natural numbers; They could mean ordinals, finite sequents, or reals. Some lowercase Latin letters, such as $f$ and $g$, may represent functions rather than natural numbers.
$\forall^0$ and $\exists^0$ mean quantifications over natural numbers, and $\forall^1$ and $\exists^1$ mean those over real numbers.

\modified{
There are two notions of definability of reals: For a complexity class $\Gamma$, we say a real $R$ is \index{Gamma@$\Gamma$-definable} \emph{$\Gamma$-definable} if there is a $\Gamma$-formula $\phi(m)$ such that $R = \{n\in\bbN\mid \phi(m)\}$.
A real $R$ is \index{Gamma@$\Gamma$-singleton}\emph{$\Gamma$-singleton} (or \index{Implicitly $\Gamma$-definable}\emph{implicitly $\Gamma$-definable}) if there is a $\Gamma$-formula $\phi(X)$ such that $R$ is the unique real satisfying $\phi(R)$. We mostly use $\Gamma$-singleton reals in this paper.
}

\subsection{\texorpdfstring{$\Sigma^1_2$-\added{singleton}}{Sigma 1 2-singleton} reals}
In recursion-theoretic context, $\Sigma^1_2$-singleton real is `reducible' to a \added{$\Pi^1_1$-singleton real} in the following manner:
\begin{proposition}[$\Pi^1_1\mhyphen\CA_0$] \label{Proposition: Implicit Sigma 1 2 reduces to Implicit Pi 1 1}
    If $A$ is a $\Sigma^1_2$-singleton real, then there is a \added{$\Pi^1_1$-singleton real} $B$ such that $A$ is primitive recursive in $B$.
\end{proposition}
\begin{proof}
    Suppose that $\exists^1 Z\phi(Z,X)$ is a $\Sigma^1_2$-formula \added{witnessing $A$ being a $\Sigma^1_2$-singleton real} for a $\Pi^1_1$-formula $\phi$. By $\Pi^1_1$-uniformization theorem, we can find a $\Pi^1_1$-formula $\hat{\phi}$ uniformizing $\phi$. Then consider the $\Pi^1_1$-formula
    \begin{equation*}
        \psi(Y) \equiv \phi((Y)_0,(Y)_1),
    \end{equation*}
    where $(Y)_i = \{n\in\bbN\mid \lag i,n\rag\in Y\}$. If $\hat{\phi}(B,A)$ holds for some real $B$, then $B\oplus A$ is a \added{singleton real witnessed by} $\psi$. Clearly, $A$ is primitive recursive in $B\oplus A$.
\end{proof}

The previous theorem says there is no significant recursion-theoretic difference between $\Sigma^1_2$-singleton reals and \added{$\Pi^1_1$-singleton} reals. Thus, examining the recursion-theoretic behavior for \added{$\Pi^1_1$-singleton} reals is sufficient for examining that for $\Sigma^1_2$-singleton reals.
However, we stick to $\Sigma^1_2$-definable reals since \added{$\Sigma^1_2$-definable reals are more proof-theoretically natural objects.}
Suzuki proved that every $\Delta^1_2$-real is hyperarithmetical in a \added{$\Pi^1_1$-singleton}  real, hinting \added{at} the connection between $\Sigma^1_2$-statements and the intermediate pointclasses. The proof is not hard, so let us include its proof for completeness.
\begin{theorem}[Suzuki \cite{Suzuki1964Delta12}, $\Sigma^1_2\mhyphen\AC_0$]
    For $X\subseteq \bbN^k$, $X$ is $\Delta^1_2$ iff there is a \added{$\Pi^1_1$-singleton} real $R$ such that $X\in \Delta^1_1[R]$.
\end{theorem}
\begin{proof}
    Suppose that $\phi$ is a $\Pi^1_1$-formula witnessing $R$ is \added{a $\Pi^1_1$-singleton} , so $R$ is the unique real satisfying $\phi(R)$. If $\psi$ is a $\Sigma^1_1$-formula satisfying
    \begin{equation*}
        \vec{n}\in X \iff \psi(\vec{n},R),
    \end{equation*}
    for every $\vec{n}\in \bbN^k$, then we have
    \begin{equation*}
        \vec{n}\in X\iff \added{\exists^1} Y[\phi(Y)\land \psi(\vec{n},Y)],
    \end{equation*}
    so $X$ is $\Sigma^1_2$. By applying the same argument to $\lnot X$, we can derive $X$ is $\Delta^1_2$.
    Conversely, suppose that $X\in \Delta^1_2$, and let $\chi(\vec{n},i)$ be a $\Sigma^1_2$-formula defining the characteristic function for $X$:
    \begin{equation*}
        [\vec{n}\in X \lr \chi(\vec{n},1)] \land [\vec{n}\notin X \lr \chi(\vec{n},0)].
    \end{equation*}
    Then we have
    \begin{equation*}
        Z=X\iff \forall^0 \vec{n} [(\vec{n}\in Z\to \chi(\vec{n},1)) \land (\vec{n}\notin Z\to \chi(\vec{n},0))].
    \end{equation*}
    Now let $\chi(\vec{n},i)\equiv \exists^1 Y\psi(\vec{n},i,Y)$ for some $\Pi^1_1$-formula $\psi$. Then by $\Sigma^1_2\mhyphen \AC_0$, we have
    \begin{equation*}
        Z=X \iff \exists^1 Y \forall^0\vec{n} \bigl[\bigl(\vec{n}\in Z\to \psi(\vec{n},1,(Y)_{\vec{n}})\bigr) \land \bigl(\vec{n}\notin Z\to \psi(\vec{n},0,(Y)_{\vec{n}})\bigr)\bigr].
    \end{equation*}
    The right-hand side has the form $\exists^1 Y \phi(Y)$ for a $\Pi^1_1$-formula $\phi$. By $\Pi^1_1$-uniformization that is provable over $\Pi^1_1\mhyphen\CA_0$, we can find a $\Pi^1_1$-uniformization $\hat{\phi}$ of $\phi$. Now let $R$ be the unique real satisfying $\hat{\phi}(R)$. Then we have
    \begin{equation*}
        Z=X \iff \forall^0\vec{n} \bigl[\bigl(\vec{n}\in Z\to \psi(\vec{n},1,(R)_{\vec{n}})\bigr) \land \bigl(\vec{n}\notin Z\to \psi(\vec{n},0,(R)_{\vec{n}})\bigr)\bigr].
    \end{equation*}
    Hence
    \begin{equation*}
        \vec{n}\in X \iff \psi\bigl(\vec{n},1,(R)_{\vec{n}}\bigr) \iff \lnot \psi\bigl(\vec{n},0,(R)_{\vec{n}}\bigr),
    \end{equation*}
    so $X$ is $\Delta^1_1[R]$.
\end{proof}

But can we get a result showing a more proof-theoretic connection between $\Sigma^1_2$-sentences and the intermediate pointclasses? 
The following lemma says every $\Sigma^1_2$-\emph{sentence} is $\Sigma^0_2[R]$-sentence for some implicitly $\Pi^1_1$-definable real $R$:
\begin{lemma}[$\Pi^1_1\mhyphen \CA_0$] \label{Lemma: Sigma 1 2 as Sigma 0 2 with a real}
    If $\sigma$ is a $\Sigma^1_2$-sentence, then we can find a \added{$\Pi^1_1$-singleton} real $R$ such that $\sigma$ is equivalent to a $\Sigma^0_2[R]$-sentence.
\end{lemma}
\begin{proof}
    Suppose that $\sigma\equiv \exists X \phi(X)$ for a $\Pi^1_1$-formula $\phi$. Then we have
    \begin{equation*}
        \exists^1 X\phi(X)\iff \exists^1 M [\text{$M$ is a $\beta$-model}\land M\vDash \exists^1 X\phi(X)].
    \end{equation*}
    The right-hand-side of the above formula takes the form $\exists^1 M \psi(M)$ for a $\Pi^1_1$-formula $\psi$. Let $\hat{\psi}$ be the uniformization of $\psi$, and let $R$ be a real satisfying $\hat{\psi}(R)$. Then we have 
    \begin{equation} \label{Formula: Equivalence-Sigma12-Sigma02[R]}
        \exists^1 X\phi(X)\iff R\vDash \exists^1 X\phi(X),
    \end{equation}
    and the right-hand-side is $\Sigma^0_2[R]$.
\end{proof}

\autoref{Lemma: Sigma 1 2 as Sigma 0 2 with a real} itself has no practical information since $R$ exists only when $\sigma$ holds, so the equivalence \eqref{Formula: Equivalence-Sigma12-Sigma02[R]} does not make sense if $\sigma$ fails. 
However, \autoref{Lemma: Sigma 1 2 as Sigma 0 2 with a real} illustrates the possibility of decomposing $\Pi^1_2$ into $\Pi^0_2[R]$ or $\Pi^1_1[R]$ for a \added{$\Pi^1_1$-singleton} real $R$, which motivates the connection between proof-theoretic \added{dilators} and $R$-recursive well-orders.

\begin{definition}
    Suppose that $T$ is a theory and $\phi(X)$ is a $\Sigma^1_2$-formula \added{witnessing} $R$ \added{being a $\Sigma^1_2$-singleton}. Furthermore, suppose that $T$ proves $\phi(X)$ holds for at most one real; That is,
    \begin{equation*}
        T\vdash \forall^1 X,Y [\phi(X)\land \phi(Y)\to X=Y].
    \end{equation*}
    Then let us define $T[R] := T + (\exists^1 X \phi(X))$.
\end{definition}
The reason we assume the uniqueness assertion
\begin{equation*}
    \forall^1 X,Y [\phi(X)\land\phi(Y)\to X=Y]
\end{equation*}
is a $T$-theorem that is usually provable over a weak theory like $\ACA_0$. Also relatedly, in the context of our proposed main theorem, we can add a `true $\Sigma^1_2$-sentence' to the theory $T$ by looking at $|T|_{\Pi^1_2}(\alpha)$ for a sufficiently large $\alpha$, but adding a true $\Pi^1_2$-sentence is impossible in this manner. 

\subsection{\texorpdfstring{$\Pi^1_1[R]$}{Π¹₁[R]}-proof theory} \label{Subsection: Pi 1 1 R proof theory}
In this subsection, we define the proof-theoretic ordinal augmented with a $\Sigma^1_2$-singleton $R$.
We want to make this notion syntactic, so instead of using a semantic notion of a real $R$, we use its $\Sigma^1_2$-definition to formulate $\Gamma[R]$-formulas and $\Pi^1_1[R]$ proof-theoretic ordinal.

\begin{definition} \label{Definition: Gamma r hierarchy} 
    Let $\phi(X)$ be a $\Sigma^1_2$-formula \added{witnessing} $R$ \added{being a $\Sigma^1_2$-singleton}.
    For a lightface class $\Gamma$ of formulas, the \emph{$\Gamma[R]$-formula} is a formula of the form $\forall^1 X [\phi(X)\to \psi(X)]$ for a $\Gamma$-formula $\psi$.
\end{definition}
Note that if $T$ proves $\exists!^1 X\phi(X)$, then $T$ proves $\forall^1 X [\phi(X)\to \psi(X)]$ is equivalent to $\exists^1 X [\phi(X)\land\psi(X)]$.

\begin{definition}  \label{Definition: Pi 1 1 r PTO}
    Let $T$ be a theory proving $\exists!^1 X\phi(X)$ for a $\Sigma^1_2$-formula $\phi(X)$  \added{witnessing} $R$ \added{being a $\Sigma^1_2$-singleton}.
    We say $T$ is \emph{$\Gamma[R]$-sound} if every $T$-provable $\Gamma[R]$-sentence is true; That is, if $T\vdash \forall^1 X[\phi(X)\to\psi(X)]$ for a $\Gamma$-formula $\psi(X)$, then $\psi(R)$ holds.
    
    For a $\Pi^1_1[R]$-sound theory $T$ such that $T\vdash \exists!^1 X\phi(X)$, let us define
    \begin{equation*} 
        |T|_{\Pi^1_1[R]} := \sup\{\operatorname{otp}(\alpha) \mid \text{$\alpha$ is $R$-recursive and }T\vdash \WO(\alpha)\}.
    \end{equation*}
    However, the expression $\WO(\alpha)$ is not expressible in the language of second-order arithmetic.
    Hence we rephrase the above definition as follows: Every $R$-recursive linear order takes the form $\alpha'(R)$ for some recursive function $\alpha'$. 
    Thus we can take
    \begin{multline*}
        |T|_{\Pi^1_1[R]} := \sup \bigl\{\operatorname{otp}(\alpha'(R)) \bigm| \text{$\alpha'$ is a recursive function and } \\ T\vdash \forall^1 X [\phi(X)\to \WO(\alpha'(X))]\bigr\}.
    \end{multline*}
    The above definition works when $T$ is $\Pi^1_1[R]$-sound, otherwise $\alpha'(R)$ may not even be a linear order for some $\alpha'$.
\end{definition}

From the remaining part of this subsection, let us check the properties of $\Pi^1_1[R]$-proof-theoretic ordinal. \emph{Throughout this subsection, we always assume that $\phi(X)$ is a $\Sigma^1_2$-formula \added{witnessing} $R$ \added{being a $\Sigma^1_2$-singleton real} and $T$ proves $\exists!^1 X \phi(X)$.}

Let us recall the precise statement of the Kleene normal form theorem:
\begin{theorem}[Kleene, $\ACA_0$] \pushQED{\qed}
    Let \added{$\chi(X,n)$} be a $\Pi^1_1$-formula. Then we can effectively find a primitive recursive linear order $\alpha(X,n)$ such that $\added{\chi(X,n)}\lr \WO(\alpha(X,n))$ for every real $X$ and a natural number $n$. \qedhere 
\end{theorem}

For later use, let us prove the following version of the fixed point theorem:
\begin{lemma} \label{Lemma: Kleene Fixed point theorem}
    For a natural number $k$, we can find another natural number $e$ such that
    \begin{equation*}
        \ACA_0 \vdash \forall^1 X \forall^0 n [\{\overline{e}\}^X(n) \simeq \{\{\overline{k}\}^X(\overline{e})\}^X(n)].
    \end{equation*}
\end{lemma}
\begin{proof}
    By $s$-$m$-$n$-theorem for oracles, we can find a primitive recursive function $h$ (which does not depend on the choice of $X$) such that 
    \begin{equation*}
        \ACA_0\vdash \forall^1 X\forall^0 m,n \bigl[\{h(m)\}^X(n) \simeq \{\{m\}^X(m)\}^X(n)\bigr].
    \end{equation*}
    Let us find $m$ such that $\{k\}^X\circ h$ is $\{m\}^X$, uniformly on $X$.
    Furthermore, we have
    \begin{equation*}
        \ACA_0 \vdash \forall^1 X \forall^0 n[\{\overline{k}\}^X(h(n))\simeq \{\overline{m}\}^X(n)].
    \end{equation*}
    Then we have
    \begin{equation*}
        \ACA_0 \vdash \forall^1 X\forall^0 n[
        \{h(\overline{m})\}^X(n) \simeq 
        \{\{\overline{m}\}^X(\overline{m})\}^X(n) \simeq
        \{\{\overline{k}\}^X(h(\overline{m}))\}^X(n)].
    \end{equation*}
    Thus $m=h(e)$ is the desired natural number.
\end{proof}

The following effective version of $\Sigma^1_1[R]$-boundedness theorem has an important role in establishing properties of $\Pi^1_1[R]$-proof-theoretic ordinal. Its proof follows from a variation of that by Rathjen \cite[Lemma 1.1]{Rathjen1991ParametersinBRBI}, and let us define a preliminary notion for the proof:
\begin{definition}
    Consider the $\Pi^0_2[X]$-formula
    \begin{equation*}
        \RecLO(X,e) \equiv \text{The relation $\prec^X_e := \{\lag i,j\rag\mid \{e\}^X(\lag i,j\rag) = 0\}$ is a linear order.}
    \end{equation*}
\end{definition}
Let us recall that the linear order $\alpha(X,n)$ in the statement of the Kleene normal form theorem is primitive recursive, meaning that there is a primitive recursive function $f$ such that $f(X,n,\cdot)$ computes the characteristic function for $\alpha(X,n)$. Hence we can find a natural number $e$, which only depends on $\alpha$ and not on $X$ and $n$, such that 
\begin{equation} \label{Formula: Index for a p.r. linear order}
    \forall^0 i,j [\{e\}^X(n)(\lag i,j\rag) = 0 \iff i <_{\alpha(X,n)} j].
\end{equation}
This fact is also provable \added{in} $\ACA_0$.

\begin{lemma}[Effective {$\Sigma^1_1[R]$}-boundedness, $\ACA_0$] \label{Lemma: Effective Sigma 11R boundedness}
    Let $T$ be a theory extending $\ACA_0$.
    Suppose that $\psi(X,n)$ is a $\Sigma^1_1[X]$-formula such that
    \begin{equation} \label{Formula: Statement-Sigma 11R bddness 00}
        T \vdash \forall^1 R \bigl[\phi(R) \to \bigl[\forall^0 n [\psi(R,n)\to \RecLO(R,n)\land \WO(\prec^R_n)]\bigr]\bigr].
    \end{equation}
    (Informally, it says $T$ proves $\{n\in\bbN\mid \psi(R,n)\}$ is a set of indices for $R$-recursive well-orders of the form $\prec^R_n$.) 
    Then we can find a natural number $e$ such that
    \begin{equation*}
        T \vdash \forall^1 R \phi(R)\to [\RecLO(R,\overline{e}) \land \WO(\prec^R_{\overline{e}}) \land \lnot\psi(R, \overline{e})].
    \end{equation*}
\end{lemma}
\begin{proof}
    By Kleene normal form theorem and the previous discussion about \eqref{Formula: Index for a p.r. linear order}, we can find a natural number $k$ that only depends on $\psi$ such that
    \begin{equation} \label{Formula: Pf-Sigma 11R bddness 00}
        \ACA_0\vdash\forall^1 X\forall^0 n\ \RecLO(X,\{\overline{k}\}^X(n)) \land \bigl[\lnot\psi(X,n) \lr \added{\WO(\prec^\calX_{\{\overline{k}\}^X(n)})}\bigr].
    \end{equation}
    By \autoref{Lemma: Kleene Fixed point theorem}, we can find a natural number $e$ such that
    \begin{equation*}
        \ACA_0\vdash\forall^1 X \forall^0 n \bigl[\{\overline{e}\}^X(n) \simeq \bigl\{\{\overline{k}\}^X(\overline{e})\bigr\}^X(n)\bigr].
    \end{equation*}
    We claim that $e$ is the desired natural number: \emph{Reasoning over $T$}, observe that $\{\overline{e}\}^X$ and $\{\{\overline{k}\}^X(\overline{e})\}$ are extensionally the same, so $\RecLO(X,\overline{e})$ iff $\RecLO(X,\{\overline{k}\}^X(\overline{e}))$, and the latter holds.
    Hence $\RecLO(X,\overline{e})$ holds for every choice of $X$.

    Now we claim that $\WO(\prec^R_{\overline{e}})$ holds. Suppose the contrary assume that $\lnot\WO(\prec^R_{\overline{e}})$ holds, which is equivalent to $\lnot\WO\bigl(\prec^R_{\{\overline{k}\}^X(\overline{e})}\bigr)$.
    Then by \eqref{Formula: Pf-Sigma 11R bddness 00}, we have $\psi(R,\overline{e})$.
    However, we assumed that \eqref{Formula: Statement-Sigma 11R bddness 00} holds, so $\WO(\prec^R_{\overline{e}})$, a contradiction.
    $\lnot\psi(R,\overline{e})$ follows from a part of the previous argument.
\end{proof}

Then we can prove the following, which has essentially the same proof in \cite{Rathjen1999Realm}:
\begin{proposition} \label{Proposition: Pi 1 1 R PTO and arithmetical in R WO length}
    Let $T$ be a $\Pi^1_1[R]$-sound extension of $\ACA_0+ \added{\exists!^1 X\phi(X)}$. Then 
    \begin{equation*}
        |T|_{\Pi^1_1[R]} = \sup\{\operatorname{otp}(\alpha)\mid \alpha\text{ is arithmetical-in-$R$ well-order such that }T\vdash\WO(\alpha)\}.
    \end{equation*}
\end{proposition}
\begin{proof}
    Let us clarify the meaning of arithmetical-in-$R$ linear order in our context: $\alpha$ is arithmetical-in-$R$ means there is a first-order formula $\psi(n,X)$ such that $\alpha=\{n\mid \psi(n,R)\}$.
    Then the statement `$\alpha$ is a linear order' is also arithmetical-in-$R$, whose formal statement is
    \begin{multline*}
        \forall^1 R\bigl[\phi(R)\to \forall^0 i,j,k \psi(\lag i,i\rag, R)\land \bigl[\psi(\lag i,j\rag),R)\lor \psi(\lag j,i\rag, R)\bigr] \\  \land \bigl[\psi(\lag i,j\rag,R) \land \psi(\lag j,k\rag,R)\to\psi(\lag i,k\rag,R)\bigr] \bigr].
    \end{multline*}
    Similarly, the statement `$\alpha$ is a well-order' is $\Pi^1_1[R]$.
    
    Let $\Emb(\alpha,\beta)$ be the $\Sigma^1_1$-statement saying `There is an embedding from a linear order $\alpha$ to another linear order $\beta$.'
    Now let $\alpha(R)$ be an arithmetical-in-$R$ linear order over $\bbN$ such that $T\vdash \WO(\alpha(R))$.
    Then observe that 
    \begin{equation*}
        T \vdash \forall^1 R \phi(R)\to \bigl[ \forall^0 n\bigl[\RecLO(R,n)\land \Emb\bigl(\prec^R_n,\alpha(R)\bigr) \to \WO(\prec^R_n)\bigr]\bigr],
    \end{equation*}
    so \autoref{Lemma: Effective Sigma 11R boundedness} applied to the formula $\psi(R,n)\equiv \RecLO(R,n)\land \Emb(\prec^R_n,\alpha(R))$, we can find a natural number $e$ such that
    \begin{equation*}
        T \vdash \forall^1 R \phi(R)\to [\RecLO(R,\overline{e})\land \WO(\prec^R_{\overline{e}}) \land \lnot \Emb(\prec^R_{\overline{e}},\alpha(R))].
    \end{equation*}
    The formula in the right-hand-side of \added{the turnstile} is $\Pi^1_1[R]$. Hence by $\Pi^1_1[R]$-soundness of $T$, $\forall^1 R \phi(R)\to \lnot \Emb(\prec^R_{\overline{e}},\alpha(R))$ is a true $\Pi^1_1[R]$-statement.
    Hence $\operatorname{otp}(\alpha(R)) \le \operatorname{otp}(\prec^R_e) \le |T|_{\Pi^1_1[R]}$.\footnote{The second inequality is strict since we can prove $\alpha<|T|_{\Pi^1_1[R]}$ implies $\alpha+1<|T|_{\Pi^1_1[R]}$.}
\end{proof}

\begin{proposition}
    Let $T$ be a $\Pi^1_1[R]$-sound $\Sigma^1_1[R]$-definable theory. Then $|T|_{\Pi^1_1[R]}<\delta(R)$, where $\delta(R)$ is the supremum of \added{ the ordertype of $R$-recursive well-orders}.
\end{proposition}
\begin{proof}
    Let 
    \begin{equation*}
        X=\{e\in \bbN \mid \RecLO(R,e)\land T\vdash \ulc \forall^1 Y [\phi(Y)\to \WO(\prec^Y_{\overline{e}})] \urc\}.
    \end{equation*}
    Then $X$ is $\Sigma^1_1[R]$, so the usual $\Sigma^1_1[R]$-boundedness theorem implies $\sup\{\operatorname{otp}(\prec^R_e)\mid e\in X\} < \delta(R)$.
\end{proof}

We prove Kriesel's theorem for $\Pi^1_1[R]$-proof-theoretic ordinal, which states adding the true $\Sigma^1_1[R]$-sentence to a theory does not change the $\Pi^1_1[R]$-proof-theoretic ordinal of the theory. 
\cite[Proposition 2.24]{Rathjen1999Realm} uses L\"ob's theorem to prove the given proposition for $R=\emptyset$. We use a different way akin to the first equivalence in \cite{Walsh2022incompleteness}, and the following variation of \cite[Lemma 5.2]{PakhomovWalsh2023omegaReflection} has a focal role in the proof:
\begin{lemma}[$\ACA_0$] \label{Lemma: Disjunction of two linear orders} 
    Let $\alpha(X)$ and $\beta(X)$ be recursive functions such that both of them are linear orders for every real $X$. Then we can define a new recursive function $(\alpha\lor\beta)(X)$ such that for every real $X$,
    \begin{enumerate}
        \item $(\alpha\lor\beta)(X)$ is a linear order.
        \item $\WO((\alpha\lor\beta)(X)) \lr \WO(\alpha(X)) \lor \WO(\beta(X))$.
        \item If $\beta(X)$ is ill-founded, then there is an embedding $\alpha(X)\to(\alpha\lor\beta)(X)$.
    \end{enumerate}
\end{lemma}
\begin{proof}
    Let $T(X)$ be the set of sequences $\lag (a_i,b_i)\mid i<n\rag$ such that each of $\lag a_i\mid i<n\rag$ and $\lag b_i\mid i<n\rag$ are decreasing sequences over $\alpha(X)$ and $\beta(X)$ respectively.
    Since both $\alpha$ and $\beta$ are recursive functions, $T$ is also a recursive function, and $T(X)$ is an $X$-recursive tree over $\bbN\times\bbN$. (We understand $\alpha(X)$ and $\beta(X)$ have the field a subset of $\bbN$, which is possible since both of them are $X$-recursive.)
    Then consider the linear order $(\alpha\lor\beta)(X)$ defined by $(T(X),\le_\KB)$; that is, 
    \begin{equation*}
        \lag (a_i,b_i) \mid i<n\rag <_\KB \lag (a'_i,b'_i) \mid i<n'\rag
    \end{equation*}
    if and only if either
    \begin{enumerate}
        \item $\lag (a_i,b_i) \mid i<n\rag$ is a proper initial segment of $\lag (a'_i,b'_i) \mid i<n'\rag$, or
        \item There is $i<n$ such that $(a_i,b_i) < (a'_i,b'_i)$, where the comparison is done under the lexicographic order over $\bbN\times\bbN$.
    \end{enumerate}
    $\le_\KB$ is $\Delta^0_1$-definable, so $(\alpha\lor\beta)$ is a recursive function. The remaining part of the proof is identical to that of \cite[Lemma 5.2]{PakhomovWalsh2023omegaReflection}, so we omit the proof.
\end{proof}

\begin{proposition} \label{Proposition: Adding true Sigma 1 1 R sentence does not change the Pi 1 1 R PTO}
    Let $T$ be a $\Pi^1_1[R]$-sound extension of $\ACA_0+ \exists!^1 X\phi(X)$. Then $|T|_{\Pi^1_1[R]} = |T+\sigma|_{\Pi^1_1[R]}$ for a true $\Sigma^1_1[R]$-formula $\sigma$.
\end{proposition}
\begin{proof}

    Let $\sigma$ be $\exists^1 X [\phi(X)\to \psi(X)]$ for some $\Sigma^1_1$-formula $\psi(X)$ with the only free variable $X$. By the Kleene normal form theorem, we have a primitive recursive $\alpha(X)$ such that
    \begin{equation*}
        \ACA_0\vdash \forall^1 X [\lnot\psi(X) \lr \WO(\alpha(X))].
    \end{equation*}
    Suppose that we have a primitive recursive well-order $\beta(X)$ such that
    \begin{equation*}
        T + \sigma \vdash \forall^1 X [\phi(X)\to \WO(\beta(X))].
    \end{equation*}
    Recall that every $\Pi^1_1[R]$-sentence has the form $\forall^1 X [\phi(X)\to \WO(\beta(X))]$ for some primitive recursive linear order $\beta$ by the Kleene normal form theorem. Then we have
    \begin{equation*}
        T \vdash \lnot \sigma \lor \forall^1 X [\phi(X)\to \WO(\beta(X))],
    \end{equation*}
    which is equivalent to
    \begin{equation*}
        T \vdash \forall^1 X \bigl[\phi(X)\to [\WO(\alpha(X))\lor \WO(\beta(X))]\bigr].
    \end{equation*}
    By \autoref{Lemma: Disjunction of two linear orders}, we have
    \begin{equation*}
        T \vdash \forall^1 X \bigl[\phi(X)\to \WO\bigl((\alpha\lor\beta)(X)\bigr)\bigr].
    \end{equation*}
    By $\Pi^1_1[R]$-soundness of $T$, the sentence $\forall^1 X \bigl[\phi(X)\to \WO\bigl((\alpha\lor\beta)(X)\bigr)\bigr]$ is true. Hence we have $\WO\bigl((\alpha\lor\beta)(R)\bigr)$. Also, the definition of $|T|_{\Pi^1_1[R]}$ implies $\operatorname{otp}\bigl((\alpha\lor\beta)(R)\bigr) < |T|_{\Pi^1_1[R]}$.
    Since $\sigma$ is true, $\psi(R)$ is also true, which is equivalent to $\lnot \WO(\alpha(R))$. This implies $\beta(R)$ embeds to $(\alpha\lor\beta)(R)$, so $\operatorname{otp}(\beta(R))\le \operatorname{otp}\bigl((\alpha\lor\beta)(R)\bigr) < |T|_{\Pi^1_1[R]}$.
\end{proof}

\subsection{Characterizing \texorpdfstring{$\Pi^1_1[R]$}{Π¹₁[R]}-proof theory} \label{Subsection: Characterizing Pi 1 1 R PT}
Walsh \cite{Walsh2023characterizations} proved that for an arithmetically definable $\Pi^1_1$-sound theories, the following are equivalent:
\begin{enumerate}
    \item Comparing their proof-theoretic ordinals.
    \item Comparing their $\Pi^1_1$-consequences modulo true $\Sigma^1_1$-statements.
    \item Comparing their $\Pi^1_1$-reflection, which is a $\Pi^1_1$-analogue of their consistency.
\end{enumerate}
We can establish a similar characterization to $\Pi^1_1[R]$-proof-theoretic ordinal for a $\Sigma^1_2$-\added{singleton} real $R$. Throughout this section, let us fix a $\Sigma^1_2$-formula $\phi(X)$ \added{witnessing} a real $R$ \added{being a $\Sigma^1_2$-singleton}.

\begin{definition}
    Let $S$, $T$ be two theories extending $\ACA_0 + \exists^1! X\phi(X)$. Furthermore, assume that they are $\Pi^1_1[R]$-sound. 
    For a sentence $\sigma$, define $T\vdash^{\Sigma^1_1[R]}\sigma$ iff there is a $\Sigma^1_1$-formula $\psi(X)$ with only free variable $X$ such that 
    \begin{equation*}
        T + (\exists^1 X \phi(X)\land\psi(X)) \vdash \sigma.
    \end{equation*}
    Also, define $S \subseteq^{\Sigma^1_1[R]}_{\Pi^1_1[R]} T$ if and only if for every $\Pi^1_1$-formula $\psi(X)$ with only free variable $X$, 
    \begin{equation*}
         S\vdash^{\Sigma^1_1[R]} \forall^1 X(\phi(X)\to\psi(X)) \implies T\vdash^{\Sigma^1_1[R]} \forall^1 X(\phi(X)\to\psi(X)).
    \end{equation*}
\end{definition}

Then we can prove the following:
\begin{theorem} \label{Theorem: First equivalence for Pi 1 1 R}
    Let $S$, $T$ be $\Pi^1_1[R]$-sound extensions of $\ACA_0 + \exists^1!X\phi(X)$. Then
    \begin{equation*}
        S \subseteq^{\Sigma^1_1[R]}_{\Pi^1_1[R]} T \iff |S|_{\Pi^1_1[R]} \le |T|_{\Pi^1_1[R]}.
    \end{equation*}
\end{theorem}

\added{\autoref{Theorem: First equivalence for Pi 1 1 R} follows from the following lemma:}
\begin{lemma} \label{Lemma: for the First equivalence for Pi 1 1 R}
    Let $T$ be $\Pi^1_1[R]$-sound extensions of $\ACA_0 + \exists^1!X\phi(X)$. If $\alpha(X)$ is a recursive function, then
    \begin{equation*}
        T\vdash^{\Sigma^1_1[R]}\forall^1 X [\phi(X)\to\WO(\alpha(X))] \iff \alpha(R) < |T|_{\Pi^1_1[R]}.
    \end{equation*}
\end{lemma}
\begin{proof}
    In one direction, suppose that $T\vdash^{\Sigma^1_1[R]}\forall^1 X [\phi(X)\to\WO(\alpha(R))]$ holds.
    Let $\psi(X)$ be a $\Sigma^1_1$-formula such that
    \begin{equation} \label{Formula: First equivalence for Pi 1 1 R - 00}
        T + \exists^1 X [\phi(X)\land \psi(X)]\vdash\forall^1 X [\phi(X)\to\WO(\alpha(R))]
    \end{equation}
    
    By Kleene normal form theorem, we can find a recursive $\beta(X)$ such that 
    \begin{equation*}
        \ACA_0 \vdash \forall^1 X \bigl[ \psi(X)\lr \lnot\WO\bigl(\beta(X)\bigr)\bigr].
    \end{equation*}
    Hence from \eqref{Formula: First equivalence for Pi 1 1 R - 00}, we have
    \begin{equation*}
        T \vdash \forall^1 X \bigl[\phi(X)\to \bigl[\WO(\beta(X))\lor \WO(\alpha(X))\bigr]\bigr].
    \end{equation*}
    By \autoref{Lemma: Disjunction of two linear orders}, we can find a new recursive function $(\alpha\lor\beta)(X)$ such that 
    \begin{equation*}
        \ACA_0\vdash \forall^1 X \bigl[\WO((\alpha\lor\beta)(X)) \lr \WO(\alpha(X)\lor\beta(X))\bigr].
    \end{equation*}
    Hence $T\vdash \forall^1 X \bigl[\phi(X)\to\WO((\alpha\lor\beta)(X))\bigr]$, so $(\alpha\lor\beta)(R) < |T|_{\Pi^1_1[R]}$.\footnote{We can see that $|T|_{\Pi^1_1[R]}$ is always a limit ordinal, so the desired inequality follows.}
    However, $\beta(R)$ is ill-founded, so $\alpha(R)\le (\alpha\lor\beta)(R)$, which proves the desired conclusion.

    Conversely, assume that $\alpha(R) < |T|_{\Pi^1_1[R]}$.
    Then we can find a recursive $\beta(X)$ such that $\alpha(R)\le\beta(R)$ and $T\vdash \forall^1X[\phi(X)\to\WO(\beta(X))]$.
    $\alpha(R)\le\beta(R)$ implies the $\Sigma^1_1[R]$-statement
    \begin{equation*}
        \sigma\equiv \exists^1 X \bigl[\phi(X)\land \exists^1 F \bigl(\text{$F$ is an embedding from $\alpha(X)$ to $\beta(X)$}\bigr)\bigr]
    \end{equation*}
    is true.
    Then we can see that $T+\sigma$ proves $\forall^1 X[\phi(X)\to \WO(\alpha(X))]$, as desired.
\end{proof}

\begin{proof}[Proof of \autoref{Theorem: First equivalence for Pi 1 1 R}]
    Suppose that $S \subseteq^{\Sigma^1_1[R]}_{\Pi^1_1[R]} T$ holds.
    If $\xi<|S|_{\Pi^1_1[R]}$, then we can find a recursive $\alpha(X)$ such that $\xi<\alpha(R)$ and $S\vdash \forall^1 X[\phi(X)\to\WO(\alpha(X))]$.
    By the assumption, we have $T\vdash^{\Sigma^1_1[R]} \forall^1 X[\phi(X)\to\WO(\alpha(X))]$, so $\alpha(R) < |T|_{\Pi^1_1[R]}$ by \autoref{Lemma: for the First equivalence for Pi 1 1 R}.
    Hence $\xi \le \alpha(R) < |T|_{\Pi^1_1[R]}$.
    Since $\xi$ is arbitrary, we have $|S|_{\Pi^1_1[R]} \le |T|_{\Pi^1_1[R]}$.

    Conversely, if $|S|_{\Pi^1_1[R]} \le |T|_{\Pi^1_1[R]}$ and $S\vdash^{\Sigma^1_1[R]} \forall^1 X [\phi(X)\to\WO(\alpha(X))]$ for some recursive $\alpha(X)$, then $\alpha(R) < |T|_{\Pi^1_1[R]}$. This implies $T\vdash^{\Sigma^1_1[R]} \forall^1 X[\phi(X)\to\WO(\alpha(X))]$.
    Kleene normal form implies every $\Pi^1_1[R]$-sentence takes the form $\forall^1 X[\phi(X)\to\WO(\alpha(X))]$ for some recursive $\alpha(X)$, so we have $S \subseteq^{\Sigma^1_1[R]}_{\Pi^1_1[R]} T$.
\end{proof}

It is possible to prove the second equivalence for $\Pi^1_1[R]$-proof theory, whose proof follows from a modification of \cite[\S3.2]{Walsh2023characterizations} or \cite[\S6.2]{Jeon2024HigherProofTheoryI}. Hence \added{let me} leave its proof to the reader:
\begin{theorem}
    Let $S$, $T$ be arithmetically definable $\Pi^1_1[R]$-sound theories extending $\ACA_0 + \exists^1 X\phi(X)$. Then 
    \begin{equation*}
        |S|_{\Pi^1_1[R]} \le T_{\Pi^1_1[R]} \iff \ACA_0 + \exists^1! X\phi(X) \vdash^{\Sigma^1_1[R]} \RFN[{\Pi^1_1[R]}](T)\to \RFN[{\Pi^1_1[R]}](S),
    \end{equation*}
    where
    \begin{equation*}
        \RFN[{\Pi^1_1[R]}](T) \equiv \forall \psi(X)\in \Pi^1_1 [T\vdash \phi(X)\to \psi(X) \implies \vDash_{\Pi^1_1}\phi(R)\to\psi(R)].
    \end{equation*}
\end{theorem}

\subsection{Dilators} \label{Subsection: Dilator}
Dilators are central objects in Girard's $\Pi^1_2$-proof theory \cite{Girard1982Logical2}. 
There are several equivalent definitions of predilators and dilators. One of them is defining predilators and dilators as functors:
\begin{definition}
    A functor $F$ from the category of linear orders $\LO$ to $\LO$ is a \emph{semidilator} if $F$ preserves direct limits and pullbacks.
    A semidilator $F$ is a \emph{predilator} if it satisfies the \emph{monotonicity condition}: For every $f,g\colon X\to Y$,
    \begin{equation*}
        \forall x\in X[f(x)\le g(x)] \implies \forall \sigma\in F(X) [F(f)(\sigma)\le F(g)(\sigma)].
    \end{equation*}
    A semidilator $F$ is a \emph{dilator} if $F(X)$ is a well-order when $X$ is a well-order.
\end{definition}

Freund \cite{Freund2021Pi14} defined a semidilator\footnote{Freund used a word \emph{prae-dilator} to mean a semidilator. We follow terminologies in the author's previous work \cite{Jeon2024HigherProofTheoryI}.} as a functor $F\colon \LO\to\LO$ equipped with a support transformation $\supp\colon F\to [\cdot]^{<\omega}$ satisfying the \emph{support condition:} For every linear order $X, Y$ and an increasing $f\colon X\to Y$,
\begin{equation*}
    \{\sigma\in F(Y) \mid \supp_Y(\sigma)\subseteq \ran(f)\} \subseteq \ran(F(f)).
\end{equation*}

A third way to define a dilator is by viewing it as a denotation system, which is adopted in \cite{Jeon2024HigherProofTheoryI}. The denotation system definition of a semidilator identifies a semidilator $F$ with its \emph{trace}
\begin{equation*}
    \operatorname{Tr}(F) = \{(n,\sigma) \mid n\in\bbN\land \sigma\in F(n)\land \supp_n(\sigma)=n\}.
\end{equation*}
These three approaches are all equivalent (see \cite[\S 3]{Jeon2024HigherProofTheoryI} for its proof), and each approach has pros and cons:
Girard's original approach of dilators-as-functors is easy to state but has no technical advantage.
Freund's definition is convenient and useful in practice, especially when we check some object we constructed is a semidilator. However, Freund's definition is inconvenient when we want to view semidilators as a denotation system.
The denotation system approach is hard to state, and less useful when checking whether certain objects are semidilators. However, it becomes useful when we need to view semidilators as denotation systems, which is apparent in the abstract construction of predilators (see \cite{Girard1982Logical2} for the details.)
We use all of these approaches with no distinction.
The fourth approach for dilators is dilators as \emph{dendroids and quasidendroids}. We will not deeply address the dendroid approach, but quasidendroids are necessary to define $\sfbeta$-proofs. We address an excerpt of the fourth approach in \autoref{Section: Quasidendroids}.

The main reason for using dilators in proof theory is that dilators represent $\Pi^1_2$-statements. It is known as Girard's completeness theorem:
\begin{theorem}[{\cite[Theorem 8.E.1]{Girard1982Logical2}}, $\ACA_0$] \label{Theorem: Girard completeness theorem for dilators}
\pushQED{\qed}
    Let $\phi$ be a $\Pi^1_2$-sentence. Then we can find a primitive recursive predilator $D$ such that $\phi \lr \Dil(D)$ holds, where $\Dil(D)$ is an abbreviation for `$D$ is a dilator.'
    \qedhere
\end{theorem}

\subsection{\texorpdfstring{$\Pi^1_2$}{Pi 1 2}-proof theory}
The standard ordinal analysis (or $\Pi^1_1$-proof theory) computes the proof-theoretic ordinal of a theory. Similarly, Girard's $\Pi^1_2$-proof theory computes the \emph{proof-theoretic dilator} of a theory.
The notion of proof-theoretic dilator is rigorously defined in \cite{AguileraPakhomov2023Pi12} as follows:
\begin{definition}
    Let $T$ be a $\Pi^1_2$-sound extension of $\ACA_0$. Then $|T|_{\Pi^1_2}$ is the dilator unique up to bi-embeddability satisfying the following conditions:
    \begin{enumerate}
        \item\label{Item: Pi12dilator} For a recursive predilator $D$, if $T\vdash \Dil(D)$ then $D$ embeds to $|T|_{\Pi^1_2}$.
        \item (Universality) If a dilator $\hat{D}$ satisfies \eqref{Item: Pi12dilator}, then $\hat{D}$ embeds $|T|_{\Pi^1_2}$ and the following diagram commutes:
        \begin{equation*}
            \begin{tikzcd}
            	{D_0} \\
            	\vdots & {|T|_{\Pi^1_2}} & {\hat{D}} \\
            	{D_i}
            	\arrow[from=2-2, to=2-3]
            	\arrow[from=1-1, to=2-2]
            	\arrow[from=3-1, to=2-2]
            	\arrow[from=1-1, to=2-3, bend left=20]
            	\arrow[from=3-1, to=2-3, bend right=20]
            \end{tikzcd}
        \end{equation*}
        Here $\{D_i\mid i<\omega\}$ is an enumeration of all $T$-provably recursive dilators, i.e., a recursive dilator $D$ such that $T\vdash \Dil(D)$.
    \end{enumerate}
\end{definition}

Aguilera and Pakhomov \cite{AguileraPakhomov2023Pi12} proved that $|T|_{\Pi^1_2}$ is unique up to bi-embeddability if exists, and if $T$ is $\Pi^1_2$-sound extension of $\ACA_0$, then the ordered sum of all recursive dilators $D$ such that $T\vdash \Dil(D)$ is (bi-embeddable with) $|T|_{\Pi^1_2}$.

One way to calculate $|T|_{\Pi^1_1}$ follows from proving the cut elimination of an infinitary system associated with $T$.
Similarly, we can calculate $|T|_{\Pi^1_2}$ for a specific $T$ by proving the cut elimination argument of an associated $\sfbeta$-theory. In particular, calculating $|\ACA_0|_{\Pi^1_2}$ follows this strategy, which will appear in Aguilera-Pakhomov's unpublished work \cite{AguileraPakhomov??Pi12PTA}.

$\Pi^1_1[R]$-consequences of a theory for a $\Sigma^1_2$-singleton real $R$ is linearly comparable modulo true $\Sigma^1_1[R]$-consequences as proved in \autoref{Subsection: Characterizing Pi 1 1 R PT}.
However, \cite{AguileraPakhomov2024Nonlinearity} found that comparing the $\Pi^1_2$-consequences of two $\Pi^1_2$-sound theories modulo true $\Sigma^1_2$-sentences is not linear:

\begin{theorem}[Aguilera-Pakhomov \cite{AguileraPakhomov2024Nonlinearity}] \pushQED{\qed}
    There is a $\Pi^1_2$-sound recursive extension $T$ of $\ATR_0$ and two $\Pi^1_2$-sentences $\phi_0$, $\phi_1$ such that $T+\phi_i$ is $\Pi^1_2$-sound for $i=0,1$, but neither $T\vdash^{\Sigma^1_2}\phi_0\to\phi_1$ nor $T\vdash^{\Sigma^1_2}\phi_1\to\phi_0$ hold. \qedhere 
\end{theorem}

\subsection{\texorpdfstring{$\Sigma^1_2$}{Sigma 1 2}-proof theory}
\autoref{Theorem: Girard completeness theorem for dilators} also entails that every $\Sigma^1_2$-sentence is equivalent to a statement of form $\lnot\Dil(D)$ for some primitive recursive predilator $D$. Predilators that are not dilators frequently appear when we analyze $\Sigma^1_2$-consequences of a theory, so we name them \emph{pseudodilators}. The ordinal information of a pseudodilator $F$ is the least ordinal $\alpha$ making $F(\alpha)$ ill-founded, which we call the \emph{climax of $F$}:
\begin{definition}
    A predilator $F$ is a \emph{pseudodilator} if $F$ is not a dilator. The \emph{climax of $F$} is the least ordinal $\alpha$ such that $F(\alpha)$ is ill-founded, and we call it $\Clim(F)$.
\end{definition}

Climaxes of pseudodilators can be used to define a $\Sigma^1_2$-version of the proof-theoretic ordinal of a theory:
\begin{definition}
    For a $\Sigma^1_2$-sound extension $T$ of $\ACA_0$, define
    \begin{equation*}
        s^1_2(T) = \sup\{\Clim(F) \mid \text{$F$ is a recursive predilator such that }T\vdash\lnot\Dil(F)\}.
    \end{equation*}
\end{definition}

Unlike the $\Pi^1_2$ case, $\Sigma^1_2$-consequences of a theory modulo true $\Pi^1_2$-sentences are linearly comparable. Even better, we have the following theorem, which is a $\Sigma^1_2$-analogue of Walsh's results in \cite{Walsh2023characterizations}:
\begin{theorem}[{\cite[\S 7]{Jeon2024HigherProofTheoryI}}] \pushQED{\qed}
    For two $\Sigma^1_2$-sound theories extending $\ACA_0$, the following holds:
    \begin{enumerate}
        \item Let us define $T\vdash^{\Pi^1_2} \phi$ if there is a true $\Pi^1_2$-sentence $\sigma$ such that $T+\sigma\vdash \phi$,
        and $S \subseteq^{\Pi^1_2}_{\Sigma^1_2} T$ if for every $\Sigma^1_2$-sentence $\phi$, $S\vdash^{\Pi^1_2} \phi \implies T\vdash^{\Pi^1_2} \phi$.
        Then we have
        \begin{equation*}
            s^1_2(S)\le s^1_2(T) \iff  S\subseteq^{\Pi^1_2}_{\Sigma^1_2} T.
        \end{equation*}
        \item Furthermore, if $S$ and $T$ are arithmetically definable $\Sigma^1_2$-sound extensions of $\Sigma^1_2\mhyphen\AC_0$, then we have
        \begin{equation*}
            s^1_2(S)\le s^1_2(T)\iff \Sigma^1_2\mhyphen \AC_0\vdash^{\Pi^1_2} \RFN[\Sigma^1_2](T)\to \RFN[\Sigma^1_2](S). \qedhere 
        \end{equation*}
    \end{enumerate}
\end{theorem}

\section{Pohlers' characteristic ordinals}
In this section, we review Pohlers' framework presented in \cite{Pohlers2022Performance}.
We may view that proof-theoretic ordinal gauges the strength of a theory over the structure $\bbN$ of natural numbers; By \autoref{Proposition: Pi 1 1 R PTO and arithmetical in R WO length}, $|T|_{\Pi^1_1}$ is equal to the supremum of ordertype of well-orders definable in $\bbN$ whose well-foundeness is a theorem of $T$. Pohlers generalized proof-theoretic ordinal by replacing $\bbN$ with an \emph{acceptable} countable structure:

\begin{definition}
    Let $\frakM$ be a countable structure in the language $\calL$. A structure is \emph{acceptable} if the following \emph{coding scheme} are first-order definable in $\frakM$:
    \begin{enumerate}
        \item A copy of $(\bbN,\le)$.
        \item A map $x_0,\cdots,x_{m-1}\mapsto \lag x_0,\cdots,x_{m-1}\rag$.
        \item The unary predicate $\mathrm{Seq}$ for the set of finite $\frakM$-sequences.
        \item A function $\operatorname{lh}$ returning a length of given finite sequence. (If $\lnot\mathrm{Seq}(x)$, $\operatorname{lh}(x)=0$.)
        \item A projection function $q(x,i)$ taking a finite sequence $x$ and $i\in\bbN$ (from the copy in $\frakM$) and returning the $i$th component of $x$. (Otherwise, return 0.)
    \end{enumerate}
    Then let us consider the new language $\calL(\frakM)$ augmented with the elements of $\frakM$ as constant symbols.
    An $\calL(\frakM)$-theory $T$ is an \emph{acceptable axiomatization of $\frakM$} if $T$ defines $\frakM$ and its coding scheme, $\frakM$-sound, and proves the following:%
    \footnote{Pohlers did not state the definition of acceptable axiomatization as a single definition. The following definition is a combination of his description in \cite[p38, Footnote 14, p50]{Pohlers2022Performance}.}
    \begin{enumerate}
        \item Every true atomic formula in the language $(\bbN,\le,\lag\cdot\rag,\mathrm{Seq},\operatorname{lh},q)$.
        \item Induction schema on the copy of $(\bbN,\le)$ for interpreted $\frakM$-formulas.
    \end{enumerate}
\end{definition}
Our choice of $T$ is first-order. However, we want to gauge $\Pi^1_1$-consequences of $T$ to gauge its strength \added{even when $T$ is a first-order arithmetic}. To cope with this situation, we assume $T$ has a unary predicate $X$ corresponding to a free set variable, and allow it to occur in the induction scheme. We call a sentence of the form $\phi(X)$ a \emph{pseudo $\Pi^1_1$-sentence.}

\begin{definition}[\added{Pohlers} {\cite[p48]{Pohlers2022Performance}}]
    A countable structure $\frakM$ over $\calL$ is \emph{strictly acceptable} if it is acceptable and 
    \begin{enumerate}
        \item There is an $\frakM$-definable well-order over $\frakM$.
        \item There is a first-order definable coding $\ulc F\urc$ for $\calL$-formulas with parameters in $\frakM$ such that $\ulc G\urc \in \operatorname{CS}(\ulc F\urc)$ is first-order expressible in $\frakM$.
        (See \cite[Definition 5.3.3]{Pohlers2009Prooftheory} or \cite[p40]{Pohlers2022Performance} for the definition of $\operatorname{CS}$.)
    \end{enumerate}
\end{definition}

\begin{definition}[\added{Pohlers} {\cite[Definition 3.15]{Pohlers2022Performance}}]
    For an acceptable structure $\frakM$ in the language $\calL$ and an acceptable axiomatization $T$ of $\frakM$, define
    \begin{itemize}
        \item $\delta^\frakM = \sup \{\operatorname{otp}(\prec) \mid \text{$\prec$ is definable in $\frakM$ with parameters and $\frakM\vDash \TI(\prec,X)$}\}$,
        \item $\delta^\frakM(T) = \sup \{\operatorname{otp}(\prec) \mid \text{$\prec$ is definable in $\frakM$ with parameters and $T\vdash \TI(\prec,X)$}\}$.
    \end{itemize}
    Here $X$ is a fresh predicate variable, and $\TI(\prec,X)$ states the transfinite induction along $\prec$ over $X$:
    \begin{equation*}
        \TI(\prec,X) \equiv \forall x[\forall y [y\prec x\to y\in X]\to x\in X]\to (\field(\prec)\subseteq X).
    \end{equation*}
\end{definition}

Pohlers proved for a strictly acceptable $\frakM$ and its acceptable axiomatization $T$ of $\frakM$, adding the first-order truth about $\frakM$ does not change the `strength' of $T$ in the sense that $\delta^\frakM(T) = \delta^\frakM(T+\Th(\frakM))$. (See \cite[Theorem 3.21]{Pohlers2022Performance}.) Then Pohlers argued that the only way to increase the strength of $T$ is by adding axioms for universes above $\frakM$. The extension of $\frakM$ what Pohlers considered is a \emph{Spector class}:
\begin{definition}
    A \emph{Spector class} $\Gamma$ is a collection of relations over $\frakM$ such that 
    \begin{enumerate}
        \item Every atomic predicate and function over $\frakM$, and their complements are in $\Gamma$. (For functions, consider their graph instead.)
        \item $\Gamma$ contains coding scheme over $\frakM$.
        \item $\Gamma$ is closed under $\cap$, $\cup$, $\exists^0$, $\forall^0$, and trivial combinatorial substitutions.\footnote{Trivial combinatorial substitution is a map that is a composition of projection maps and the tuple map; See \cite[p312]{Moschovakis2009}.}
        \item $\Gamma$ has a \emph{universal set}; That is, for each $n\in\bbN$ there is an $(n+1)$-ary relation $U\in \Gamma$ such that every $n$-ary $R\in\Gamma$ is a section of $U$.
        \item $\Gamma$ has the \emph{prewellordering property}; That is, for every $P\in\Gamma$ there is a \emph{norm} $\sigma_P\colon P\to\Ord$ such that the relations
        \begin{enumerate}
            \item $\vec{m} \le^*_P \vec{n} \iff P(\vec{m})\land [P(\vec{n})\to (\sigma(\vec{m})\le\sigma(\vec{n}))]$, and
            \item $\vec{m} <^*_P \vec{n} \iff P(\vec{m})\land [P(\vec{n})\to (\sigma(\vec{m})<\sigma(\vec{n}))]$
        \end{enumerate}
        are both in $\Gamma$.
    \end{enumerate}
\end{definition}

An example of a Spector class is the set of $\Pi^1_1$-definable sets with a parameter real $R$ over $\bbN$. It is well-known that every $\Pi^1_1$-set with parameter real $R$ is the least fixedpoint of an $R$-recursive monotone operator; If $\Gamma^\infty$ is the least fixedpoint of a monotone operator $\Gamma$, then we can `rank' $n\in \Gamma^\infty$ by looking the least $\alpha$ such that $n\in \Gamma^\alpha$. (Here we define $\Gamma^0 = \varnothing$ and $\Gamma^\alpha = \bigcup_{\xi<\alpha}\Gamma(\Gamma^\xi)$.) This rank gives the \emph{norm} of $\Gamma^\infty$. 

There are many Spector classes above a given structure $\frakM$. For example, both $\Pi^1_1$-sets and $\Sigma^1_2$-sets form Spector classes over $\bbN$. Each Spector class can be viewed as a universe above $\frakM$, and let us consider them collectively:
\begin{definition}[{\cite[Definition 3.25]{Pohlers2022Performance}}]
    For a Spector class $\Gamma$ and $R\in\Gamma$, define
    \begin{equation*}
        o(R) = \sup \{\sigma(\vec{x})+1\mid R(\vec{x})\land \text{$\sigma$ is a $\Gamma$-norm on $R$.}\}
    \end{equation*}
    For $\Lambda\subseteq \Gamma$, define $o(\Lambda) = \sup \{o(R)\mid R\in\Lambda\}$.
    We define the \emph{Spector spectrum} of $\frakM$ by the class of all Spector classes above $\frakM$, an the \emph{ordinal spectrum of $\frakM$} by $\Spec^\frakM = \{o(\Gamma)\mid \text{$\Gamma$ is a Spector class above }\frakM\}$.
    $\kappa_\frakM$ is the least element of $\Spec^\frakM$, or equivalently, $\kappa_\frakM = o(\SP_\frakM)$.
\end{definition}
Among Spector classes over $\frakM$, the most natural one we should consider is the \emph{least Spector class over $\frakM$}, which is the intersection of every Spector class over $\frakM$, which we will denote by $\SP_\frakM$. For an acceptable $\frakM$, we have $\delta^\frakM = o(\SP_\frakM)$ (cf. \cite[Theorem 3.34]{Pohlers2022Performance}), that is, the supremum of $\SP_\frakM$-norms is precisely the supremum of the ordertype of $\frakM$-definable well-orders. 
Hence we can use norms over a Spector class to describe the strength of an acceptable axiomatization of $\frakM$:
\begin{definition}[{\cite[Definition 3.29]{Pohlers2022Performance}}]
    Let $T$ be an acceptable axiomatization of a universe above $\frakM$.
    For a Spector class $\Gamma$ over $\frakM$, define
    \begin{equation*}
        o^\Gamma(T) = \sup\{\sigma_R(\vec{x})+1\mid R = \{\vec{x}\in |\frakM|^k\mid F(\vec{x})\}\in\Gamma \land T\vdash F(\vec{x})\}.
    \end{equation*}
    Then define $\Spec^\frakM(T) = \{o^\Gamma(T) \mid \Gamma\text{ is a Spector class over }\frakM\}$.
    We also define $\kappa^\frakM(T) = o^{\SP_\frakM}(T) = \min \Spec^\frakM(T)$.
\end{definition}
Again, for a strictly acceptable $\frakM$ and accpetable axiomatization $T$ of $\frakM$ proving Weak K\H onig's lemma, we have $\delta^\frakM(T) = \kappa^\frakM(T)$ \cite[Theorem 6.7.4]{Pohlers2009Prooftheory}.
We may think of a Spector class over $\frakM$ as a `measuring stick' gauging the strength of $T$ with different lengths and different scales. A larger Spector class can be thought of as a longer measuring stick with a coarser scale.
Although there is no clear limit of the Spector classes over $\frakM$, Pohlers focused only on the Spector class of the following form: For $\xi$ less than the least recursively inaccessible ordinal,
\begin{enumerate}
    \item $\SP^0_\frakM = \frakM$.
    \item $\SP^{\xi+1}_\frakM$ is the least Spector class over $(\frakM;\SP^\xi_\frakM)$.
    \item $\SP^\xi_\frakM = \bigcup_{\eta<\xi}\SP^\eta_\frakM$ if $\xi$ is a limit ordinal.%
    \footnote{For a limit $\xi$, $\SP^\xi$ is not a Spector class unless $\xi$ is recursively inaccessible.}
\end{enumerate}
We also define $\kappa^\frakM_0 = 0$, $\kappa^\frakM_{\mu+1} = o(\SP^{\mu+1}_\frakM)$, and $\kappa^\frakM_\lambda = \sup \{\kappa^\frakM_\xi\mid \xi<\lambda\}$ for a limit $\lambda$.
When we understand $\SP^\xi_\frakM$ as a structure, it is the expansion of $\frakM$ with new sets obtained from the iterated Spector operator as predicates.
Note that $\SP^\xi_\frakM$ is a \emph{first-order} structure, so for example, $\SP^\xi_\bbN$ is a structure of the first-order arithmetic with extra predicates, and not a model of second-order arithmetic.

\begin{theorem}[Pohlers {\cite[Remark 3.35]{Pohlers2022Performance}}]
    For an acceptable structure $\frakM$ and an ordinal $\mu\ge 1$ less than the least recursively inaccessible ordinal,
    \begin{equation*}
        o^{\frakM_{\mu+1}}(\Th(\frakM_\mu) + \text{Mathematical Induction}) = \varepsilon_{\kappa^{\frakM}_\mu+1}.
    \end{equation*}
    Hence if $\frakM$ is strictly acceptable, then we have
    \begin{equation*}
        \delta^{\frakM_\mu}(\Th(\frakM_\mu) + \text{Mathematical Induction}) = \varepsilon_{\kappa^{\frakM}_\mu+1}.
    \end{equation*}
\end{theorem}
In particular, if $\frakM=\bbN$, then $\bbN$ is strictly acceptable and $\kappa^\bbN_\mu = \omega^\CK_\mu$ for $\mu$ less than the least recursively inaccessible ordinal. (See \cite[Remark 3.33]{Pohlers2022Performance}.)
The following theorem shows the connection between Pohlers' characteristic ordinals and $\Pi^1_1[R]$-proof theory.
Its proof involves a precise definition of iterated hyperjumps, so we postpone its proof until \autoref{Section: Iterated Hyperjumps}.

\begin{proposition} \label{Proposition: Connecting Pohlers' ordinal and Pi 1 1 R PTO}
    Let $\xi$ be a successor ordinal less than the least recursively inaccessible ordinal.
    Suppose that $T$ is a strictly acceptable axiomatization of $\SP^{\xi}_\bbN$ containing $\Th(\SP^\xi_\bbN)$.
    Then we have 
    \begin{equation*}
        \delta^{\SP^{\xi}_\bbN}(T) = |\ACA_0 + \text{`$\HJ^\xi(\emptyset)$ exists'} +  T\restriction \HJ^\xi(\emptyset)|_{\Pi^1_1[\HJ^\xi(\emptyset)]},
    \end{equation*}
    where $T\restriction R$ is the theory obtained from $T$ by restricting its language to that of first-order arithmetic with an extra unary predicate for $R$.
\end{proposition}
We will also see later that the $\xi$th iterated hyperjump $\HJ^\xi(\emptyset)$ of the empty set for $\xi$ less than the least recursively inaccessible ordinal is a $\Pi^1_1$-singleton. 

\section{Quasidendroids} \label{Section: Quasidendroids}
Defining $\sfbeta$-logic requires a `functorial' proof tree for well-orders. Girard \cite{Girard1981Dilators} formulated functorial trees under the name \emph{quasidendroid}. However, quasidendroids correspond to dilators and not predilators, and we also need a pre- notion for quasidendroids.
Hence, we define the prequasidendroid first, and this makes our definition slightly different from that of Girard since Girard's definition does not fit in the pre-notion of quasidendroids. 

\begin{definition} \label{Definition: alpha prequasidendroid}
    Let $X$ and $\alpha$ be linear orders. An \emph{$\alpha$-prequasidendroid over $X$} is a subset $D$ of a tree over $X\sqcup\alpha$, the disjoint union of $X$ and $\alpha$, satisfying the following:
    \begin{enumerate}
        \item $D$ does \emph{not} contain an initial segment of an element of $D$.
    \end{enumerate}
    The next condition for the definition requires new \added{terminology}.
    We call elements of $X$ the \emph{term part of $D$}, and elements of $\alpha$ the \emph{parameter part of $D$}. 
    We use underlined elements $\underline{\xi}, \underline{\eta}$ to distinguish parameter parts from term parts. 
    
    Now let $D^*$ be the closure of $D$ under an initial segment.
    The next condition says the `successor' of a node in $D^*$ must be either all underlined or all non-underlined:
    \begin{enumerate}
        \setcounter{enumi}{1}
        \item For every $\sigma\in D^*$, it cannot be the case that there are $x\in X$ and $\xi\in\alpha$ such that $\sigma^\frown\lag x\rag\in D^*$ and $\sigma^\frown \lag\underline{\xi}\rag\in D^*$.
    \end{enumerate}
    
    We say $D$ is an \emph{$\alpha$-quasidendroid} if $D$ additionally satisfies the following:
    \begin{enumerate}
        \setcounter{enumi}{2}
        \item \label{Item: Def-alpha prequasidendroid cond 2} (Local well-foundedness) For each $\sigma\in D^*$, the set $\{x\in X \mid\exists \xi\in\alpha [\sigma^\frown\lag(x,\xi)\rag\in D^*]\}$ is well-ordered.
        \item There is no infinite branch over $D^*$.
    \end{enumerate}
\end{definition}

Girard had not defined the `pre-' notion for the quasidendroid. We define it since an `ill-founded prequasidendroid' has a more significant rule in this paper. Local well-foundedness in our definition of quasidendroid is important\added{;} Otherwise, the `linearized' quasidendroid can have an infinite decreasing sequence even when there is no infinite branch. Girard only considered the case when the term part $X$ is a set of ordinals. But we let $X$ be a general linear order, which should allow every predilator to be expressed by a prequasidendroid.

We can get a new $\alpha$-prequasidendroid from a $\beta$-prequasidendroid for $\beta>\alpha$ by \emph{mutilation}:

\begin{definition}
    Let $D$ be a $\beta$-prequasidendroid over $X$ and $f\colon\alpha\to\beta$ be an increasing map. Let us construct the \emph{mutilation} ${}^fD$ as follows: Starting from $D$,
    \begin{enumerate}
        \item Remove $\sigma\in D$ if $\sigma$ has a component of the form $\underline{\xi}$ for some $\xi\notin \ran f$.
        \item Replace every $\underline{f(\eta)}$ with $\underline{\eta}$.
    \end{enumerate}
    If we extend $f$ to a function from $X\sqcup \alpha$ to $X\sqcup \beta$ by letting $f(x)=x$ for $x\in X$, then we can rephrase the previous construction of ${}^D$ into the following equality:
    \begin{equation*}
        {}^fD = \{\lag \sigma_i\mid i<n\rag \mid \lag f(\sigma_i)\mid i<n\rag  \in D\},
    \end{equation*}
    We also define the \emph{mutilation function} $\sfm^D_f\colon ({}^fD)^*\to D^*$ by
    \begin{enumerate}
        \item $\sfm^D_f(\lag\rag)=\lag\rag$.
        \item $\sfm^D_f(\sigma^\frown \lag x \rag) = \sfm^D_f(\sigma)^\frown \lag x \rag$ when $x$ is a term part.
        \item $\sfm^D_f(\sigma^\frown \lag \underline{\xi}\rag) = \sfm^D_f(\sigma)^\frown \lag \underline{f(\xi)}\rag$ when $\underline{\xi}$ is a parameter part.
    \end{enumerate}
\end{definition}
It is easy to check that ${}^fD$ is an $\alpha$-prequasidendroid over $X$. 
Moreover, ${}^fD$ is an $\alpha$-quasidendroid if $D$ is $\beta$-quasidendroid, and $\sigma\in {}^f D$ iff $\sfm^D_f(\sigma)\in D$.

In the definition of the mutilation function, we see that $\sfm^D_f$ does not change the term part but changes the parameter part. We will often see modifying parameter parts of $\sigma\in D$ by an increasing function for a parameter, so let us employ the following notation: 
\begin{definition}
    Let $f\colon \alpha\to\beta$ be an increasing function between two linear orders $\alpha$ and $\beta$, and suppose that every element of $\alpha$ and $\beta$ are underlined. For a linear order $X$ whose elements are \emph{not} underlined, we define the function $\sfm_f\colon [X\sqcup \alpha]^{<\omega} \to [X\sqcup \beta]^{<\omega}$ by
    \begin{enumerate}
        \item $\sfm_f(\lag\rag)=\lag\rag$.
        \item $\sfm_f(\sigma^\frown \lag x \rag) = \sfm_f(\sigma)^\frown \lag x \rag$ when $x$ is not underlined.
        \item $\sfm_f(\sigma^\frown \lag \underline{\xi}\rag) = \sfm_f(\sigma)^\frown \lag \underline{f(\xi)}\rag$ when $\underline{\xi}$ is underlined.
    \end{enumerate}
\end{definition}

From an $\alpha$-prequasidendroid $D$ and $\beta>\alpha$, we can get a new $\beta$-prequasidendroid by \emph{instantiation}:
\begin{definition}
    For an $\alpha$-prequasidendroid $D$ over $X$ and $\beta\ge\alpha$, define its \emph{$\beta$-instantiation} $D_\beta$ as follows: Starting from $D$,
    \begin{enumerate}
        \item For each $\sigma\in D$, fix the set $\Occ(\sigma)$ of underlied elements occurring in $\sigma$. 
        \item For every increasing $f\colon \Occ(\sigma)\to \beta$, consider a copy of $\sigma$ in which every underlined element is replaced by its value by $f$.
        \item Collect every modified $\sigma$ for each $f\colon \Occ(\sigma)\to\beta$, and form the set $D_\beta$.
    \end{enumerate}
    
    Again, if we expand $f$ into a function $X\sqcup\alpha\to X\sqcup \beta$ by letting $f(x)=x$, then we can rephrase the definition of $D_\beta$ as follows:
    
    \begin{equation*}
        D_\beta = \{\lag f(\sigma_i)\mid i<n\rag \mid \lag \sigma_i\mid i<n\rag\in D\land f\colon \Occ(\lag \sigma_i\mid i<n\rag)\to\beta\text{ increasing}\}.
    \end{equation*}
\end{definition}
It is also easy to see that $D_\beta$ is a $\beta$-prequasidendroid. However, we do not know if $D_\beta$ is a $\beta$-quasidendroid when $D$ is an $\alpha$-quasidendroid since we do not know if $D_\beta$ has an infinite branch even when $D$ has no infinite branch. Note that if $D$ is locally well-founded, then so is $D_\beta$, so the presence of an infinite branch is the only reason why $D_\beta$ can fail to be a $\beta$-quasidendroid.

The instantiation should be the inverse operator of mutilation, and the following lemma gives a partial positive answer for the expectation:
\begin{lemma} \label{Lemma: Mutilation of an expansion Part 1}
    Let $D$ be an $\alpha$-prequasidendroid over $X$ and $\beta\ge\alpha$. For an increasing $f\colon \alpha\to\beta$, we have $D \subseteq {}^f (D_\beta)$.
\end{lemma}
\begin{proof}
    For $D \subseteq {}^f (D_\beta)$, let $\lag\sigma_i\mid i<n\rag\in D$. Then $\lag f(\sigma_i)\mid i<n\rag\in D_\beta$ by definition of $D_\beta$, so $\lag\sigma_i\mid i<n\rag\in {}^f (D_\beta)$.
\end{proof}

However, there is no reason in general that the reversed inclusion also holds. Even worse, there is no guarantee that the instantiation of a mutilation `restores' an $\alpha$-prequasidendroid. It requires an additional property called \emph{homogeneity}:
\begin{definition}
    An $\alpha$-prequasidendroid over $X$ is \emph{homogeneous} if for every increasing map $f,g\colon \gamma\to\alpha$, ${}^fD = {}^gD$.
\end{definition}
Note that the homogeneity of $D$ is equivalent to the following condition: For two increasing $f,g\colon\gamma\to\alpha$ and any finite sequence $\sigma$ over $X\sqcup \gamma$, we have $\sfm^D_f(\sigma)\in D$ iff $\sfm^D_g(\sigma)\in D$.
Then we can say that the mutilation of an instantiation restores a homogeneous prequasidendroid:
\begin{lemma} \label{Lemma: Mutilation of an expansion Part 2}
    Let $D$ be a homogeneous $\alpha$-prequasidendroid over $X$ and $f\colon\alpha\to\beta$. Then $D={}^f(D_\beta)$.
\end{lemma}
\begin{proof}
    We claim that ${}^f(D_\beta)\subseteq D$.
    Suppose that $\sigma\in {}^f(D_\beta)$, so $\sfm_f(\sigma)\in D_\beta$.
    By definition of $D_\beta$, there is $\tau\in D$ and an order isomorphism $h\colon \Occ(\tau)\to\Occ(\sfm_f(\sigma))$ such that $\sfm_h(\tau)=\sfm_f(\sigma)$.
    Observe that $|\Occ(\tau)| = |\Occ(\sfm_f(\sigma))| = |\Occ(\sigma)|$, and say it $m$. 
    Then fix an increasing enumeration $g\colon m \to \Occ(\tau)$ and $k\colon m\to \Occ(\sigma)$, and let $\sigma'=\sfm_k^{-1}(\sigma)$, $\tau'=\sfm_g^{-1}(\tau)$.
    Since $h\circ g$ and $f\circ k$ have the same domain and range that are finite, we have $h\circ g= f\circ k$. From $\sfm_{h\circ g}(\tau') = \sfm_{f\circ k}(\sigma')$, we get $\sigma'=\tau' \in {}^g D$. Homogeneity implies ${}^g D={}^k D$, so $\sigma'\in {}^k D$, implying $\sigma=\sfm_k(\sigma) \in D$.
\end{proof}

The next lemma says instantiating a mutilation of a homogeneous $\alpha$-prequasidendroid results in the previous $\alpha$-prequasidendroid when $\alpha=\omega$:
\begin{lemma} \label{Lemma: Expansion of a mutilation}
    If $f\colon \omega\to \beta$ is an increasing map and $D$ is a homogeneous $\beta$-predendroid, then $({}^fD)_\beta=D$.
\end{lemma}
\begin{proof}
    Let us prove $D\subseteq ({}^fD)_\beta$: Suppose that $\sigma\in D$. Take $m=|\Occ(\sigma)|$, and fix an increasing enumeration $h\colon m\to \Occ(\sigma)$.
    Then we have
    \begin{equation*}
        \sfm^{-1}_h(\sigma) \in {}^h D = {}^{f\restricts n} D \subseteq {}^f D.
    \end{equation*}
    
    For $({}^fD)_\beta \subseteq D$, suppose that $\sigma\in ({}^fD)_\beta$. 
    Thus we have $\tau \in {}^f D$ and an increasing map $g\colon \Occ(\tau)\to \beta$ satisfying $\sfm_g(\tau)=\sigma$. If $\imath \colon \Occ(\tau)\to \alpha$ is the inclusion map, then ${}^{f\circ \imath} D = {}^g D\ni \tau$. Hence $\sigma = \sfm_{f\circ \imath}(\tau) \in D$.
\end{proof}

Hence if we are given a homogeneous $\omega$-prequasidendroid, we can produce an $\alpha$-prequasidendroid for every linear order $\alpha$. It is similar to that for a predilator $F$ we can compute $F(\alpha)$ for every $\alpha$ from $F(\omega)$.
Thus let us impose a homogeneous $\omega$-prequasidendroid a special name:
\begin{definition}
    $D$ is a \emph{prequasidendroid} if $D$ is a homogeneous $\omega$-prequasidendroid. A prequasidendroid $D$ is a \emph{quasidendroid} if $D_\alpha$ is an $\alpha$-quasidendroid for every well-order $\alpha$.
\end{definition}
Note that Girard \cite{Girard1981Dilators} called our quasidendroid a \emph{strongly homogeneous quasidendroid}.

\begin{remark}
    It is convenient to view a prequasidendroid as a \emph{functorial} sequence of $\alpha$-prequasidendroid $\lag D_\alpha\mid \alpha\in \WO\rag$, where the functoriality means that ${}^f (D_\beta)=D_\alpha$ holds for every $\alpha$, $\beta$ and an increasing $f\colon \alpha\to\beta$. In one direction, the $\omega$-prequasidendroid $D_\omega$ in the functorial sequence is homogeneous, so it is a prequasidendroid.
    In the other direction, \autoref{Lemma: Expansion of a mutilation} implies we can recover $D_\alpha$ from $D_\omega$ for $\alpha\ge\omega$ by expansion, and the combination of \autoref{Lemma: Mutilation of an expansion Part 1} and \ref{Lemma: Mutilation of an expansion Part 2} implies we can obtain $D_n$ from $D_\omega$ for $n<\omega$.
    Understanding prequasidendroids as functorial sequences of $\alpha$-prequasidendroids is what prior materials for dilators and $\sfbeta$-logics take.
\end{remark}

Now let $D$ be a prequasidendroid over $X$. If $X$ is a subset of $\bbN$, then we can code $D$ as a subset of $\bbN$. This allows us to talk about $D$ being recursive:
\begin{definition}
    Let $D$ be a prequasidendroid over $\bbN$. We say $D$ is \emph{(primitive) recursive} if $D$ as a subset of $\bbN$ is (primitive) recursive.
\end{definition}

We can turn a prequasidendroid into a predilator. Girard used the \emph{linearization} $\Lin(D)$ given by $\Lin(D)(\alpha) := (D_\alpha,<_\KB)$. 
We will use a slightly different linearization construction, which is more useful in our context.
\begin{definition}
    Let $D$ be a prequasidendroid.
    Define the \emph{hull-linearization $\Lin^*(D)$ of $D$} by $\Lin^*(D)(\alpha) = (D_\alpha^*,\le_\KB)$.
    For $f\colon\alpha\to\beta$, we define $\Lin^*(D)(f) (\sigma) =\sfm^D_f(\sigma)$.
\end{definition}

Then we have the following:
\begin{proposition}
    Let $D$ be a prequasidendroid. Then $\Lin^*(D)$ is a predilator. Moreover, for a locally well-founded $D$ and a well-order $\alpha$, $D_\alpha^*$ has an infinite branch iff $\Lin^*(D)(\alpha)$ is ill-founded.
\end{proposition}
\begin{proof}
    We will define an appropriate support transformation to turn $\Lin^*(D)$ into a semidilator. For $\sigma\in D^*_\alpha$, define $\supp_\alpha(\sigma)=\Occ(\sigma)$.
    Then the direction computation gives $[f]^{<\omega}\circ \supp_\alpha = \supp_\beta\circ \Lin^*(D)(f)$.
    To check the support condition, fix $f\colon\alpha\to\beta$ and $\tau\in \Lin^*(D)(\beta)$ such that $\supp_\beta(\tau)\subseteq \ran f$.
    Then every underlined component of $\tau$ is a member of $f$, so $\sfm_f^{-1}(\tau)$ is well-defined. 
    Clearly $\tau = \sfm_f(\sfm_f^{-1}(\tau))\in \ran(\Lin^*(D)(f))$.

    To see $\Lin^*(D)$ is a predilator, let us check the monotonicity condition: That is, for two increasing $f,g\colon \alpha\to\beta$ such that $f\le g$ pointwise, we have $\Lin^*(D)(f)\le \Lin^*(D)(g)$ pointwise.
    However, one can see that for $\sigma\in \Lin^*(D)(\alpha)$ and $i<|\sigma|$, the $i$th component of $\sfm_f(\sigma)$ is less than or equal to $\sfm_g(\sigma)$, so $\sfm_f(\sigma) \le_\KB \sfm_g(\sigma)$.

    The last condition follows from the fact that if $T$ is a tree over a well-order, then $T$ has an infinite branch if and only if $(T,\le_\KB)$ is ill-founded.
\end{proof}

As a corollary, we have
\begin{corollary} \pushQED{\qed}
    Let $D$ be a prequasidendroid. Then $D$ is a quasidendroid iff $\Lin^*(D)$ is a dilator. Furthermore, if $D$ is locally well-founded and $\Lin^*(D)$ is not a dilator, then $\Clim(\Lin^*(D))$ is the least ordinal such that $D^*_\alpha$ has an infinite branch. \qedhere
\end{corollary}

\section{\texorpdfstring{$\sfbeta$}{Beta}-logic}
In this section, we review Girard's $\sfbeta$-logic, which is the right framework for $\Pi^1_2$-proof theory. This part is included not only because of the reader's convenience but also because arguments in this section will be reused to compute the `rank' of a $\Sigma^1_2$-singleton real.
Most results in this section are not new and appeared in Girard's \cite{Girard1982Logical2}. However, many parts of the definitions underwent non-trivial modifications to admit their formulations over a fragment of second-order arithmetic.

\subsection{Defining \texorpdfstring{$\sfbeta$}{Beta}-logic}
Girard defined the notion of $\sfbeta$-logic \cite{Girard1982Logical2} to analyze $\Pi^1_2$-consequences of a theory. To see how $\sfbeta$-logic works, let us recall what $\omega$-logic is: $\omega$-logic provides a satisfactory framework for $\Pi^1_1$-sentences. The semantics is based on $\omega$-models, and the syntax of $\omega$-logic uses infinitary deduction for quantifiers over natural numbers.
Likewise, $\sfbeta$-logic uses $\sfbeta$-models for its semantics. The $\sfbeta$-model of arithmetic is $\Pi^1_1$-correct $\omega$-models, and every $\sfbeta$-model admits its set-theoretic translation (See \cite[Ch. VI]{Simpson2009} for more details.) Hence we can regard $\sfbeta$-models as models with ordinals.
We want to liberate the notion of $\sfbeta$-model to \added{languages} other than the language of arithmetic, and axiomatic systems different from $\ACA_0$ or $\ATR_0$. Thus we define $\sfbeta$-models as models with well-orders, which is slightly different from, but will turn out to be equivalent to \added{$\sfbeta$-models} in second-order arithmetic.
Handling $\sfbeta$-logic over second-order arithmetic (Especially, over systems strictly weaker than $\ATR_0$) is inconvenient, but often unavoidable.
 
\begin{definition}
    A \emph{$\sfbeta$-language} is a language with a distinguished type $\Ord$ with a distinguished predicate $\le^\Ord$, and only terms of type $\Ord$ are variables.
    When we speak of a formula of a $\sfbeta$-language, we always assume that it has no free variables of type $\Ord$.

    A \emph{$\sfbeta$-theory} is a theory over a $\sfbeta$-language. For a $\sfbeta$-language $\calL$, a \emph{$\sfbeta$-model of $\calL$} is a model $M$ in a usual sense but $\Ord^M$ is a well-order, and $(\le^\Ord)^M$ is the linear order of the well-order $\Ord^M$. $\Ord^M=0$ is possible. 
\end{definition}
Girard's original definition for a $\sfbeta$-model requires $\Ord^M$ to be ordinal. Our definition is slightly more general, and this definition allows us \added{to formulate} $\sfbeta$-logic in second-order arithmetic.%
\footnote{Vauzeilles \cite{Vauzeilles1988CutElimOmegaLogic} called the modified version of $\sfbeta$-logic an \emph{$\Omega$-logic}, but we will call the modified version simply by $\sfbeta$-logic.}
Like other logics, we can define the semantic implication:
\begin{definition}
    A formula $\phi$ of a $\sfbeta$-theory $T$ is \emph{$\sfbeta$-valid in $T$} if every closed instantiation is valid over a $\sfbeta$-model $M$ of $T$. We write it by $T\vDash^\sfbeta \phi$.
\end{definition}

The following example shows the connection between Girard's $\sfbeta$-models and $\sfbeta$-models in second-order arithmetic:
\begin{example}[{\added{Girard \cite[Example 10.1.3.]{Girard1982Logical2}}}]
    Consider the language $\calL$ extending the language of second-order arithmetic by adding a new type $\Ord$, and new binary relation symbol $\le^\Ord$ and ternary relation symbol $\Emb$ of type $\bbR\times\bbN\times\Ord$.
    The theory $\RCA_0^\sfbeta$ is an extension of $\RCA_0$ with the following axioms:
    \begin{enumerate}
        \item $\forall^1 X\forall^0 n\forall^\Ord \xi,\eta [\Emb(X,n,\xi)\land \Emb(X,n,\eta) \to \xi=\eta]$.
        \item $\forall^1 X\forall^0 n \exists^\Ord\xi [\Emb(X,n,\xi)]$.
        \item $\forall^1 X \forall^0 n,m\forall^\Ord\xi,\eta [\WO(X)\land n<_X m \land \Emb(X,n,\xi)\land \Emb(X,m,\eta)]\to \eta\nleq^\Ord \xi$.
    \end{enumerate}
    These axioms describe `$\Emb$ is a predicate describing a way to embed a well-order into $\Ord$.'
    Then, any $\sfbeta$-model of second-order arithmetic induces a Girard-styled $\sfbeta$-model of $\RCA_0^\sfbeta$: Suppose that $M$ is a $\Pi^1_1$-correct $\omega$-model of arithmetic.  Then define
    \begin{enumerate}
        \item $\Ord^M = \{\operatorname{otp}(\alpha)\mid \alpha\in M\land  M\vDash \WO(\alpha) \}$.
        \item $\Emb^M(X,n,\xi)$ iff either ($\lnot\WO(X)\land \xi=0$) or ($\WO(X)$ and $\xi$ is isomorphic to $\{i\in \field(X)\mid i<_X n\}$).
    \end{enumerate}
    It gives a Girard-styled $\sfbeta$-model for $\RCA^\sfbeta_0$. Note that the previous construction works over a sufficiently strong set theory (more precisely, primitive recursive set theory plus Axiom Beta). If we stick to second-order arithmetic, we may take $\Ord^M$ by the ordered sum of all well-orders coded in $M$. Even in that case, we need $\ATR_0$ as a metatheory to define $\Emb^M$ correctly.
    
    Conversely, if $N$ is a Girard-styled $\sfbeta$-model for $\RCA^\sfbeta_0$, then it is a $\Pi^1_1$-correct $\omega$-model: $\RCA_0$ proves $\bbN$ with the usual order is a well-order, so $\bbN^M$ is isomorphic to a well-order. Since $\RCA_0$ also proves $\bbN$ has no limit point, the same also holds for $\bbN^M$, which implies $\bbN^M\cong M$.
    It is clear that for a recursive linear order $\alpha$ if $M\vDash\WO(\alpha)$, then $\alpha$ is a well-order. Thus by the Kleene normal form theorem, $M$ is $\Pi^1_1$-correct.
\end{example}

So far, we have defined the semantics for a $\sfbeta$-logic. Then what is the syntax of $\sfbeta$-logic? 
\modified{In an earlier version of the draft, we followed Girard's sequent-calculus definition \cite[Definition 10.1.6]{Girard1982Logical2}. In this version, we will employ Tait calculus. We introduce the syntax of $\sfbeta$-logic for a fixed well-order $\alpha$, then introduce how they are functorially combined.}
\begin{definition}
    Let $\calL$ be a $\sfbeta$-language and $\alpha$ be a well-order. The language $\calL[\alpha]$ is obtained from $\calL$ by adding constants $\overline{\xi}$ of type $\Ord$ for $\xi\in\alpha$. The new language is not technically a $\sfbeta$-language since it has constant symbols of type $\Ord$, and we do not allow free variables of type $\Ord$. \modified{Since we will work with Tait Calculus, we introduce negated atomic formulas in $\calL$ as new atomic symbols.}

    \modified{The Tait calculus} $\LK_\alpha$ is defined as follows: The sequents are made of formulas of $\calL[\alpha]$, and we have the following axioms and inference rules, where $\tau$ is a type other than $\Ord$.
    \newline
    \begin{center}
        \AxiomC{}
        \RightLabel{Id}
        \UnaryInfC{$\Gamma, A, \lnot A$}
        \DisplayProof
        \\[2em]
        \AxiomC{}
        \RightLabel{Ax (When $\zeta\le\xi\in\alpha$)}
        \UnaryInfC{$\Gamma, \overline{\zeta} \le^\Ord \overline{\xi}$}
        \DisplayProof
        \hspace{2em}
        \def\fCenter{\mbox{ $\vdash$ }}
        \AxiomC{}
        \RightLabel{Ax (When $\xi<\zeta\in\alpha$)}
        \UnaryInfC{$\Gamma, \overline{\zeta} \nleq^\Ord \overline{\xi}$}
        \DisplayProof
        \\[1em]
        \AxiomC{$\Gamma,A\land B, A$}
        \AxiomC{$\Gamma,A\land B, B$}
        \RightLabel{$\land$}
        \BinaryInfC{$\Gamma,A\land B$}
        \DisplayProof
        \hspace{2em}
        \def\fCenter{\mbox{ $\vdash$ }}
        \AxiomC{$\Gamma,A\lor B,A$}
        \RightLabel{$\lor$}
        \UnaryInfC{$\Gamma, A\lor B$}
        \DisplayProof
        \hspace{2em}
        \def\fCenter{\mbox{ $\vdash$ }}
        \AxiomC{$\Gamma,A\lor B,B$}
        \RightLabel{$\lor$}
        \UnaryInfC{$\Gamma, A\lor B$}
        \DisplayProof
        \\[2em]
        \AxiomC{$\Gamma, \forall^\tau x A(x),A(v)$}
        \RightLabel{$\forall^\tau$}
        \UnaryInfC{$\Gamma, \forall^\tau x A(x)$}
        \DisplayProof
        \hspace{2em}
        \AxiomC{$\Gamma, \exists^\tau x A(x), A(t)$}
        \RightLabel{$\exists^\tau$}
        \UnaryInfC{$\Gamma, \exists^\tau x A(x)$}
        \DisplayProof
        \\[2em]
        \AxiomC{$\cdots \Gamma, \forall^\Ord x A(x), A(\overline{\zeta}) \cdots$ (every $\zeta\in\alpha$)}
        \RightLabel{$\forall^\Ord$}
        \UnaryInfC{$\Gamma, \forall^\Ord x A(x)$}
        \DisplayProof
        \hspace{2em}
        \AxiomC{$\Gamma,\exists^\Ord x ,A(\overline{\zeta})$}
        \RightLabel{$\exists^\Ord$}
        \UnaryInfC{$\Gamma, \exists^\Ord x A(x)$}
        \DisplayProof
        \\[2em]
        \AxiomC{$\Gamma, A$}
        \AxiomC{$\Gamma,\lnot A$}
        \RightLabel{Cut}
        \BinaryInfC{$\Gamma$}
        \DisplayProof
        \vspace{1em}
    \end{center}
    \added{
    We view sequents as sets of formulas, so the order of formulas in a sequent is ignored.
    In $\forall^\tau$, the \emph{eigenvariable} $v$ which applies the rule must not appear free in $\Gamma$.}
    We call a proof tree generated by the previous rules an \emph{$\alpha$-preproof}. If an $\alpha$-preproof  has no infinite branch, then we call the $\alpha$-preproof an \emph{$\alpha$-proof}.
\end{definition}
We have only defined $\sfbeta$-proof system $\LK_\alpha$ for predicate logic, but we will see a modification of $\LK$ by adding additional rules (\added{like} $\omega$-rules or rules for an inductive definition for a fixed arithmetical positive operator) into $\LK_\alpha$.  
However, the previous notion of $\LK_\alpha$ and $\alpha$-proof have a fault in that both $\LK_\alpha$ and an $\alpha$-proof are not \added{syntactic objects}. For example, we expect that an $\alpha$-proof should have a `recursive' structure, but the current definition of an $\alpha$-preproof is a mixture of arbitrary well-orders $\alpha$ and syntactic objects. Also, we expect $\LK_\alpha$ be defined uniformly on $\alpha$, but we have not stated what the uniformity is. 
The latter issue is not a big problem if we restrict ourselves to $\LK_\alpha$ or its trivial modification (like $T+\LK_\alpha$, obtained from $\LK_\alpha$ by adding statements of $T$ as axioms), but \added{it} becomes problematic if we consider the expansion of $\LK_\alpha$ by adding new inference rules.
To resolve these issues, we provide a formal definition of a \emph{$\sfbeta$-proof system}, which generalizes the system $\lag \LK_\alpha\mid \alpha\in \WO\rag$ and other formal systems we will introduce.

\begin{definition}
    Let $\calL$ be a $\sfbeta$-language and $\alpha$ \added{be} a linear order. An \emph{$\alpha$-proof system} $\frakT$ consists of the \emph{set of inference rules} $\mathsf{Rule}^\frakT = \mathsf{Rule}$, the \emph{base linear order} $\bfB^\frakT = \bfB$, and an index set $I_\sfR$ for each $\sfR\in \mathsf{Rule}$.
    Each inference rule takes the form
    \begin{equation} \label{Formula: General form of an inference rule}
    \begin{mathprooftree}
        \def\fCenter{\mbox{ $\vdash$ }}
        \AxiomC{$\lag \Gamma_\iota \mid \iota\in I_\sfR\rag$}
        \RightLabel{$\sfR$}
        \UnaryInfC{$\Gamma$}
    \end{mathprooftree}
    \end{equation}
    
    We also require that $I_\sfR$ is either a (possibly empty) constant subset of $\bfB$ \added{that only depends on the choice of $\sfR$} or equal to $\alpha$.
    Typically $\bfB=\bbN$ and $I_\sfR$ is an initial segment of either $\bfB$ or equal to $\alpha$. We call $I_\sfR$ the \emph{arity of $\sfR$}, and we say \emph{$\sfR$ has an ordinal arity} when $I_\sfR=\alpha$.
    
    A \emph{$\frakT$-preproof} is a pair $(D,\phi)$ of a tree $D\subseteq [\alpha\cup \bfB]^{<\omega}$ with a labeling $\phi$ satisfying the following:
    For each $\sigma\in D$, $\phi(\sigma)=(\sfR,\Gamma)$ for some $\sfR\in\mathsf{Rule}$ and $\Gamma$ is a sequent such that
    \begin{equation*}
        I_\sfR = \{\iota\mid \sigma^\frown\lag\iota\rag\in D\}
    \end{equation*}
    and if $\phi(\sigma^\frown\lag\iota\rag) = (\sfR_\iota, \Gamma_\iota)$ then $\lag\Gamma_\iota\mid \iota\in I_\sfR\rag$ and $\Gamma$ satisfy \eqref{Formula: General form of an inference rule}.
\end{definition}

When we define an $\alpha$-proof system, we do not build just a single $\alpha$-proof system but a family of $\alpha$-proof systems for each well-order $\alpha$. The family of proof systems shares the same $\sfbeta$-language, and is somehow functorial. We will define the precise meaning of being functorial, and let us start with a preliminary notion:

\begin{definition}
    Let $\calL$ be a $\sfbeta$-language. For two well-orders $\alpha$, $\gamma$ and an increasing map $f\colon \alpha\to\gamma$, let us define $\sfm_f(\psi)$ for an $\calL[\alpha]$-formula $\psi$ by replacing every occurance of $\overline{\xi}$ in $\psi$ for $\xi\in \alpha$ with $\overline{f(\xi)}$.
\end{definition}
We extend the notion $\sfm_f$ to finite sequences of formulas. For a sequence $\Gamma = A_0,\cdots,A_{m-1}$ of formulas, $\sfm_f(\Gamma) = \sfm_f(A_0),\cdots,\sfm_f(A_{m-1})$.

\begin{definition}
    Let $\frakT$ be an $\omega$-proof system (not to be confused with $\omega$-logic.) We say $\frakT$ is a \emph{$\sfbeta$-proof system} if $\bfB$ is well-founded and moreover the following holds: 
    \begin{enumerate}
        \item If $I_\sfR=\omega$, then the premises of $\sfR$ must have the form $\Gamma'[\overline{\iota}/x]$ for $\iota\in I_\sfR=\omega$, where $\Gamma'$ is a sequent over the language $\calL[\omega]$ with an ordinal variable $x$.

        \item If $I_\sfR = \omega$, then $\Occ(\Gamma_\iota)\subseteq \Occ(\Gamma)\cup\{\underline{\iota}\}$, where the sequents are given as \eqref{Formula: General form of an inference rule} \added{and $\Occ(\Gamma)$ is the set of ordinals occurring in a formula in $\Gamma$.}
        
        \item If $I_\sfR\subseteq \bfB$, then $\Occ(\Gamma_\iota)\subseteq \Occ(\Gamma)$ for each $\iota\in I_\sfR$, where the sequents are given as \eqref{Formula: General form of an inference rule}.
    \end{enumerate}
\end{definition}

\begin{remark}
    We will see later the case when $\bbN\subseteq\bfB$ and $\alpha=\omega$.
    In this case, denoting both natural numbers as members of $\bfB$ and as members of the ordinal parameter set $\omega$ by $n$ will be confusing.
    To avoid confusion, we always \added{underline} the parameter elements when \added{they appear} in the \added{formal definition of} $\sfbeta$-logic\added{, like, when we label sequents by members of $\alpha$. A comparable notation is the underlined parameters in quasidendroids. Overlined ordinals will still denote the corresponding constant symbols in the expanded language occurring in formulas.}
\end{remark}

The next lemma easily follows from the induction on $\sigma$, which roughly states \added{that} in a $\sfbeta$-proof system, the parameter ordinals of a sequent in a certain position of a preproof only \added{depend} on the parameters in the initial sequent and the position itself.
\begin{lemma} \pushQED{\qed} \label{Lemma: Formula: Beta proof system parameter restriction}
    Let $\frakT$ be a $\sfbeta$-proof system and $\pi=(D,\phi)$ be a $\frakT$-preproof. For $\sigma\in D$ with $\phi(\sigma) = (\sfR_\sigma,\Gamma_\sigma)$, we have
    \begin{equation*}
        \Occ(\Gamma_\sigma) \subseteq \Occ(\Gamma_{\lag\rag}) \cup \Occ(\sigma) \qedhere 
    \end{equation*}
\end{lemma}

We can instantiate every $\sfbeta$-system to an $\alpha$-proof system uniformly:
\begin{definition}
    Let $\frakT$ be a $\sfbeta$-proof system and $\alpha$ a well-order. Let us define its \emph{instantiation} $\frakT_\alpha$ as follows: $\frakT_\alpha$ and $\frakT$ have the same set of inference rules and base linear order, but $I^{\frakT_\alpha}_\sfR=\alpha$ when $I^\frakT_\sfR = \omega$.
    The rule $\sfR$ is well-defined in this case since the definition of a $\sfbeta$-proof system requires the hypotheses of $\sfR$ \added{to} take the same form except for parameters $\iota\in I_\sfR$.
\end{definition}

We also require a homogeneity to proofs of a $\sfbeta$-proof system:
\begin{definition}
    Let $\frakT$ be an $\alpha$-proof system. For a $\frakT$-preproof $\pi=(D,\phi)$ and an increasing $f\colon \gamma\to\alpha$ if we define the \emph{mutilation} ${}^f\pi=({}^fD,{}^f\phi)$ by
    \begin{enumerate}
        \item ${}^fD = \{\sigma\mid \sfm_f(\sigma)\in D\}$,
        \item If $\sigma\in {}^fD$, then $f(\sigma)\in D$ and every overlined element (i.e., ordinal) appearing in $\phi(f(\sigma))$ is in $\ran f$. Then define
        \begin{equation*}
            \bigl({}^f \phi\bigr)(\sigma) = \sfm_f^{-1}\bigl(\phi\bigl(\sfm_f^D(\sigma)\bigr)\bigr).
        \end{equation*}
    \end{enumerate}
    In a more plain term, we construct ${}^f\pi$ as follows:
    \begin{enumerate}
        \item In a $\gamma$-preproof $\pi$, remove all premises containing overlined ordinals not in $\ran f$.
        \item Replace every $\overline{f(\xi)}$ with $\overline{\xi}$.
    \end{enumerate}
    
    We say a $\frakT$-preproof $\pi$ is \emph{homogeneous} if for every increasing $f,g\colon \gamma\to\alpha$, we have ${}^f\pi = {}^g\pi$ as $\frakT_\gamma$-preproofs.
\end{definition}

Note that the homogeneity of a preproof has the following finitary characterization:
\begin{lemma}
    Let $\frakT$ be a $\sfbeta$-proof system and $\alpha$ a well-order.
    For a $\frakT_\alpha$-preproof $\pi$, $\pi$ is homogeneous if and only if for every $n\in\bbN$ and increasing $f,g\colon n\to \alpha$, we have ${}^f\pi = {}^g\pi$.
\end{lemma}
\begin{proof}
    One direction is clear. For the remaining direction, suppose that $f,g\colon\gamma\to\alpha$ are two increasing functions.
    We claim that for $\sigma\in [\alpha\cup\bfB]^{<\omega}$,
    \begin{equation*}
        \sfm_f(\sigma)\in D \iff \sfm_g(\sigma)\in D.
    \end{equation*}
    Let $h\colon n\to \Occ(\sigma)$ be the unique increasing bijection. Then we have $\sigma\in \ran \sfm_h$ and
    \begin{equation*}
        \sfm_f(\sigma) = \sfm_{f\circ h}(\sfm_h^{-1}(\sigma))\in D \iff 
        \sfm_g(\sigma) = \sfm_{g\circ h}(\sfm_h^{-1}(\sigma))\in D.
    \end{equation*}
    It gives ${}^fD={}^gD$. Proving ${}^f\phi = {}^g\phi$ is similar, so we omit its proof.
\end{proof}

For a homogeneous $\frakT$-preproof $\pi$ for a sequent \modified{$\Gamma^*$} over $\calL$ (so $\Occ(\Gamma^*)$ is empty), we can instantiate it and get a $\frakT_\alpha$-preproof:
\begin{proposition}
    Let $\frakT$ be a $\sfbeta$-proof system and $\alpha$ a well-order.
    Suppose that \added{$\pi=(D,\phi)$} is a homogeneous $\frakT$-preproof with the conclusion \modified{$\Gamma^*$} such that $\Occ(\Gamma^*)=\varnothing$.
    Let us define the instantiation $\pi_\alpha=(D_\alpha,\phi_\alpha)$ to $\alpha$ as follows:
    \begin{enumerate}
        \item $D_\alpha$ is the set of all $\sfm_f(\sigma)$ for each $\sigma\in D$ and an increasing $f\colon \Occ(\sigma)\to \alpha$.
        \item For the same $\sigma\in D$ and $f$, we define
        \begin{equation*}
            \phi_\alpha(\sfm_f(\sigma)) = \sfm_f(\phi(\sigma)).
        \end{equation*}
    \end{enumerate}
    Then $\pi_\alpha$ is a well-defined $\frakT_\alpha$-preproof and is homogeneous.
\end{proposition}
\begin{proof}
    First, we need to ensure that $\phi_\alpha$ is well-defined.
    Note that \added{for $\sigma\in D$ and $f\colon \Occ(\sigma)\to \alpha$, $\sfm_f(\phi(\sigma))$} is well-defined by \autoref{Lemma: Formula: Beta proof system parameter restriction}: The lemma guarantees every ordinal parameter \added{occurring in the} second component of $\phi(\sigma)$ belongs to the domain of $f$.
    \added{We claim that the definition of $\phi_\alpha$ does not depend on the choice of $f$ and $\sigma$ in the following sense: If}
    $\sigma,\tau\in D$, $f\colon \Occ(\sigma)\to\alpha$, and $g\colon \Occ(\tau)\to \alpha$ satisfy $\sfm_f(\sigma) = \sfm_g(\tau)$ then $\sfm_f(\phi(\sigma)) = \sfm_g(\phi(\tau))$.
    
    Let $e\colon \lvert\ran f\cup \ran g\rvert\to \ran f\cup \ran g$ be the enumeration function and $h\colon \lvert\ran f\rvert\to \lvert\ran f\cup \ran g\rvert$, $k\colon \lvert\ran g\rvert\to \lvert\ran f\cup \ran g\rvert$ be such that $f=e\circ h$ and $g=e\circ k$.
    Then $\sfm_f(\sigma)=\sfm_g(\tau)$ implies $\sfm_h(\sigma)=\sfm_k(\tau)$.
    Also by homogeneity, we have ${}^{\Id}\phi = {}^h \phi = {}^k\phi$, so
    \begin{equation*}
        \sfm_h(\phi(\sigma)) = \sfm_h\bigl(\sfm_h^{-1}\bigl(\phi(\sfm_h(\sigma)\bigr)\bigr) = \phi(\sfm_h(\sigma)) = \phi(\sfm_k(\tau)) = \sfm_k\bigl(\sfm_k^{-1}\bigl(\phi(\sfm_k(\tau)\bigr)\bigr) = \sfm_k(\phi(\tau)),
    \end{equation*}
    which implies $\sfm_f(\phi(\sigma)) = \sfm_g(\phi(\tau))$.

    Then let us prove that $\pi_\alpha$ is a $\frakT_\alpha$-preproof.
    Suppose that we are given \added{$\sigma\in D$, $f\colon\Occ(\sigma)\to \alpha$,} and $\phi(\sigma)=(\sfR,\Gamma)$.
    Then we have $\phi(\sfm_f(\sigma)) = (\sfR,\sfm_f(\Gamma\vdash\Delta))$.
    If $I_\sfR\subseteq \bfB$, then $\lag\sfm_f(\Gamma_\iota)\mid \iota\in I_\sfR\rag$ and $\sfm_f(\Gamma)$ clearly satisfy \eqref{Formula: General form of an inference rule}.
    Otherwise, $I_\sfR$ over $\frakT_\alpha$ changes, so we need additional steps to verify \eqref{Formula: General form of an inference rule}.
    In this case, $\Gamma_\iota$ takes the form $\Gamma'[\overline{\iota}/x]$, and $\sigma^\frown\lag\underline{\iota}\rag\in D$ for each $\iota\in\omega$.
    We claim that for every $\xi\in\alpha$, $\sfm_f(\sigma)^\frown\lag \underline{\xi}\rag\in D_\alpha$ and
    \begin{equation*}
        \phi_\alpha(\sfm_f(\sigma)^\frown\lag \underline{\xi}\rag) = (\sfR',\sfm_f(\Gamma')[\overline{\xi}/x]).
    \end{equation*}
    for some $\sfR'$.
    If $\xi\in \ran f$ so there is $k\in\omega$ such that $f(k)=\xi$, then $\sfm_f(\sigma^\frown \lag \underline{k}\rag) = \sfm_f(\sigma)^\frown\lag \underline{\xi}\rag \in D_\alpha$ and 
    \begin{equation*}
        \phi_\alpha(\sfm_f(\sigma^\frown \lag \underline{\xi}\rag)) = 
        \sfm_f(\phi_\alpha(\sigma^\frown \lag \underline{k}\rag)) = 
        \sfm_f(\sfR', \Gamma'[\overline{k}/x]\vdash \Delta'[\overline{k}/x]) = 
        (\sfR', \sfm_f(\Gamma')[\overline{\xi}/x]).
    \end{equation*}
    \added{If $\xi\notin \ran f$, let $g\colon \lvert\Occ(\sigma)\rvert+1\to \ran f\cup \{\xi\}$ be the increasing enumeration.}
    Then by homogeneity, ${}^{g^{-1}\circ f} D = {}^{\Id_{\Occ(\sigma)}}D\ni \sigma$, so $\sfm_{g^{-1}\circ f}(\sigma)\in D$.
    Moreover, by the homogeneity of $\pi$, we also have $\phi(\sigma) = ({}^{g^{-1}\circ f}\phi)(\sigma)$. It implies
    \begin{equation*}
        \phi(m_{g^{-1}\circ f}(\sigma)) = \sfm_{g^{-1}\circ f}(\phi(\sigma)) = \bigl(\sfR, \sfm_{g^{-1}\circ f}(\Gamma)\bigr).
    \end{equation*}
    In particular, we have $\sfm_{g^{-1}\circ f}(\sigma)^\frown \lag \underline{k}\rag \in D$ for every $k\in\omega$.
    Hence if $k$ be the natural number such that $g(k)=\xi$, then
    \begin{equation*}
        \sfm_f(\sigma)^\frown\lag\underline{\xi}\rag = \sfm_g\bigl(\sfm_{g^{-1}\circ f}(\sigma)^\frown \lag \underline{k}\rag\bigr) \in D_\alpha
    \end{equation*}
    Moreover,
    \begin{equation*}
        \phi_\alpha(\sfm_f(\sigma)^\frown\lag\underline{\xi}\rag) = \phi_\alpha(\sfm_g(\sfm_{g^{-1}\circ f}(\sigma)^\frown \lag \underline{k}\rag))
        = \sfm_g(\phi(\sfm_{g^{-1}\circ f}(\sigma)^\frown \lag \underline{k}\rag).
    \end{equation*}
    By homogeneity again, we have
    \begin{itemize}
        \item $\phi\bigl(\sfm_{g^{-1}\circ f}(\sigma)\bigr) = \sfm_{g^{-1}\circ f}\bigl(\phi(\sigma)\bigr) = \bigl(\sfR,\sfm_{g^{-1}\circ f}(\Gamma)\bigr)$, and 
        \item $\phi\bigl(\sfm_{g^{-1}\circ f}(\sigma^\frown\lag \underline{j}\rag)\bigr) = \sfm_{g^{-1}\circ f}\bigl(\phi(\sigma^\frown\lag \underline{j}\rag)\bigr) = \bigl(\sfR',\sfm_{g^{-1}\circ f}(\Gamma')\bigl[\underline{g^{-1}(f(j))}/x\bigr]\bigr)$ for $j\in \dom f$ and some $\sfR'$,
    \end{itemize}
    meaning that $\phi\bigl(\sfm_{g^{-1}\circ f}(\sigma)^\frown \lag \underline{k}\rag\bigr)=\bigl(\sfR'',\sfm_{g^{-1}\circ f}(\Gamma')[\underline{k}/x]\bigr)$ for some $\sfR''$. Hence we have
    \begin{equation*}
        \phi_\alpha\bigl(\sfm_f(\sigma)^\frown\lag\underline{\xi}\rag\bigr) = 
        \sfm_g\bigl(\phi\bigl(\sfm_{g^{-1}\circ f}(\sigma)^\frown \lag \underline{k}\rag\bigr)\bigr)
        = \bigl(\sfR'',\sfm_{f}(\Gamma')[\underline{\xi}/x]\bigr)
    \end{equation*}
    finishing the proof of that $\pi_\alpha$ is a $\frakT_\alpha$-preproof. 
    For homogeneity, observe that for $f\colon n\to \alpha$, we have ${}^f(D_\alpha) = \{\sigma\in D\mid \Occ(\sigma)\subseteq n\}$, which does not depend on $f$. Moreover, $\bigl({}^f(\phi_\alpha)\bigr)(\sigma) = \phi(\sigma)$, finishing the proof.
\end{proof}
Now let us define the provability relation for $\sfbeta$-logic:
\begin{definition}
    Let $\frakT$ be a $\sfbeta$-proof system.
    We say an $\frakT$-preproof $\pi = (D,\phi)$ is a \emph{$\sfbeta$-preproof from $\frakT$} if it is homogeneous.
    Furthermore, $\pi$ is a \emph{$\sfbeta$-proof} if, in addition, the instantiation $\pi_\alpha$ is an $\alpha$-proof (so $\pi_\alpha$ has no infinite branch) for every infinite well-order $\alpha$. We say an $\omega$-preproof $\pi$ is \emph{(primitive) recursive} if there is a (primitive) recursive code for $\pi$.
    
    We say a sequent $\Gamma$ over a $\sfbeta$-language $\calL$ is \emph{$\sfbeta$-provable} if there is an $\sfbeta$-proof from $\LK$ for the sequent $\Gamma$.
    For an $\calL$-formula $\phi$ and an $\calL$-theory $T$, we say \emph{$\phi$ is $\sfbeta$-provable from $T$} if there is a $\sfbeta$-proof for $\phi$ from the system $\LK+T$, whose definition is the same with $\LK$ but has extra sequents $\vdash A$ for each $A\in T$ as axioms.
    We write $T\vdash^\sfbeta \phi$ if \modified{the sequent consisting of the single sentence} $\phi$ is $\sfbeta$-provable from $T$.
\end{definition}

Note that a $\sfbeta$-preproof itself does not have a quasipredendroid structure, but we can extract a quasipredendroid from it:
\begin{definition}
    Let $\frakT$ be a $\sfbeta$-proof system and $\pi=(D,\phi)$ be a $\sfbeta$-preproof. Let us define the associated prequasidendroid $\Den(\pi)$ as follows:
    It is a prequasidendroid over the linear order $\bfB\cup\{\infty\}$, where $\infty$ is a new element larger than any elements of $\bfB$. We define $\Den(\pi)$ by
    \begin{equation*}
        \Den(\pi) = \{\lag\sigma(0),\infty,\sigma(1),\infty,\cdots, \sigma(m-1),\infty\rag\mid \sigma\in D^*\land |\sigma|=m\}.
    \end{equation*}
\end{definition}

A simpler definition should be
\begin{equation*}
    \Den(\pi) = \{\sigma^\frown\lag\infty\rag \mid \sigma\in D^*\},
\end{equation*}
but it does not work since if $\phi(\sigma)=(\sfR,\Gamma)$ for some $\sigma\in D^*$ and an ordinal arity $\sfR$, then every immediate successor of $\sigma$ is of the form $\sigma^\frown\lag \underline{m}\rag$, so adding $\sigma^\frown\lag\infty\rag$ could violate the definition of a prequasidendroid.

\begin{proposition}
    Let $\frakT$ be a $\sfbeta$-proof system and $\pi=(D,\phi)$ be a $\sfbeta$-preproof. Then $\Den(\pi)$ is a prequasidendroid. Also, for each linear order $\alpha$,
    \begin{equation} \label{Formula: Proof quasidendroid instantiation}
        \Den(\pi)_\alpha = \{\lag\sigma(0),\infty,\sigma(1),\infty,\cdots, \sigma(m-1),\infty\rag\mid \sigma\in D^*_\alpha\land |\sigma|=m\}.
    \end{equation}
    
    Moreover, $\pi$ is a $\sfbeta$-proof if and only if $\Den(\pi)$ is a quasidendroid.
\end{proposition}
\begin{proof}
    $\Den(\pi)$ being a prequasidendroid is easy to check. For \eqref{Formula: Proof quasidendroid instantiation}, observe that we defined $D_\alpha$ by the set of all $\sfm_f(\sigma)$ for $\sigma\in D$ and an increasing $f\colon\Occ(\sigma)\to\alpha$, and
    \begin{equation*}
        \Occ(\lag\sigma(0),\infty,\cdots,\sigma(m-1),\infty\rag) = \Occ(\sigma).
    \end{equation*}
    
    For the moreover part, observe that $\pi_\alpha$ has an infinite branch if and only if $D_\alpha$ has an infinite branch, which is equivalent to that $\Den(\pi)_\alpha$ has an infinite branch.
\end{proof}

\subsection{Completeness of \texorpdfstring{$\sfbeta$}{Beta}-logic} \label{Subsection: completeness of beta logic}
We have defined the semantics and the syntax of $\sfbeta$-logic. To ensure they cohere, we need to prove the soundness and completeness of $\sfbeta$-logic,

\begin{proposition}[Soundness]
    Let $\phi$ be a sentence over $\calL$ and $T$ be a $\calL$-theory for a $\sfbeta$-language $\calL$.
    If there is an $\alpha$-proof for $T\vdash \phi$, then every $\sfbeta$-model $M$ of $T$ such that $\Ord^M=\alpha$ satisfies $\phi$. As a corollary, if $T\vdash^\sfbeta \phi$, then $T\vDash^\sfbeta \phi$. 
\end{proposition}
\begin{proof}
    We can prove that if there is an $\alpha$-proof for a sequent $\Gamma$, and if $M$ is a $\sfbeta$-model with $\Ord^M=\alpha$, then $M\vDash \bigvee\Gamma$, which follows from an induction on the proof tree.
    If $T\vdash^\sfbeta \phi$, then we can find a finite $T_0\subseteq T$ and an $\alpha$-proof for $T_0\vdash\phi$ for each ordinal $\alpha$. Thus for every $\sfbeta$-model $M\vDash T$, we have $M\vDash \phi$.
\end{proof}

Completeness follows from Sch\"utte-style argument for the completeness of $\omega$-logic:
\begin{theorem}[Completeness, Girard {\cite[Theorem 10.1.22]{Girard1982Logical2}}]
     If $T\vDash^\sfbeta \phi$, then $T\vdash^\sfbeta\phi$.
\end{theorem}

The completeness will follow from the preproof property:
\begin{theorem}[Preproof property of $\sfbeta$-logic, $\ACA_0$] \label{Theorem: Preproof property of beta logic}
    Let $\calL$ be a primitive recursive language, and $\Gamma$ be a non-empty sequent of $\calL$-sentences.
    Then we can construct a primitive recursive $\sfbeta$-preproof $\pi$ from $\LK$ such that
    \begin{enumerate}
        \item If the $\sfbeta$-preproof $\pi$ is a $\sfbeta$-proof, then $\bigvee \Gamma$ is $\sfbeta$-valid.
        \item Conversely, if $\pi$ is not a $\sfbeta$-proof and $\pi_\alpha$ is ill-founded, then we can recursively build a $\sfbeta$-model $M$ for $\lnot\bigvee \Gamma$ such that $\Ord^M \subseteq \alpha$ from an infinite branch of $\pi_\alpha$.
    \end{enumerate}
\end{theorem}

\modified{
Let us fix a primitive enumeration $\epsilon$ of the $\calL[\omega]$-formulas such that each $\calL[\omega]$-formula occurs infinitely often in $\epsilon$. We will view $\overline{n}$ as an ordinal variable, which represents the $n$th ordinal occurring in the sequent. We also fix, for each type $\tau$ in $\calL$, a primitive recursive enumeration $\epsilon_\tau$ of terms of type $\tau$ in $\calL[\omega]$. 

We will construct a preproof recursively from the last sequent.
We will view sequents as lists of formulas during the construction, so obtaining the desired preproof requires turning the list into a set of formulas. 
For each sequent $\Gamma$ in the proof tree, we assign a dagger to a formula in $\Gamma$. In the initial case, the leftmost formula takes the dagger. We will call the daggered formula the \emph{leading formula.}
We also assign the \emph{stage number} $\varsigma_\Gamma\in\bbN$ to each sequent. The stage number will indicate the number of construction steps we took to reach the current sequent. In the initial case, $\varsigma_\Gamma = 0$. The construction goes as follows:
\begin{enumerate}
    \item If $\Gamma$ is an axiom, then we apply the corresponding rules for the axiom.

    \item Suppose that we are not in the first case.
    The construction is divided into two steps: First, we introduce further sequents on the top, which depend on the leading formula. Second, we introduce new formulas to the sequent by applying Cut. The choice of the formula will depend on the stage number. Then we move the dagger to the next formula in the sequent, and increase the stage number by 1.
    Let us examine the details of the construction. The first step goes as follows:
    \begin{enumerate}
        \item Suppose that the leading formula is atomic. We do nothing in this case. 

        \item Suppose that the leading formula is of the form $A\land B$. Then we introduce sequents as follows:
        \begin{prooftree}
            \AxiomC{$\Gamma, A$}
            \AxiomC{$\Gamma, B$}
            \RightLabel{$\land$}
            \BinaryInfC{$\Gamma$}
        \end{prooftree}

        \item Suppose that the leading formula is of the form $A\lor B$. Then we introduce sequents as follows:
        \begin{prooftree}
            \AxiomC{$\Gamma, A, B$}
            \RightLabel{$\lor$}
            \UnaryInfC{$\Gamma, A$}
            \RightLabel{$\lor$}
            \UnaryInfC{$\Gamma$}
        \end{prooftree}

        \item Suppose that the leading formula is of the form $\forall^\tau x A(x)$ for some $\tau\neq\Ord$. If $v$ is the G\"odel-number least free variable not occurring in $\Gamma$, then we introduce sequents as follows:
        \begin{prooftree}
            \AxiomC{$\Gamma, A(x)$}
            \RightLabel{$\forall^\tau$}
            \UnaryInfC{$\Gamma$}
        \end{prooftree}

        \item Suppose that the leading formula is of the form $\exists^\tau x A(x)$ for some $\tau\neq\Ord$. Let us fix the bijection $h_\Gamma\colon m_\Gamma\to \Occ(\Gamma)$ for some (unique) $m_\Gamma\in \bbN$.
        Let $s_0,\cdots,s_{n-1}$ be the $\calL[\omega]$-terms of type $\tau$ such that $s_i$ is the $<\varsigma_\Gamma$th term in $\epsilon_\tau$, and $\Occ(s_i)\subseteq m_\Gamma$. If we take $t_i = \sfm_{h_\Gamma}(s_i)$, then the new sequents are
        \begin{prooftree}
            \AxiomC{$\Gamma, A(t_0), A(t_1),\cdots,A(t_{n-1})$}
            \RightLabel{$\exists^\tau$}
            \UnaryInfC{$\vdots$}
            \noLine
            \UnaryInfC{$\Gamma, A(t_0), A(t_1)$}
            \RightLabel{$\exists^\tau$}
            \UnaryInfC{$\Gamma, A(t_0)$}
            \RightLabel{$\exists^\tau$}
            \UnaryInfC{$\Gamma$}
        \end{prooftree}
        $n=0$ is possible. We do nothing in this case.
        
        \item Suppose that the leading formula is of the form $\forall^\Ord x A(x)$. Then we introduce sequents as follows:
        \begin{prooftree}
            \AxiomC{$\cdots\ \Gamma, A(\overline{\xi})\ \cdots $ ($\xi<\alpha)$}
            \RightLabel{$\forall^\Ord$}
            \UnaryInfC{$\Gamma$}
        \end{prooftree}
        
        \item Suppose that the leading formula is of the form $\exists^\Ord x A(x)$. If $\Occ(\Gamma) = \{\xi_0<\cdots<\xi_{m-1}\}$, then we introduce sequents as follows:
        \begin{prooftree}
            \AxiomC{$\Gamma, A(\overline{\xi}_0), A(\overline{\xi}_1),\cdots,A(\overline{\xi}_{n-1})$}
            \RightLabel{$\exists^\Ord$}
            \UnaryInfC{$\vdots$}
            \noLine
            \UnaryInfC{$\Gamma, A(\overline{\xi}_0), A(\overline{\xi}_1)$}
            \RightLabel{$\exists^\Ord$}
            \UnaryInfC{$\Gamma, A(\overline{\xi}_0)$}
            \RightLabel{$\exists^\Ord$}
            \UnaryInfC{$\Gamma$}
        \end{prooftree}
        Again, $m=0$ is possible, and we do nothing in this case.
    \end{enumerate}
\end{enumerate}
This finishes the first step. In the second stage, let $\Gamma'$ be one of the topmost sequents we obtain in the first step. 
Let us choose the least $\ge\varsigma$-th formula $C$ in $\epsilon$ such that, if $h_{\Gamma'}\colon m_{\Gamma'}\to \Occ(\Gamma')$ is an order isomorphism, then $\Occ(C)\subseteq m_{\Gamma'}$. Then we introduce $\sfm_{h_{\Gamma'}}(C)$ and its negation by Cut. 
That is, if the leading formula of $\Gamma$ was $\forall^\Ord x A(x)$, then we apply the second step as follows:
\begin{prooftree}
    \AxiomC{$\Gamma, A(\overline{\xi}), \sfm_{h_{\Gamma, A(\overline{\xi})}}(C)$}
    \AxiomC{$\Gamma, A(\overline{\xi}), \lnot \sfm_{h_{\Gamma, A(\overline{\xi})}}(C)$}
    \RightLabel{Cut}
    \BinaryInfC{$\cdots\ \Gamma, A(\overline{\xi})\ \cdots $ ($\xi<\alpha)$}
    \RightLabel{$\forall^\Ord$}
    \UnaryInfC{$\Gamma$}
\end{prooftree}
}
The following lemma shows the homogeneity of $\pi_\alpha$ and the functoriality of $\lag \pi_\alpha\mid \alpha\in \WO\rag$:
\begin{lemma} \label{Lemma: Functoriality of a preproof in the preproof property}
    For every $f\colon \alpha\to\beta$, we have ${}^f(\pi_\beta) = \pi_\alpha$.
\end{lemma}
\begin{proof}
    \modified{
    We prove the following by induction on the sequent from the last sequent: Suppose that $\Occ(\Gamma)\subseteq\ran f$. Then
    \begin{itemize}
        \item $\sfm_f^{-1}(\Gamma)$ occurs in $\pi_\alpha$ of the same stage number.
        \item If $\Gamma$ has conclusion $\Delta$ by inference rule $\sfR$ in $\pi_\beta$, then $\sfm_f^{-1}(\Gamma)$ has conclusion $\sfm_f^{-1}(\Delta)$ by the same inference rule $\sfR$.
        \item The position of the leading formula is preserved under $\sfm_f^{-1}$: That is, if the $i$th formula is the leading formula of $\Gamma$ in $\pi_\beta$, then  the $i$th formula is the leading formula of $\sfm^{-1}_f(\Gamma)$ in $\pi_\alpha$.
    \end{itemize}
    The inductive argument will happen portion-wise on the construction of $\pi_\beta$, so it handles more than one sequents simultaneously. However, every sequent we handle in the portion will have the same stage number. 
    We need the case division by the last inference rule, and we will only examine a few cases.
    \begin{enumerate}
        \item Suppose that $\Gamma$ of stage number $\varsigma$ in $\pi_\beta$ has the leading formula $A\lor B$, so the portion of the proof takes the form
        \begin{prooftree}
            \AxiomC{$\Gamma,A,B,C$}
            \AxiomC{$\Gamma,A,B,\lnot C$}
            \RightLabel{Cut}
            \BinaryInfC{$\Gamma, A, B$}
            \RightLabel{$\lor$}
            \UnaryInfC{$\Gamma, A$}
            \RightLabel{$\lor$}
            \UnaryInfC{$\Gamma$}
        \end{prooftree}
        where
        \begin{itemize}
            \item $C = \sfm_{h_{\Gamma,A,B}}(C_0)$, 
            \item $h_{\Gamma,A,B}\colon m_{\Gamma,A,B}\to \Occ(\Gamma,A,B)$ is the order isomorphism, and 
            \item $C_0$ is the least $\ge\varsigma$th formula in $\epsilon$ with $\Occ(C_0)\subseteq m_{\Gamma,A,B}$. 
        \end{itemize}
        We are assuming that $\Occ(\Gamma)\subseteq\ran f$. Since $A\lor B$ occurs in $\Gamma$ and $\Occ(C)\subseteq\Occ(\Gamma,A,B)$ by definition, we have that $\Occ(\Gamma,A,B,C)\subseteq \Occ(\Gamma)$. In particular, we have $m_{\Gamma} = m_{\Gamma,A,B}=m_{\Gamma,A,B,C}$.
        
        Hence $\Gamma':= \sfm^{-1}_f(\Gamma)$ also occurs in $\pi_\alpha$ with the same stage number $\varsigma$, and the same leading formula, the consequence, and the following inference rule modulo $\sfm^{-1}_f$ in $\pi_\alpha$.
        In particular, the leading formula of $\Gamma'$ in $\pi_\alpha$ is $\sfm^{-1}_f(A\lor B)$.
        Hence, the following portion of the proof in $\pi_\alpha$ is the following:
        \begin{prooftree}
            \AxiomC{$\Gamma',A',B',C'$}
            \AxiomC{$\Gamma',A',B',\lnot C'$}
            \RightLabel{Cut}
            \BinaryInfC{$\Gamma', A', B'$}
            \RightLabel{$\lor$}
            \UnaryInfC{$\Gamma', A'$}
            \RightLabel{$\lor$}
            \UnaryInfC{$\Gamma'$}
        \end{prooftree}
        where
        \begin{itemize}
            \item $A' = \sfm^{-1}_f(A)$ and $B' = \sfm^{-1}_f(B)$,
            \item $C' = \sfm_{h_{\Gamma',A',B'}}(C_1)$, 
            \item $h_{\Gamma',A',B'}\colon m_{\Gamma',A',B'}\to \Occ(\Gamma',A',B')$ is the order isomorphism, and 
            \item $C_1$ is the least $\ge\varsigma$th formula in $\epsilon$ with $\Occ(C_1)\subseteq m_{\Gamma',A',B'}$. 
        \end{itemize}
        We also have $m_{\Gamma'}=m_{\Gamma',A',B'} = m_{\Gamma',A',B',C'}$.
        Since $m_\Gamma=m_{\Gamma'}$, we get $C_0=C_1$.
        From $\sfm^{-1}_f(\Gamma)=\Gamma'$, we have $h_\Gamma=f\circ h_{\Gamma'}$, and so $\sfm_f^{-1}(C)=C'$. This finishes the proof for this case.
        
        \item Now suppose that $\Gamma$ of stage number $\varsigma$ in $\pi_\beta$ has the leading formula $\exists^\tau x A(x)$. The portion of the proof takes the form
        \begin{prooftree}
            \AxiomC{$\Gamma, A(t_0), A(t_1),\cdots,A(t_{n-1}),C$}
            \AxiomC{$\Gamma, A(t_0), A(t_1),\cdots,A(t_{n-1}),\lnot C$}
            \BinaryInfC{$\Gamma, A(t_0), A(t_1),\cdots,A(t_{n-1})$}
            \RightLabel{$\exists^\tau$}
            \UnaryInfC{$\vdots$}
            \noLine
            \UnaryInfC{$\Gamma, A(t_0), A(t_1)$}
            \RightLabel{$\exists^\tau$}
            \UnaryInfC{$\Gamma, A(t_0)$}
            \RightLabel{$\exists^\tau$}
            \UnaryInfC{$\Gamma$}
        \end{prooftree}
        where
        \begin{itemize}
            \item $t_i = \sfm_{h_\Gamma}(s_i)$, where $s_0,\cdots,s_{n-1}$ are the $\calL[\omega]$-terms of type $\tau$ such that $s_i$ is the $<\varsigma_\Gamma$th term in $\epsilon_\tau$, and $\Occ(s_i)\subseteq m_\Gamma$. 
            \item $C$ is defined as in the previous case.
        \end{itemize}
        By the same argument as before, we get 
        \begin{equation*}
            \Occ(\Gamma,A(t_0),\cdots,A(t_{n-1}),C)\subseteq\Occ(\Gamma).
        \end{equation*}
        Hence we have $m_{\Gamma,A(t_0),\cdots,A(t_{n-1}),C} =m_\Gamma$. By the inductive hypothesis, $\Gamma':=\sfm^{-1}_f(\Gamma)$ in $\pi_\alpha$ has the same consequence and the following inference rule, stage number $\varsigma$, and the leading formula modulo $\sfm^{-1}_f$.
        Hence, the following portion of the proof in $\pi_\alpha$ is the following:
        \begin{prooftree}
            \AxiomC{$\Gamma', A'(t'_0), A'(t'_1),\cdots,A'(t'_{n'-1}),C'$}
            \AxiomC{$\Gamma', A'(t'_0), A'(t'_1),\cdots,A'(t'_{n'-1}),\lnot C'$}
            \BinaryInfC{$\Gamma', A'(t'_0), A'(t'_1),\cdots,A'(t'_{n'-1})$}
            \RightLabel{$\exists^\tau$}
            \UnaryInfC{$\vdots$}
            \noLine
            \UnaryInfC{$\Gamma, A'(t'_0), A'(t'_1)$}
            \RightLabel{$\exists^\tau$}
            \UnaryInfC{$\Gamma, A'(t'_0)$}
            \RightLabel{$\exists^\tau$}
            \UnaryInfC{$\Gamma'$}
        \end{prooftree}
        where
        \begin{itemize}
            \item $t'_i = \sfm_{h_{\Gamma'}}(s_i)$, where $s'_0,\cdots,s'_{n'-1}$ are the $\calL[\omega]$-terms of type $\tau$ such that $s'_i$ is the $<\varsigma_\Gamma$th term in $\epsilon_\tau$, and $\Occ(s'_i)\subseteq m_{\Gamma'}$. 
            \item $C'$ is defined as in the previous case.
        \end{itemize}
        We again have $C' =\sfm^{-1}_f(C)$. We also have $n=n'$ and $s_i=s_i'$ for $i<n$ since $s_0',\cdots,s_{n'-1}'$ are determined by $\epsilon_\tau$ and $\varsigma$. This shows $t'_i = \sfm^{-1}_f(t_i)$.
        
        \item Now suppose that $\Gamma$ of stage number $\varsigma$ in $\pi_\beta$ has the leading formula $\forall^\Ord x A(x)$, so the portion of the proof takes the form
        \begin{prooftree}
            \AxiomC{$\Gamma,A(\overline{\xi}),C_\xi$}
            \AxiomC{$\Gamma,A(\overline{\xi}),\lnot C_\xi$}
            \RightLabel{Cut}
            \BinaryInfC{$\cdots\ \Gamma, A(\overline{\xi})\ \cdots \quad (\xi<\beta)$}
            \RightLabel{$\forall^\Ord$}
            \UnaryInfC{$\Gamma$}
        \end{prooftree}
        where $C_\xi = \sfm_{h_{\Gamma,A(\overline{\xi})}}(\hat{C}_\xi)$, and $\hat{C}_\xi$ is the least $\ge\varsigma$th formula in $\epsilon$ with $\Occ(\hat{C}_\xi)\subseteq m_{\Gamma,A(\overline{\xi})}$.
        Note that $\Occ(\Gamma)\subseteq\ran f$ does not imply $\Occ(\Gamma,A(\overline{\xi})) = \Occ(\Gamma)\cup \{\xi\}\subseteq\ran f$. It is not problematic since we ignore $\xi$ with $\xi\notin \ran f$ during mutilation.
        Also, note that
        \begin{equation*}
            m_{\Gamma,A(\overline{\xi})} = 
            \begin{cases}
                m_\Gamma, & \text{if }\xi \in \Occ(\Gamma),\\
                m_{\Gamma}+1, & \text{otherwise.}
            \end{cases}
        \end{equation*}
        By the inductive hypothesis, $\Gamma':=\sfm^{-1}_f(\Gamma)$ has the same stage number in $\pi_\alpha$ with $\Gamma$ in $\pi_\beta$. Also, $\Gamma'$ in $\pi_\alpha$ also has the same leading formula and the following inference rule and consequences modulo $\sfm^{-1}_f$.
        Hence, the corresponding portion of $\pi_\alpha$ is
        \begin{prooftree}
            \AxiomC{$\Gamma',A'(\overline{\eta}),C'_\eta$}
            \AxiomC{$\Gamma',A'(\overline{\eta}),\lnot C'_\eta$}
            \RightLabel{Cut}
            \BinaryInfC{$\cdots\ \Gamma', A'(\overline{\eta})\ \cdots \quad (\eta<\alpha)$}
            \RightLabel{$\forall^\Ord$}
            \UnaryInfC{$\Gamma'$}
        \end{prooftree}
        where 
        \begin{itemize}
            \item $A' = \sfm^{-1}_f(A)$,
            \item $C'_\eta = \sfm_{h_{\Gamma',A'(\overline{\eta})}}(\hat{C}_\eta')$, and
            \item $\hat{C}'_\eta$ is the least $\ge\varsigma$th formula in $\epsilon$ with $\Occ(\hat{C}'_\eta)\subseteq m_{\Gamma',A'(\overline{\eta})}$.
        \end{itemize}
        We can see that if $\xi=f(\eta)$, then $h_{\Gamma,A(\overline{\xi})}=f\circ h_{\Gamma',A'(\overline{\eta})}$ and $m_{\Gamma,A(\overline{\xi})}=m_{\Gamma',A'(\overline{\eta})}$. This shows $\hat{C}'_\eta = \hat{C}_\xi$, and the remaining parts are straightforward. \qedhere 
    \end{enumerate}
    }
\end{proof}

By soundness, if the resulting $\sfbeta$-preproof is a $\sfbeta$-proof, then $T\vDash^\sfbeta \bigvee\Gamma$. 
Conversely, \added{suppose} that $\pi_\alpha$ is ill-founded. Then for each $n$ we can find a sequent $\Gamma_n$ such that $\Gamma_n$ is a hypothesis at level $n$ and $\Gamma_{n+1}$ is a hypothesis of the portion of proof above $\Gamma_n$. \added{Let $\calX = \bigcup_{n<\omega} \Occ(\Gamma_n)$.}
\emph{Throughout the remaining part of the proof, we will fix an infinite branch $\lag\Gamma_n\mid n\in\bbN\rag$.}
Then we can show the following:
\begin{lemma} \label{Lemma: Every formula appears in an infinite branch}
    \added{If $C$ is a formula in $\calL[\alpha]$ whose ordinal parameters belong to $\calX$, then there is $n$ such that precisely one of $C$ or $\lnot C$ appears in $\Gamma_n$ of an infinite branch.}
\end{lemma}
\begin{proof} \modified{
    Note that both $C$ and $\lnot C$ can occur in $\Gamma_n$, otherwise the branch terminates. Now we claim that one of $C$ or $\lnot C$ must appear in an infinite branch. Pick a large $k$ such that $\Occ(C)\subseteq \ran h_{\Gamma_k}$, then choose a large $n>k$ such that 
    \begin{itemize}
        \item $\varsigma_{\Gamma_{n-1}} <\varsigma_{\Gamma_n}$, and
        \item $\sfm^{-1}_{h_{\Gamma_n}}(C)$ is the $\varsigma_{\Gamma_n}$th formula of $\epsilon$.
    \end{itemize}
    The first condition implies that $\Gamma_{n+1}$ and a few sequents after $\Gamma_{n+1}$ are introduced by the construction scheme of $\pi_\alpha$. From $\Occ(C)\subseteq \ran h_{\Gamma_n}$, the construction of $\pi_\alpha$ indicates that the portion of $\pi_\alpha$ after $\Gamma_n$ is the following:
    \begin{prooftree}
        \AxiomC{$\Gamma_n',C$}
        \AxiomC{$\Gamma_n',\lnot C$}
        \RightLabel{Cut}
        \BinaryInfC{$\vdots$}
        \UnaryInfC{$\Gamma_n$}
    \end{prooftree}
    We have either $\Gamma_n',C = \Gamma_{n+1}$ or $\Gamma_n',\lnot C = \Gamma_{n+1}$, as desired.}
\end{proof}

\modified{Now let us construct a desired model from the infinite branch $\lag \Gamma_n\mid n<\omega\rag$ as follows:
\begin{enumerate}
    \item For a sort $\tau$, the set of elements of sort $\tau$ is the set of all terms.
    \item $\Ord^M = \calX\subseteq\alpha$.
    \item If $f$ is a function symbol, we define $f^M(t_0,\cdots,t_{n-1}) = f(t_0,\cdots,t_{n-1})$.
    \item $M\vDash A$ iff $\lnot A$ occurs in $M$.
\end{enumerate}
The following lemma shows that $M$ is really a model:
\begin{lemma} \label{Lemma: Preproof property gives a model}
For each $n\in\bbN$, let $A$ be the leading formula of $\Gamma_n$. 
    \begin{enumerate}
        \item If $A \equiv B\land C$, then either $B$ or $C$ is in $\Gamma_{n+1}$.
        \item If $A \equiv B\lor C$, then both $B$ and $C$ are in $\Gamma_{n+1}$.
        \item If $A \equiv \forall^\tau x B(x)$ for $\tau\neq \Ord$, then there is a free variable $x$ of type $\tau$ such that $B(x)$ is in $\Gamma_{n+1}$.
        \item If $A \equiv \exists^\tau x B(x)$ for $\tau\neq \Ord$, then for every term $t \in \calL_\calX$ of type $\tau$, there is $m>n$ such that $B(t)$ occurs in $\Gamma_m$.
        \item If $A\equiv \forall^\Ord x B(x)$, then there is $\xi\in\calX$ such that $B(\overline{\xi})$ is in $\Gamma_{n+1}$.
        \item If $A\equiv \exists^\Ord x B(x)$, then for every $\xi\in\calX$ there is $m>n$ such that $B(\overline{\xi})$ occurs in $\Gamma_m$.
    \end{enumerate}
\end{lemma}
\begin{proof}
    The cases for $\land$, $\lor$ and $\forall^\tau$ are easy to handle. The case for $\exists^\Ord$ is analogous to the case for $\exists^\tau$ with $\tau\neq\Ord$, so we only examine this case.
    For a given term $t$ of type $\tau$ with $\Occ(t)\subseteq \calX$, let us find a sufficiently large $n$ in a way that
    \begin{enumerate}
        \item $\varsigma_{\Gamma_{n-1}} < \varsigma_{\Gamma_n}$,
        \item $\Occ(t) \subseteq \Occ(\Gamma_n)$, 
        \item The leading formula of $\Gamma_n$ is $\exists^\tau x B(x)$.
        \item $\sfm^{-1}_{h_{\Gamma_n}}(t)$ is the $<\varsigma_{\Gamma_n}$-th term in $\epsilon$, where $h_{\Gamma_n}\colon m_{\Gamma_n}\to \Occ(\Gamma_n)$ is the unique isomorphism for some $m_{\Gamma_n}<\omega$.
\end{enumerate}
        Then the construction of $\pi_\alpha$ guarantees that there is $n'>n$ such that $\Gamma_{n'}$ contains $B(t)$.
\end{proof}}
\added{Hence $M$ is a model of $\lnot\bigvee \Gamma$}. It finishes the proof of the preproof property of $\sfbeta$-logic.

\section{Genedendrons} \label{Section: Genedendrons}
In this section, we introduce the notion of \emph{genedendrons}.

\begin{definition}
    A \emph{semi-genedendron} is a pair $(D,\varrho)$ satisfying the following:
    \begin{enumerate}
        \item $D$ is a prequasidendroid.
        \item $\varrho = \lag\varrho_\alpha\colon D_\alpha^*\to \bbN^{<\omega}\mid \alpha\in\WO\rag$ is a functorial continuous map taking a finite branch of $D_\alpha^*$ and returning a finite sequence in the following sense:
        \begin{enumerate}
            \item If $\sigma,\tau\in D_\alpha^*$, $\sigma\subseteq\tau$, then $\varrho_\alpha(\sigma)\subseteq \varrho_\alpha(\tau)$.
            \item For an increasing map $f\colon\gamma\to\alpha$, if $\sigma\in D(\gamma) = {}^f D(\alpha)$, then $\varrho_\alpha(f(\sigma)) = \varrho_\gamma(\sigma)$.
        \end{enumerate}
        We say $\varrho$ the \emph{extraction scheme}.
    \end{enumerate}
    We say $(D,\varrho)$ is \emph{locally well-founded} if $D$ is locally well-founded.
    A semi-genedendron $(D,\varrho)$ is \emph{genedendron}\footnote{Ancient Greek `Birth tree' or `Generating tree'} if for every ordinal $\alpha$ and an infinite branch $b$ of $D_\alpha^*$, $\varrho_\alpha(b):=\bigcup_{n\in\bbN} \varrho_\alpha(b\restricts n)$ is constant. If $\varrho_\alpha(b)$ has a constant value $R$, we say \emph{$(D,\varrho)$ generates $R$.}%
    \footnote{There is a small prickle in the definition of a genedendron: The only examples of genedendrons appearing in this paper come from functorial trees (e.g. \cite[Definition 4.1]{GirardNormann1992Embeddability}) and the author does not think there are other important examples. It might be better to define genedendrons in terms of functorial trees, but we will stick to prequasidendroids to avoid developing the theory of functorial trees.}
\end{definition}

Note that the definition of a (semi-)genedendron does not include that $D_\alpha$ has an infinite branch for some $\alpha$. The following lemma says, like a prequasidendroid, a semi-genedendron is determined by its small fragment:

\begin{lemma}
    Every semi-genedendron $(D,\varrho)$ is determined by $(D_\omega,\varrho_\omega)$.
\end{lemma}
\begin{proof}
    $D$ being a prequasidendroid implies $D$ is determined by $D_\omega$ by extension.
    For $\varrho$, observe that for every $\sigma\in D_\alpha^*$ there is an increasing map $f\colon n\to \alpha$ and $\tau\in D_n^*$ such that $f[\tau]=\sigma$.
    By the definition, we have $\varrho_n(\tau) = \varrho_\alpha(\sigma)$.
\end{proof}
The previous lemma says the definition of semi-genedendron and genedendron can be formulated over a fragment of second-order arithmetic. Of course, we should replace ordinals with well-orders in the definition of a genedendron to fit the definition into the second-order arithmetic.
Also, we can talk about a (semi-)genedendron being recursive: A (semi-)genedendron $(D,\varrho)$ is recursive if $(D_\omega,\varrho_\omega)$ is recursive.

The next proposition says genedendron is a way to formulate $\Sigma^1_2$-definable reals:
\begin{proposition}[$\ACA_0$] \label{Proposition: Every Sigma 1 2 definable real is generated by a genedendron}
    Let $\phi(X)$ be a $\Sigma^1_2$-formula defining a real $R$. Then we can find a recursive locally well-founded genedendron $(D,\varrho)$ such that there is a well-order $\alpha$ and an infinite branch $B$ over $D(\alpha)$ such that $R = \varrho_\alpha(B)$.
\end{proposition}
\begin{proof}
    Its proof follows from a usual Shoenfield tree construction:
    Let $\psi(Y,X)$ be a $\Pi^1_1$-formula such that $\phi(X) \equiv \exists^1 Y\psi(Y,X)$. Then we can effectively associate a recursive tree $T$ over $\omega\times\omega\times\omega$ such that
    \begin{equation*}
        \psi(Y,X) \iff \forall^1 Z (\lag Z,Y,X\rag\notin [T]).
    \end{equation*}
    Here $[T]$ is the set of infinite branches of $T$.
    Fix a standard enumeration $\lag \bfs_i\mid i\in\bbN\rag$ of $\bbN^{<\omega}$, and for a fixed well-order $\alpha$ define
    \begin{multline*}
        \hat{T}(\alpha) = \{\lag \vec{t},\vec{u},\vec{\zeta}\rag\in \bbN^{<\omega}\times\bbN^{<\omega}\times\alpha^{<\omega} \mid |\vec{t}|=|\vec{u}|=|\vec{\zeta}| \\
        \land \forall i,j<|\vec{t}| [\lag \bfs_i,\vec{t}\restricts |\bfs_i|, \vec{u}\restricts |\bfs_i| \rag\in T\land \bfs_i\supsetneq \bfs_j\lr \zeta_i<\zeta_j]\}.
    \end{multline*}
    Then let us turn $\hat{T}(\alpha)$ to an $\alpha$-prequasidendroid $D_\alpha$ over $(\bbN^2,\le_\mathsf{lex})$ as follows: 
    \begin{equation*}
        D_\alpha = \{\lag (t_0,u_0),\zeta_0,\cdots, (t_{m-1},u_{m-1}),\zeta_{m-1}, (0,0)\rag \mid \lag \vec{t},\vec{u},\vec{\zeta}\rag\in \hat{T}(\alpha),\  |\vec{t}|=m\}.
    \end{equation*}
    We claim that for every increasing map $f\colon\alpha\to\beta$, we have ${}^f D_\beta = D_\alpha$, which proves the homogeneity of $D_\beta$: It follows from proving the equivalence
    \begin{equation*}
        \lag \vec{t},\vec{u},\vec{\zeta}\rag\in \hat{T}(\alpha) \iff \lag \vec{t},\vec{u},f(\vec{\zeta})\rag\in \hat{T}(\beta)
    \end{equation*}
    for $\vec{t},\vec{u}\in \bbN^{<\omega}$ and $\vec{\zeta}\in \alpha^{<\omega}$, which is immediate from the definition of $\hat{T}(\alpha)$.
    
    Then define $\varrho_\alpha \colon D_\alpha^*\to \bbN^{<\omega}$ by
    \begin{itemize}
        \item $\varrho_\alpha(\lag (t_0,u_0),\zeta_0,\cdots, (t_{m-1},u_{m-1}),\zeta_{m-1}\rag) = \lag u_i\mid i<m\rag$, and
        \item $\varrho_\alpha(\lag (t_0,u_0),\zeta_0,\cdots, (t_{m-1},u_{m-1}),\zeta_{m-1}, (t_m,u_m)\rag) = \lag u_i\mid i\le m\rag$.
    \end{itemize}
    $\varrho_\alpha$ does not depend on the ordinal part of the input, and it implies $(D,\varrho)$ is a semi-genedendron.

    Now we claim that $(D,\varrho)$ is a genedendron: Suppose that there is an infinite branch
    \begin{equation*}
        B = \lag (t_0,u_0),\zeta_0, (t_1,u_1),\zeta_1, \cdots \rag
    \end{equation*}
    of $D_\alpha$. 
    It means $\lag\vec{t},\vec{u},\vec{\zeta}\rag$ is an infinite branch of $\hat{T}(\alpha)$.
    Hence by definition of $\hat{T}(\alpha)$, so the map $\bfs_i\mapsto \zeta_i$ is a function ranking the tree $T_{\vec{t},\vec{u}} = \{s\in\bbN^{<\omega}\mid \lag s,\vec{t}\restricts|s|,\vec{u}\restricts|s| \rag\in T\}$. Hence $\forall^1 Z \lag Z,\vec{t},\vec{u}\rag \notin [T]$, which is equivalent to $\psi(\vec{t},\vec{u})$.
    Thus $\exists^1 Y\psi(Y,\vec{u})$, so $R=\vec{u}$. By the definition of $\varrho_\alpha$, $\bigcup_{n\in\bbN} \varrho_\alpha(b\restricts n) = \vec{u}=R$, as desired.
    
    Lastly, it is clear that $\hat{T}(\omega)$ is recursive, so is $D_\omega$. $\varrho_\omega$ is $\Delta^0_1$-definable, so recursive. The local well-foundedness of $D_\alpha$ follows from that the term set $\bbN^2$ is well-ordered. 
\end{proof}

\autoref{Proposition: Every Sigma 1 2 definable real is generated by a genedendron} shows every $\Sigma^1_2$-\added{singleton} real is generated by a genedendron. It is easy to see that every real given by a genedendron is \added{a} $\Sigma^1_2$-\added{singleton}. 
However, the genedendron provided in the proof of \autoref{Proposition: Every Sigma 1 2 definable real is generated by a genedendron} has little practical use. Genedendrons we can compute by hand come from $\sfbeta$-preproofs over a suitable language and axioms. 

Observe that for a locally well-founded prequasidendroid $D$, $D(\alpha)$ has an infinite branch iff $\Lin^*(D)(\alpha)$ is ill-founded. Hence the following is immediate:

\begin{lemma} \pushQED{\qed}
    Let $(D,\varrho)$ be a locally well-founded semi-genedendron. Then $\Clim(\Lin^*(D))$ is the least ordinal $\alpha$ such that $D(\alpha)$ has an infinite branch. \qedhere 
\end{lemma}

We will frequently use the climax of the linearization of a genedendron to `rank' a genedendron. Thus for a notational easiness, we use $\Clim(D)$ to denote $\Clim(\Lin^*(D))$ for a genedendron $(D,\varrho)$.
We can use the climax of a genedendron to rank a $\Sigma^1_2$-\added{singleton} real as follows:
\begin{definition} \label{Definition: Sigma 1 2 altitude}
    Let $R$ be a $\Sigma^1_2$-definable real. Let us define the \emph{$\Sigma^1_2$-altitude of $R$} by
    \begin{equation*}
        \Alt_{\Sigma^1_2}(R) = \min \{\Clim(D) \mid \text{$(D,\varrho)$ is a locally well-founded recursive genedendron computing $R$}\}.
    \end{equation*}
\end{definition}

Let us compute $\Alt_{\Sigma^1_2}$ for some reals:

\begin{example}
    Suppose that $R$ is an arithmetical real, that is, there is an arithmetical formula $\phi(n)$ such that
    \begin{equation*}
        \forall^0 n [R(n)=0\lr \lnot\phi(n)]\land [R(n)=1 \lr \phi(n)].
    \end{equation*}
    We claim that $\Alt_{\Sigma^1_2}(R)=0$. Consider the $\omega$-logic with the new extra unary predicate $X$, and consider the Sch\"utte-styled $\omega$-preproof for the sentence
    \begin{equation} \label{Formula: Example - arithmetical real rank}
        \exists^0 n [(\phi(n)\lor n\in X) \land (\lnot\phi(n)\lor n\notin X)].
    \end{equation}
    Then the $\omega$-preproof has an infinite branch $b$, which encodes an $\omega$-model satisfying the negation of \eqref{Formula: Example - arithmetical real rank}. We can view the $\omega$-preproof as a constant prequasidendroid $D$, and define $\varrho_\alpha$ by
    \begin{equation*}
        \varrho_\alpha(\sigma) = \{(n,0) \mid \ulcorner \underline{n}\notin X\urcorner \in \sigma\} \cup \{(n,1) \mid \ulcorner \underline{n}\in X\urcorner \in \sigma\}.
    \end{equation*}
    Then $\varrho_0(b)=r$. Since $D=D_\alpha$ always has an infinite branch, $\Clim(D)=0$.
\end{example}

\subsection{Another example: Computing the \texorpdfstring{$\Sigma^1_2$}{Sigma 1 2}-altitude of the hyperjump of $\emptyset$} \label{Subsection: Hyperjump 0 - Sigma 1 2 alt}
In this subsection, we prove that $\rank_{\Sigma^1_2}(\HJ(\emptyset)) = \omega_1^\CK$, where $\HJ(X)$ is the hyperjump of $X$. 
For technical convenience, we assume that $\HJ(X)$ has a form of characteristic function.

One direction follows from an easy model-theoretic argument:
\begin{lemma}
    $\Alt_{\Sigma^1_2}(\HJ(\emptyset)) \ge \omega_1^\CK$.
\end{lemma}
\begin{proof}
    Suppose that $M$ is a $\sfbeta$-model of second-order arithmetic decoded from $\HJ(\emptyset)$ and $(D,\varrho)$ is a locally well-founded recursive genedendron generating $\HJ(\emptyset)$. Then $(D,\varrho)\in M$, and $M$ also thinks $(D,\varrho)$ is a locally well-founded genedendron. However, $M$ thinks $(D,\varrho)$ does not generate any real (otherwise $\HJ(\emptyset)\in M$.) 
Hence $M$ thinks $D_\alpha$ has no infinite branch for every well-order $\alpha\in M$, and every well-order in $M$ has ordertype $<\omega_1^\CK$. This shows $\Clim(D)\ge \omega_1^\CK$. Since $(D,\varrho)$ is arbitrary, we have the desired lower bound.
\end{proof}

$\Alt_{\Sigma^1_2}(\HJ(\emptyset)) \le \omega_1^\CK$ is mildly challenging to establish. We will construct a locally well-founded recursive genedendron $(D,\varrho)$ generating $\HJ(\emptyset)$ such that $\Clim(\Lin^*(D))=\omega_1^\CK$, and here is where we cast $\sfbeta$-logic.

Let us fix a $\Pi^0_1$ formula $P(X,x)$ defining a monotone operator whose fixpoint is $\HJ(\emptyset)$. That is, if we let $\Phi(X) = \{x\in\bbN \mid P(X,x)\}$, and define $I_\Phi^\alpha(X) = \bigcup_{\beta<\alpha} \Phi(I_\Phi^\beta(X))$, then $\HJ(\emptyset)$ is the characteristic function for $I_\Phi^{\omega_1^\CK}(\emptyset)$.

Now let us consider a variant of the inductive logic $T$, \added{which is a variant of Girard's one} \cite[\S 11.3]{Girard1982Logical2}.
The language $\calL$ of $T$ comprises the language of \emph{first order arithmetic} $\{0,S,+,\cdot,\le\}$ \added{of type} $\bbN$, plus a new type $\Ord$ and a binary relation $\le^\Ord$, a binary relation $I_\Phi$ of type $\Ord\times\bbN$, and a unary relation $\overline{I}_\Phi$ of type $\bbN$. For an ordinal $\alpha$, the axioms and the rules for $T[\alpha]$ are the rules of $\sfbeta$-logic plus the following, 
\begin{enumerate}\modified{
    \item The axioms and rules for $\omega$-logic, that is, 
    \begin{center}
        \begin{longtable}{c c}
            \AxiomC{}
            \RightLabel{=}
            \UnaryInfC{$\Gamma,t=t$}
            \DisplayProof
            &
            \\[1em] 
            
            \AxiomC{}
            \RightLabel{$\neq_0$}
            \UnaryInfC{$\Gamma,t\neq u, \lnot A(t), A(u)$}
            \DisplayProof
            &
            \AxiomC{}
            \RightLabel{$\neq_1$ ($n\neq m$)}
            \UnaryInfC{$\Gamma,S^n0 \neq S^m0 $}
            \DisplayProof
            \\[2em]

            \AxiomC{$\cdots \Gamma, A(S^n0) \cdots$ ($n\in\bbN$)}
            \RightLabel{$\forall^\bbN$}
            \UnaryInfC{$\Gamma, \forall^\bbN x A(x)$}
            \DisplayProof
            &
            \AxiomC{$\Gamma, A(t)$}
            \RightLabel{$\exists^\bbN$}
            \UnaryInfC{$\Gamma,  \exists^\bbN x A(x)$}
            \DisplayProof
        \end{longtable}
        \vspace{-1.8em}
    \end{center}
    Here $t,u$ are $\bbN$-terms.
    
    \item The rules for the `operator iteration' $I_{\Phi}$:
    \begin{center}
        \begin{longtable}{c c}
            \AxiomC{$\cdots \Gamma, \overline{\xi}\le^\Ord\overline{\eta}, \lnot P(I_{\Phi}(\overline{\eta},\cdot),t) \cdots$ ($\eta\in\alpha$)}
            \RightLabel{$\lnot I_{\Phi}$}
            \UnaryInfC{$\Gamma, \lnot I_{\Phi}(\overline{\xi},t)$}
            \DisplayProof
            &
            \AxiomC{$\Gamma, P(I_{\Phi}(\overline{\eta},\cdot),t)$}
            \RightLabel{$I_{\Phi}$ ($\eta<\xi$)}
            \UnaryInfC{$\Gamma, I_{\Phi}(\overline{\xi},t)$}
            \DisplayProof
        \end{longtable}
        \vspace{-1.8em}
    \end{center}
    Here $I_{\Phi}$ can be applied only when $\eta<\xi\in\alpha$, and $\lnot I_{\Phi}$ requires proofs for every $\eta<\xi\in\alpha$. 
    $P(I_{\Phi}(\overline{\eta},\cdot),t)$ is a formula obtained from $P(\added{X,t})$ by replacing $X(x)$ with $I_{\Phi}(\overline{\eta},x)$.

    \item The rules for the `operator fixedpoint' $\overline{I}_\Phi$:
    \begin{center}
        \begin{longtable}{c c}
            \AxiomC{$\cdots \Gamma, \lnot P(I_{\Phi}(\overline{\eta},\cdot),t) \cdots$ ($\eta\in\alpha$)}
            \RightLabel{$\lnot \overline{I}_{\Phi}$}
            \UnaryInfC{$\Gamma, \lnot \overline{I}_{\Phi}(t)$}
            \DisplayProof
            &
            \def\fCenter{\mbox{ $\vdash$ }}
            \AxiomC{$\Gamma, P(I_{\Phi}(\overline{\eta},\cdot),t)$}
            \RightLabel{$\overline{I}_{\Phi}0$ ($\eta\in\alpha$)}
            \UnaryInfC{$\Gamma, \overline{I}_{\Phi}(t)$}
            \DisplayProof
            \\[2em]
            &
            \def\fCenter{\mbox{ $\vdash$ }}
            \AxiomC{$\Gamma, P(\overline{I}_{\Phi}(\cdot),t)$}
            \RightLabel{$\overline{I}_{\Phi}$1}
            \UnaryInfC{$\Gamma, \overline{I}_{\Phi}(t)$}
            \DisplayProof
        \end{longtable} 
        \vspace{-1.8em}
    \end{center}}
\end{enumerate}

Now let again $M$ be the $\sfbeta$-model (\`a la Girard) decoded from $\HJ(\emptyset)$. We can turn $M$ to an $\calL$-structure as follows: We drop the second-order part and $\Emb^M$ of $M$. Then define
\begin{equation*} \textstyle 
    I_\Phi^\xi = \bigcup_{\eta<\xi} \Phi(I_\Phi^\eta),
\end{equation*}
where $\Phi(X) = \{x\in\bbN\mid P(X,x)\}$, and interpret $I_\Phi(\xi,n)$ iff $n\in I_\Phi^\xi$. Also, define $\overline{I}_\Phi(n)$ iff $n\in I_\Phi^{\omega_1^\CK}$.
Then we can prove the following lemma by induction on a proof:
\begin{lemma} \label{Lemma: Inductive logic for omega1CK soundness}
    Every sequent $\omega_1^\CK$-provable from $T$ is valid over $M$ in the following sense: If $\Gamma$ has an $\omega_1^\CK$-proof, then $M\vDash \bigvee\Gamma$.
\end{lemma}

Now we modify the proof of the $\sfbeta$-completeness theorem to construct a genedendron: The essential part of the proof for the preproof property, giving a primitive recursive $\sfbeta$-preproof for the given sequent.
We add the following into the construction of $\pi_\alpha$:
\begin{enumerate} \modified{
    \item Suppose that the leading formula of $\Gamma$ is $\forall^\bbN x B(x)$. The portion (before the Cut rule) is
    \begin{prooftree}
        \AxiomC{$\cdots \Gamma, B(S^n0) \cdots$ for all $n\in\bbN$}
        \RightLabel{$\forall^\bbN$}
        \UnaryInfC{$\Gamma$}
    \end{prooftree}

    \item Suppose that the leading formula of $\Gamma$ is $\exists^\bbN x B(x)$. The portion (before the Cut rule) is
    \begin{prooftree}
        \AxiomC{$\Gamma, B_0(0),\cdots,B_0(S^{\varsigma-1}0)$}
        \RightLabel{$\exists^\bbN$}
        \UnaryInfC{$\Gamma, \fCenter B(0), B(1)$}
        \RightLabel{$\exists^\bbN$R}
        \UnaryInfC{$\vdots$}
        \noLine
        \UnaryInfC{$\Gamma, \fCenter B(0)$}
        \RightLabel{$\exists^\bbN$}
        \UnaryInfC{$\Gamma$}
    \end{prooftree}
    Here, $\varsigma = \varsigma_\Gamma$ is the stage number of $\Gamma$.

    \item If the leading formula of $\Gamma$ is $I_\Phi(\overline{\xi},t)$, then the portion is
    \begin{prooftree}
        \AxiomC{$\Gamma, P(I_{\Phi}(\overline{\eta}_0,\cdot),t),\cdots,P(I_{\Phi}(\overline{\eta}_{q'-1},\cdot),t)$}
        \RightLabel{$I_{\Phi}$}
        \UnaryInfC{$\vdots$}
        \noLine
        \UnaryInfC{$\Gamma, P(I_{\Phi}(\overline{\eta}_0,\cdot),t)$}
        \RightLabel{$I_{\Phi}$}
        \UnaryInfC{$\Gamma$}
    \end{prooftree}
    Here, $\eta_0,\cdots,\eta_{q-1}$ is the increasing enumeration of $\{\eta\in \Occ(\Gamma)\mid \eta<\xi\}$.

    \item If the leading formula of $\Gamma$ is $\lnot I_\Phi(\overline{\xi},t)$, then the portion is
    \begin{prooftree}
        \AxiomC{$\cdots \Gamma, \overline{\xi}\le^\Ord\overline{\eta}, \lnot P(I_{\Phi}(\overline{\eta},\cdot),t) \cdots$ ($\eta\in\alpha$)}
        \RightLabel{$\lnot I_{\Phi}$}
        \UnaryInfC{$\Gamma$}
    \end{prooftree}

    \item If the leading formula of $\Gamma$ is $\overline{I}_\Phi(t)$, then the portion is
    \begin{prooftree}
        \AxiomC{$\Gamma,P(I_{\Phi}(\overline{\eta}_0,\cdot),t),\cdots,P(I_{\Phi}(\overline{\eta}_{q-1},\cdot),t), P(\overline{I}_\Phi,t)$}
        \RightLabel{$\overline{I}_{\Phi}1$}
        \UnaryInfC{$\Gamma, P(I_{\Phi}(\overline{\eta}_0,\cdot),t),\cdots,P(I_{\Phi}(\overline{\eta}_{q'-1},\cdot),t)$}
        \RightLabel{$\overline{I}_{\Phi}0$}
        \UnaryInfC{$\vdots$}
        \noLine
        \UnaryInfC{$\Gamma, P(I_{\Phi}(\overline{\eta}_0,\cdot),t)$}
        \RightLabel{$\overline{I}_{\Phi}0$}
        \UnaryInfC{$\Gamma$}
    \end{prooftree}
    Here, $\eta_0,\cdots,\eta_{q-1}$ is the increasing enumeration of $\{\eta\in \Occ(\Gamma)\mid \eta<\xi\}$.

    \item If the leading formula of $\Gamma$ is $\lnot \overline{I}_\Phi(t)$, then the portion is
    \begin{prooftree}
        \AxiomC{$\cdots \Gamma, \lnot P(I_{\Phi}(\overline{\eta},\cdot),t) \cdots$ ($\eta\in\alpha$)}
        \RightLabel{$\lnot I_{\Phi}$}
        \UnaryInfC{$\Gamma$}
    \end{prooftree}}
\end{enumerate}

Now let $\pi_\alpha$ be the $\alpha$-preproof for the empty sequent generated by the previous scheme. We can see that ${}^f \pi_\beta=\pi_\alpha$ for every increasing $f\colon\alpha\to\beta$, which guarantees the homogeneity of the resulting $\alpha$-preproof $\pi_\alpha$. Hence, we get a $\sfbeta$-preproof $\lag\pi_\alpha\mid \alpha\in\Ord\rag$ for the empty sequent. 
We need to check that the proof for the preproof property also works for the new proof system. \modified{However, we can easily see that the proof of \autoref{Lemma: Functoriality of a preproof in the preproof property}, \ref{Lemma: Every formula appears in an infinite branch}, and \ref{Lemma: Preproof property gives a model} also works for the current system.}

We are ready to define a genedendron for $\HJ(\emptyset)$.
\emph{Throughout the remaining subsection, let us fix $\pi$ by the $\sfbeta$-preproof for the sequent $S0=0$.}
Then consider the following extraction scheme $\varrho_\alpha$: If $\pi_\alpha = (D_\alpha,\phi_\alpha)$, then
\begin{equation*}
    \varrho_\alpha(\sigma) = \{(n,1) \mid \lnot \overline{I}_{\Phi}(S^n0)\text{ \added{appears in} }\proj_1(\phi_\alpha(\sigma))\} 
    \cup \{(n,0) \mid  \overline{I}_{\Phi}(S^n0)\text{ \added{appears in } }\proj_1(\phi_\alpha(\sigma))\}
\end{equation*}
That is, $\varrho_\alpha(\sigma)$ is a finite partial function obtained from the $\sigma$th sequent of $\pi_\alpha$, by \added{checking how} $\overline{I}_{\Phi}$ \added{is placed in a sequent}, which asserts `$n$ is in $I_\Phi^{<\alpha}=\bigcup_{\xi\in\alpha} I_{\Phi}^\xi(\emptyset)$.' Then $n\in I_\Phi^{<\alpha}$ if the sentence occurs in the antecedent, and $n\notin I_\Phi^{<\alpha}$ otherwise. \added{We understand ``$\lnot\overline{I}_\Phi(S^n0)$'' is in $\Gamma$'' as $n\in \overline{I}_\Phi$ since we catch a countermodel generated from an infinite branch of a preproof.}
It defines a semi-genedendron, which is clearly locally well-founded. Since $\pi$ is recursive, so is $(D,\varrho)$. 

\begin{lemma}
    $(D,\varrho)$ is a genedendron generating $\HJ(\emptyset)$.
\end{lemma}
\begin{proof}
    Suppose that $D_\alpha$ has an infinite branch, which gives an infinite branch $\lag\Gamma_n\mid n\in\bbN\rag$ of $\pi_\alpha$.
    By the proof of the preproof property of the $\sfbeta$-logic, the infinite branch gives a $\sfbeta$-model $N$ of height $\le\alpha$ that satisfies $\lnot\bigvee\Gamma_n$. 
    Furthermore, the $\bbN$-part of $N$ is the standard set of natural numbers since we can match natural numerals to standard natural numbers.

    Now let $(I_\Phi^\xi)^N := \{m\mid N\vDash I_\Phi(\overline{\xi},S^m0)\}$ and $\overline{I}_\Phi^N := \{m\mid N\vDash \overline{I}_\Phi(S^m0)\}$. Let us make some observations for $(I_\Phi^\xi)^N$:
    \begin{itemize}
        \item If $\lnot I_\Phi(\overline{\xi}, S^m0)$ occurs in an infinite branch, then $P(I_\Phi(\overline{\eta},\cdot), S^m0)$ occurs in an  infinite branch for some $\eta<\xi$. This shows $(I_\Phi^\xi)^N \subseteq \bigcup_{\eta<\xi} \Phi((I_\Phi^\eta)^N)$.
        \item Now suppose that $I_\Phi(\overline{\xi}, S^m0)$ occurs in an infinite branch. If $X_n$ is the set of ordinals occurring in the $n$th sequent of the infinite branch, then $\bigcup_{n\in\bbN} X_n = \Ord^N$.
        Moreover, for a fixed stage number $\varsigma$, $I_\Phi(\overline{\xi}, S^m0)$ occurs in a sequent with the stage number $\varsigma$.

        Fix $\eta<\xi\in\alpha$ such that $\eta\in \bigcup_{n<\omega} X_n$. We can choose a large $n$ such that $\eta\in X_n$ and the leading formula of the $n$th sequent of the infinite branch is $I_\Phi(\overline{\xi}, S^m0)$.
        Then the definition of $\pi_\alpha$ adds a formula $I_\Phi(\overline{\eta},S^m0)$ into later sequents of the infinite branch. 
        This shows 
        \begin{equation*} \textstyle
            m\notin (I^\xi_\Phi)^N \implies m\notin\bigcup_{\eta<\xi} \Phi((I^\eta_\Phi)^N),
        \end{equation*}
        or in other words, $\bigcup_{\eta<\xi} \Phi((I_\Phi^\eta)^N) \subseteq (I_\Phi^\xi)^N$.
    \end{itemize}

    Hence for each $\xi\in \Ord^N$, we have $\bigcup_{\eta<\xi} \Phi((I_\Phi^\eta)^N) = (I_\Phi^\xi)^N$. By a similar reasoning, we also have $\overline{I}_\Phi^N = \bigcup_{\xi\in\Ord^N}\Phi((I_\Phi^\xi)^N)$.
    Now we claim that $\Phi(\overline{I}_\Phi^N)\subseteq \overline{I}_\Phi^N$ by showing
    \begin{equation*}
        m \notin \overline{I}_\Phi^N \implies m\notin \Phi(\overline{I}_\Phi^N).
    \end{equation*}
    Suppose that $\overline{I}_\Phi(S^m0)$ occurs in a sequent of the infinite branch. Then the definition of $\pi_\alpha$ guarantees that $P(\overline{I}_\Phi(\cdot),S^m0)$ also occurs in another sequent. 
    It implies $N\nvDash P(\overline{I}_\Phi(\cdot),S^m0)$ holds. Since $N$ is an $\omega$-model, it implies $\lnot P(\overline{I}_\Phi^N,m)$, i.e., $m\notin \Phi(\overline{I}_\Phi^N)$.
    
    Combining all the facts we have proven, we have $I^\xi_\Phi = (I^\xi_\Phi)^N$ for all $\xi\in \Ord^N$, $\overline{I}_\Phi^N = I_\Phi^{\Ord^N}$, and $\overline{I}_\Phi^N$ is closed under $\Phi$. 
    This shows the characteristic function for $\overline{I}_\Phi^N$ is $\HJ(\emptyset)$.
\end{proof}

\begin{lemma}
    $D_{\omega_1^\CK}$ is ill-founded, so $\Clim(D)\le\omega_1^\CK$.
\end{lemma}
\begin{proof}
    Suppose not, assume that $D_{\omega_1^\CK}$ is well-founded.
    Then $\pi_{\omega_1^\CK}$ is an $\omega_1^\CK$-proof for the empty sequent. Now let $M$ be the $\sfbeta$-model decoded from $\HJ(\emptyset)$.
    By \autoref{Lemma: Inductive logic for omega1CK soundness}, $M$ satisfies $S0=0$, a contradiction.
\end{proof}

Hence $(D,\varrho)$ is a genedendron of climax $\omega_1^\CK$ generating $\HJ(\emptyset)$.
As a final remark, proving $(D,\varrho)$ being a genedendron is a theorem of $\ACA_0$, where $(D,\varrho)$ is what is constructed from the $\sfbeta$-preproof property.
Proving $\Clim(D)=\omega_1^\CK$ is not precisely a theorem of $\ACA_0$ since its statement uses ordinals. What we can prove from $\ACA_0$ with the aforementioned arguments is the following: 
\begin{enumerate}
    \item If there is a (coded) $\sfbeta$-model $M$ of second-order arithmetic, and if $\alpha$ is a well-order embedding every well-order in $M$, then $D_\alpha$ is ill-founded for every recursive genedendron $(D,\varrho)$ generating $\HJ(\emptyset)$. (In happens when $\alpha$ is the ordered sum of every well-order in $M$.)
    \item If $(D,\varrho)$ is the recursive genedendron provided by the $\sfbeta$-preproof property, and if $D_\alpha$ is ill-founded for some well-order $\alpha$, then $\varrho_\alpha(B)$ is the hyperjump of $\emptyset$ for every infinite branch $B$ of $D_\alpha$.
\end{enumerate}

The reader might wonder why the resulting ordinal analysis in \autoref{Subsection: Hyperjump 0 - Sigma 1 2 alt} does not produce $\omega_2^\CK$ but $\omega_1^\CK$ since the $\Sigma^1_2$-proof-theoretic ordinal of $\ACA_0$ is $\omega_1^\CK$ while that of $\ACA_0 + `\text{$\HJ(\emptyset)$ exists'}$ is $\omega_2^\CK$.
The reason is that the fixed point $\overline{I}_\Phi$ of the theory in \autoref{Subsection: Hyperjump 0 - Sigma 1 2 alt} appears as a `proper class' and not a set, and there is no way to produce new `sets' from the fixed point $\overline{I}_\Phi$.

\section{Connecting the proof-theoretic dilator and \texorpdfstring{$\Pi^1_1[R]$}{Pi 1 1 R}-proof theory}
In this section, we prove the main theorem of this \added{chapter}. Proving the main theorem requires a variant of \emph{implicational predilator} defined \added{by Aguilera-Pakhomov} \cite[Lemma 8]{AguileraPakhomov2023Pi12} for two linear orders over $\bbN$. The implication predilator is recursive when the two linear orders are also recursive. We modify this construction to find a recursive predilator for two $R$-recursive linear orders for a $\Sigma^1_2$-singleton real $R$.

Let us recall the construction of an implicational predilator: For two linear orders $A$ and $B$ and a parameter linear order $X$ (which should be a well-order), we construct a tree $T(X)$ simultaneously building
\begin{enumerate}
    \item An infinite decreasing sequence over $B$, and
    \item An embedding $A\to X$.
\end{enumerate}
The infinite decreasing sequence tries to witness $\lnot\WO(B)$, and the embedding tries to see $\WO(A)$ by embedding $A$ to a well-order $X$. That is, the tree tries to simultaneously verify $\WO(A)\land \lnot\WO(B)$. If $\WO(A)\to\WO(B)$ holds, then $T(X)$ will be well-founded for every well-order $X$, so we can turn $T(X)$ \added{into} a well-order. In general, we can turn $T(X)$ \added{into} a linear order uniformly, which gives a predilator.

Now let $A$ and $B$ be the two $R$-recursive linear orders. We can find two recursive functions $a(X)$ and $b(X)$ such that $A=a(R)$, $B=b(R)$. $a(X)$ or $b(X)$ itself may not be a linear order for some $X$, but $a(R)$ and $b(R)$ must be linear orders.
The idea is \added{to construct} an implicational dilator and \added{the $\Sigma^1_2$-singleton} real $R$ simultaneously. We do not (and cannot) use the whole $R$ as a parameter to ensure the resulting predilator is recursive. But we can still use the partial information for $R$ to compute $a(R)$ and $b(R)$.
By \autoref{Proposition: Every Sigma 1 2 definable real is generated by a genedendron}, we may assume that $R$ is generated by a locally well-founded recursive genedendron $(D,\varrho)$.

\begin{proposition}[$\ACA_0$] \label{Proposition: Implicational predilator with a parameter real}
    Let $(D,\varrho)$ be a locally well-founded recursive genedendron generating a real $R$, and $a(Z),b(Z)$ be two recursive functions such that $a(R)$, $b(R)$ are linear orders. 
    Then we can construct a recursive predilator $\Imp_{(D,\varrho),a,b}$ such that the following holds: For every well-order $X$,
    \begin{enumerate}[label=({\roman*})]
        \item \label{Item: ImpDil-Prop1} $\Imp_{(D,\varrho),a,b}$ is a dilator iff ``if there is a well-order $\alpha$ such that $D_\alpha$ is ill-founded, then for every infinite branch $B\subseteq D_\alpha$, we have $\WO(a(R))\to \WO(b(R))$ for $R=\varrho_\alpha(B)$.''
        This statement is $\Pi^1_2$, and we write it informally as
        \begin{equation*}
            \lnot\Dil(\Lin^*(D)) \to [\WO(a(R)) \to\WO(b(R))].
        \end{equation*}
        (Recall that $\lnot\Dil(\Lin^*(D))$ iff the real generated by $(D,\varrho)$ exists.)
        \item \label{Item: ImpDil-Prop2} If $\Imp_{(D,\varrho),a,b}$ is not a dilator, $\Imp_{(D,\varrho),a,b}(X)$ is illfounded iff \added{both} $a(R)$ \added{and} $\Clim(D)$ embed into $X$.\footnote{We understand $\Clim(D)\le X$ as the assertion `$D(X)$ is ill-founded.'}
        \item \label{Item: ImpDil-Prop3} If $X\ge a(R), \Clim(D)$, then there is an embedding $b(R)\to \Imp_{(D,\varrho),a,b}(X)$.
    \end{enumerate}
\end{proposition}
\begin{proof}
    Fix a well-order $X$. The tree $T(X)$ is the set of triples $\lag \sigma,f,g\rag$ satisfying the following:
    \begin{enumerate}
        \item Every component of $\lag \sigma,f,g\rag$ is a finite sequence of the same length.
        \item $\sigma\in D_X$.
        \item \label{Item: ImplDil-Cond3} For every $i<|f|$, $b(\varrho_X(\sigma))$ does not decide%
        \footnote{$b(Z)$ is a recursive function, so some facts about $b(Z)$ are computable from a proper initial segment of $b(Z)$. $\varrho_X(\sigma)$ is a finite sequence, so `$b(\varrho_X(\sigma))$ does not decide \added{$P$}' means $\varrho_X(\sigma)$ is not long enough to decide \added{$P$}.}
        $f(i)\notin \field (b(\varrho_X(\sigma)))$.
        \item For every $i\mathrel{\added{<}}j<|f|$, $b(\varrho_X(\sigma))$ does not decide
        \begin{equation*}
            f(i),f(j)\in \field(b(\varrho_X(\sigma))) \to f(i) \ngtr_{b(\varrho_X(\sigma))} f(j).
        \end{equation*}
        \item $g\in X^{<\omega}$.
        \item For $i,j<|g|$, if $g(i)<_X g(j)$, then $a(\varrho_X(\sigma))$ does not decide
        \begin{equation*}
            i,j\in \field(a(\varrho_X(\sigma))) \land i \ge_{a(\varrho_X(\sigma))} j.
        \end{equation*}
    \end{enumerate}
    Here, `$s$ decides $P_s$' for a recursive predicate $P$ means \added{that} a partial oracle $s$ has enough information to decide the statement `$P_s$ is valid', and $P_s$ is really valid. 
    
    Informally, $\sigma\in D_X$ encodes partial information $\varrho_X(\sigma)$ for a real $R$.
    $f$ tries to construct an infinite decreasing sequence over $b(R)$, and $g$ tries to construct an embedding from $a(R)$ to $X$.
    We understand $a(\varrho_X(\sigma))$ and $b(\varrho_X(\sigma))$ as $R$-oracle Turing machines with partial information $\varrho_X(\sigma)$ for the oracle $R$, so the computation may not be finished until they get further information about the oracle.
    
    Conditions for $f$ and $g$ take the double-negated form; For example, we have the current form of \eqref{Item: ImplDil-Cond3} instead of \added{the positive form} `For every $i<|f|$, $b(\varrho_X(\sigma))$ decides $f(i)\in \field(b(\varrho_X(\sigma)))$.'
    \added{However, we take the double-negated form in place of the positive form, because} the previous `positive \added{form}' can fail not only because $f(i)\in \field(b(R))$, but also because $\varrho_X(\sigma)$ does not have enough information about $R$.
    To exclude the pathology from the latter case, we exclude $\lag \sigma,f,g\rag$ only when we clearly know $f(i)\in \field(b(R))$ fails from the given information $\varrho_X(\sigma)$.
    The real $R$ is uniquely determined, so the further branching along the tree $T(X)$ gives more information about $R$, and that $b(R)$ being $R$-recursive guarantees we can eventually determine $f(i)\in \field(b(\varrho_X(\sigma)))$ from a sufficiently large $\varrho_X(\sigma) \subseteq R$.

    Then let us impose a variant of Kleene-Brouwer order over $T(X)$, given by $\lag\sigma,f,g\rag \le_\KB \lag \sigma',f',g'\rag$ if and only if $\lag \sigma(0),f(0),g(0),\cdots,\sigma(|\sigma|-1),f(|\sigma|-1),g(|\sigma|-1)\rag$ is Kleene-Brouwer less than or equal to $\lag \sigma'(0),f'(0),g'(0),\cdots,\sigma'(|\sigma'|-1),f'(|\sigma'|-1),g'(|\sigma'|-1)\rag$.
    $\sigma(i)$ and $\sigma'(i)$ are compared under the underlying order of $D$, and all other components are compared under the order over $\bbN$. Let us define $\Imp_{(D,\varrho),a,b}(X)$ by the resulting linearization. Then we can see that $\Imp_{(D,\varrho), a,b}$ admits the support transformation, so $\Imp_{(D,\varrho), a,b}$ is a semidilator, and it is easy to see that $\Imp_{(D,\varrho),a,b}$ is a predilator.

    Now let us check the equivalent condition for $\Imp_{(D,\varrho),a,b}$ being a dilator: 
    If $\Lin^*(D)$ is a dilator, then $D_X$ is well-founded for every well-order $X$, so $\Imp_{(D,\varrho),a,b}(X)$ cannot have an infinite branch for every well-order $X$. Hence $\Imp_{(D,\varrho),a,b}$ is a dilator in this case.
    Now suppose that $\Lin^*(D)$ is not a dilator, and divide the cases:
    \begin{enumerate}
        \item Suppose that $b(R)$ is a well-order. Then $\Imp_{(D,\varrho),a,b}(X)$ as a tree cannot have an infinite branch for a well-ordered $X$, so $\Imp_{(D,\varrho),a,b}$ is a dilator.
        \item Suppose that $b(R)$ is ill-founded. Then $\Imp_{(D,\varrho),a,b}(X)$ as a tree has an infinite branch if and only if $D_X$ has an infinite branch and there is an embedding $a(R)\to X$. Put it in other words,
        \begin{equation} \label{Formula: When D(X) is ill founded}
            \lnot \WO(\Imp_{(D,\varrho),a,b}(X)) \iff \lnot \WO(\Lin^*(D_X)) \land a(R)\le X.
        \end{equation}
        Hence $\Imp_{(D,\varrho),a,b}$ is a dilator iff
        \begin{equation*}
            \forall^1 X [\WO(X)\to \lnot[a(R)\le X\land \lnot\WO(\Lin^*(D_X))]],
        \end{equation*}
        which is equivalent to
        \begin{equation} \label{Formula: Should be equivalent to lnot WO a(R)}
            \forall^1 X [[\WO(X) \land \lnot\WO(\Lin^*(D_X))]  \to a(R)\nleq X].
        \end{equation}
        We claim that \eqref{Formula: Should be equivalent to lnot WO a(R)} is equivalent to $\lnot \WO(a(R))$: Fix a well-order $\alpha$ such that $\Lin^*(D_\alpha)$ is ill-founded. If $a(R)$ is well-founded, then $x=\alpha + a(R)$ witnesses the negation of \eqref{Formula: Should be equivalent to lnot WO a(R)}.
        Conversely, if the negation of \eqref{Formula: Should be equivalent to lnot WO a(R)} holds, then $a(R)$ is embedded into some well-order, so $a(R)$ is a well-order too.
    \end{enumerate}
    It finishes the proof of \ref{Item: ImpDil-Prop1}. 
    For \ref{Item: ImpDil-Prop2}, observe that if $\Imp_{(D,\varrho),a,b}$ is not a dilator, then $D_X$ has an infinite branch, $b(R)$ is ill-founded, and \eqref{Formula: When D(X) is ill founded} holds.
    We can see that \eqref{Formula: When D(X) is ill founded} is precisely \ref{Item: ImpDil-Prop2}.

    For \ref{Item: ImpDil-Prop3}, we need to construct an embedding $e\colon b(R)\to \Imp_{(D,\varrho),a,b}(X)$. Let us fix an infinite branch $\lag\sigma_i\mid i\in\bbN\rag$ of $D_X$. Then $\bigcup_{n\in\bbN} \varrho_X(\lag\sigma_i\mid i<n\rag)=R$. Also, we can fix an embedding $g\colon a(R)\to X$ and extend $g$ to a function of domain $\bbN$ in a silly manner (for example, by taking $g(i)$ to be the $a(R)$-least element when $i\notin \field(a(R))$.)

    Let us enumerate $b(R) = \{\beta_i\mid i\in \bbN\}$ in an $\bbN$-increasing way. Then define $e\colon b(R)\to \Imp_{(D,\varrho),a,b}(X)$ inductively as follows: Assume that $e(\beta_j)$ is defined for every $j<i$.
    \begin{enumerate}
        \item Suppose that $\beta_i = \max_{b(R)}\{\beta_j\mid j\le i\}$. Then define
        \begin{equation*}
            e(\beta_i) = \lag \lag\sigma_0\rag, \lag \beta_i\rag, g\restricts 1\rag.
        \end{equation*}
        \item Otherwise, let $\beta_l = \min_{b(R)}\{\beta_j>_{b(R)} \beta_i\mid j<i\}$, and suppose
        \begin{equation*}
            e(\beta_l) = \lag \sigma\restricts m, f, g\restricts m\rag.
        \end{equation*}
        Then define
        \begin{equation*}
            e(\beta_i) = \lag \sigma\restricts (m+1), f^\frown\lag\beta_i\rag, g\restricts (m+1)\rag.
        \end{equation*}
    \end{enumerate}
    We can see that $e$ is an order-preserving map as follows: 
    Observe that the order between $e(\beta_i)$ and $e(\beta_j)$ only depends on the second component. 
    Furthermore, the second component of $e(\beta_i)$ is an $b(R)$-decreasing sequence ending with $\beta_i$ indexed by natural numbers less than or equal to $i$.
    
    Now we claim by induction on $i$ that $e\restricts\{\beta_0,\cdots,\beta_i\}$ is order-preserving:
    If $\beta_i = \max_{b(R)}\{\beta_j\mid j\le i\}$, then clearly $\beta_j\le\beta_i$ and $e(\beta_j)\le e(\beta_i)$ for $j<i$ by definition of $e(\beta_i)$.
    Otherwise, let $\beta_l = \min_{b(R)}\{\beta_j>_{b(R)} \beta_i\mid j<i\}$.
    Then we have $\beta_i<\beta_l$ and $e(\beta_i)<e(\beta_l)$.
    For $j<i$, if $\beta_i<\beta_j$, then $\beta_l\le \beta_j$, so $e(\beta_i)<e(\beta_l)\le e(\beta_j)$ by the induction hypothesis.
    If $\beta_i\ge\beta_j$, then $\beta_l>\beta_j$, so $e(\beta_l)>e(\beta_j)$.
    This means that the following holds: If $f_l$ and $f_j$ are the second component of $e(\beta_l)$ and $e(\beta_j)$ respectively, then either
    \begin{enumerate}
        \item there is a least $m<|f_l|,|f_j|$ such that $f_l(m)>f_j(m)$, or
        \item $f_j$ is a proper extension of $f_l$.
    \end{enumerate}
    If $f_i$ is the second component of $e(\beta_i)$, then $f_i$ is a proper extension of $f_l$. Hence if the first case holds, then $e(\beta_i)>e(\beta_j)$ as $m$ is the least number such that $f_l(m)=f_i(m)>f_j(m)$.
    In the second case, suppose that $f_j = f_l^\frown \lag \beta_k,\cdots,\beta_j\rag$ for some $k<j$. Then $\beta_k>\beta_i$ never holds; Otherwise, it contradicts the choice of $\beta_l$.
    Hence $\beta_k\le \beta_i$, so $f_j\le f_i$ under the lexicographic order.
\end{proof}

\begin{theorem} \label{Theorem: Main theorem}
    Let $T$ be a $\Pi^1_2$-sound theory extending $\ACA_0$ and $(D,\varrho)$ be a recursive locally well-founded genedendron generating $R$. Furthermore, assume that $T$ proves $(D,\varrho)$ is a locally well-founded genedendron.
    If $\alpha$ is an $R$-recursive linear order such that $D_\alpha$ is ill-founded, then 
    \begin{equation*}
        |T|_{\Pi^1_2}(\alpha) = \lvert T[R] + \WO(\alpha)\rvert_{\Pi^1_1[R]}.
    \end{equation*}
\end{theorem}
\begin{proof}
    Since $\alpha$ is $R$-recursive, we can find a recursive function $a$ such that $\alpha=a(R)$. For one direction of the inequality, let us decompose $|T|_{\Pi^1_2}$ by the sum of recursive dilators $\sum_{n\in\bbN} D_n$. Since $T\vdash \Dil(D_i)$ for each $i$, we have
    \begin{equation*}
        T[R] + \WO(a(R)) \vdash \WO((D_0+\cdots D_n)(a(R))).
    \end{equation*}
    Hence $(D_0+\cdots+ D_n)(a(R)) \le \lvert T[R] + \WO(a(R))\rvert_{\Pi^1_1[R]}$, so taking the supremum to $n$ gives
    \begin{equation*}
        |T|_{\Pi^1_2}(a(R)) \le \lvert T[R] + \WO(a(R))\rvert_{\Pi^1_1[R]}.
    \end{equation*}
    For the remaining direction, suppose that $b$ is a recursive function such that $b(R)$ is a well-order and
    \begin{equation*}
        T[R] + \WO(a(R)) \vdash \WO(b(R)).
    \end{equation*}
    Hence we have
    \begin{equation*}
        T \vdash \lnot \Dil(D) \to (\WO(a(R))\to \WO(b(R))).
    \end{equation*}
    Since $T$ extends $\ACA_0$, $T$ proves \autoref{Proposition: Implicational predilator with a parameter real}. Also, we assumed that $T$ proves $(D,\varrho)$ is a locally well-founded genedendron. 
    Hence we have $T \vdash \Dil(\Imp_{(D,\varrho),a,b})$, so $\Imp_{(D,\varrho),a,b}\le |T|_{\Pi^1_2}$. Since $D_{a(R)}$ is ill-founded, we have $a(R)\ge\Clim(D)$, so by \ref{Item: ImpDil-Prop3} of \autoref{Proposition: Implicational predilator with a parameter real}, we have $b(R) \le \Imp_{(D,\varrho),a,b}(a(R))$. In sum, we have
    \begin{equation*}
        b(R)\le \Imp_{(D,\varrho),a,b}(a(R)) \le |T|_{\Pi^1_2}(a(R)).
    \end{equation*}
    Since every $T[R]$-provably $R$-recursive well-order takes the form $b(R)$ for some recursive $b$, we have
    \begin{equation*}
        \lvert T[R]+\WO(a(R))\rvert_{\Pi^1_1[R]}\le |T|_{\Pi^1_2}(a(R)). \qedhere 
    \end{equation*}
\end{proof}

The following corollary follows by applying \autoref{Theorem: Main theorem} to the genedendron we found in \autoref{Subsection: Hyperjump 0 - Sigma 1 2 alt}:
\begin{corollary} \pushQED{\qed}
    \begin{equation*}
        \lvert\ACA_0 + \text{$\HJ(\emptyset)$ exists}\rvert_{\Pi^1_1[\HJ(\emptyset)]} = \lvert \ACA_0 \rvert_{\Pi^1_2}(\omega_1^\CK) = \varepsilon_{\omega_1^\CK+1}. \qedhere 
    \end{equation*}
\end{corollary}

\section{Final remarks}
The main result in this paper implies there is a systematic way to produce proof-theoretic information about $\Pi^1_1[R]$-consequences for a $\Sigma^1_2$-\added{singleton} real $R$ by combining $\Pi^1_2$- and $\Sigma^1_2$-information of a theory. In this last section, we discuss what we can observe from the results in this paper and directions for future work.

We may view $|T|_{\Pi^1_2}(\alpha)$ as a `section' of $|T|_{\Pi^1_2}$, and we may ask if sections of $|T|_{\Pi^1_2}$ gives useful information about $|T|_{\Pi^1_2}$. For example, we may ask if the sections of a proof-theoretic dilator are enough to determine the $\Pi^1_2$-consequence comparison modulo true $\Sigma^1_2$-statements:
\begin{question}
    Let $S$, $T$ be $\Pi^1_2$-sound r.e. extensions of $\ACA_0$. Define
    \begin{equation*}
        S \subseteq^{\Sigma^1_2}_{\Pi^1_2} T \iff \text{For every $\Pi^1_2$-sentence $\sigma$,}\ S\vdash^{\Sigma^1_2}\sigma\implies T\vdash^{\Sigma^1_2}\sigma,
    \end{equation*}
    where $T\vdash^{\Sigma^1_2}\sigma$ means `$\sigma$ is proavable from $T$ with the set of true $\Sigma^1_2$-sentences.' Do we have the following equivalence?  
    \begin{equation*}
        S \subseteq^{\Sigma^1_2}_{\Pi^1_2} T \iff \forall (D,\varrho) [|S|_{\Pi^1_2}(\Clim (D)) \le |T|_{\Pi^1_2}(\Clim (D))],
    \end{equation*}
    where $(D,\varrho)$ ranges over $\ACA_0$-provably recursive locally well-founded genedendrons.
\end{question}

The example illustrated in \autoref{Subsection: Hyperjump 0 - Sigma 1 2 alt} can be understood as a $\Sigma^1_2$-ordinal analysis for Peano arithmetic plus the unary predicate for the hyperjump of $\emptyset$ (equivalently, Kleene's $\calO$.) We obtained a genedendron by taking the $\sfbeta$-preproof of the false sequent obtained from the $\sfbeta$-preproof property. It resembles obtaining the proof-theoretic ordinal or dilator from cut elimination, so we may understand the role of the $\sfbeta$-preproof property in $\Sigma^1_2$-proof theory is like that of the cut-elimination in $\Pi^1_1$- and $\Pi^1_2$-proof theory.

We may expect the same for the $\omega$-logic: For a given $\omega$-proof system for some theory $T$ with the preproof property, we can generate an $\omega$-model of $T$ from an infinite branch of a preproof of the empty sequent. It should provide the framework for the $\Sigma^1_1$-proof theory.
Then we may ask which algebraic object captures the behavior of the ill-founded $\omega$-preproof: In the $\Pi^1_1$- and $\Pi^1_2$- cases, each of them are proof-theoretic ordinal and proof-theoretic dilator respectively. In the $\Sigma^1_2$-case, it is the climax of a $\sfbeta$-preproof of the empty sequent. In the $\Sigma^1_1$-case, pseudo-well-orders might behave as the algebraic object as foreshadowed in \cite[\S2.4]{TowsnerWalsh2024Classification}.

One may ask if we can generalize the main theorem. One natural direction is pursuing a proof-theroetic meaning of $|T|_{\Pi^1_2}(\alpha)$ for a larger $\alpha$. There are two main obstacles to this way of generalization. One is related to the limit of genedendrons and pseudodilators. We associated some ordinals below $\delta^1_2$ with a recursive genedendron $(D,\varrho)$, which generates a $\Sigma^1_2$-singleton real and the climax of $D$ is the ordinal.
To extract a proof-theoretic meaning of $|T|_{\Pi^1_2}(\alpha)$ for larger $\alpha$, we should associate $\alpha$ to a more sophisticated object than genedendrons. One of the reasons is that the underlying pseudodilator of a recursive genedendron must have the climax less than $\delta^1_2$, so any recursive genedendron cannot be associated with ordinals larger than $\delta^1_2$.
The appropriate object should come from $\sfbeta_n$-logic or its generalizations.
The other obstacle, which is more serious than the previous one, is that the proof of \autoref{Theorem: Main theorem} relies on the fact that every $\Pi^1_1[R]$-statement for a $\Sigma^1_2$-singleton real $R$ is $\Pi^1_2$, so the proof-theoretic dilator already sees $\Pi^1_1[R]$-statements. For a more complex $R$, there is no guarantee that the proof-theoretic dilator sees $\Pi^1_1[R]$-statements.
We might need to consider recursive genedendrons `derived from' a higher object associated with $\alpha>\delta^1_2$ simultaneously to get the right generalization.

The other way of a generalization is asking the proof-theoretic meaning of $|T|_{\Pi^1_n}(P)$ for an $(n-1)$-ptyx $P$ with reasonable definability constraint (like, implicit $\Sigma^1_n$-definability.) As an example, the author conjectures the following:
\begin{conjecture}
    Working over $\ZFC$ with the existence of $0^\sharp$, let $T$ be a $\Pi^1_3$-sound extension of $\ACA_0$. If we view $0^\sharp$ as a dilator (cf. \cite{AguileraFreundRathjenWeiermann2022}), we have
    \begin{equation*}
        |T|_{\Pi^1_3}(0^\sharp) = |T + \exists 0^\sharp|_{\Pi^1_1[0^\sharp]}.
    \end{equation*}
\end{conjecture}

$\Sigma^1_2$-altitude appears as a byproduct of the main result, but it should be able to get its own attention. $\Sigma^1_2$-altitude defines a way to rank $\Sigma^1_2$-singleton reals and is computable by hand for some special cases. We may ask if the $\Sigma^1_2$-altitude has recursion-theoretic or set-theoretic characterization, and the following conjecture suggests one possibility:
\begin{conjecture}
    Let $R$ be a $\Sigma^1_2$-singleton real. $\Alt_{\Sigma^1_2}(R)$ is equal to the least height of a transitive model $M$ of $\ATR_0^\mathsf{set}$ on which $R$ is a $\Sigma_1$-definable class in the language of set theory.
\end{conjecture}

\appendix
\section{Iterated hyperjumps} \label{Section: Iterated Hyperjumps}

The main aim of this section is to provide the $\Pi^1_1$-singleton definition of the $\xi$-th iterated hyperjump of the empty set. We formulate the iterated hyperjump of the empty set as an ordinal notation system.

\subsection{Inductive definition}
\begin{definition}
    Let $\Gamma\colon \calP(\omega)\to\calP(\omega)$ be an operator (which is usually definable) that is \emph{inclusive}, i.e., $X\subseteq \Gamma(X)$ for every set $X\subseteq\omega$.
    We define $\lag \Gamma^\xi \mid \xi\in\Ord\rag$ recursively by $\Gamma^\xi = \bigcup_{\eta<\xi}\Gamma(\Gamma^\eta)$.
    The \emph{closure point of $\Gamma$}, which is denoted by $|\Gamma|$, is the least ordinal $\lambda$ such that $\Gamma^\lambda = \Gamma^{\lambda+1}$. We say $\Gamma^\lambda$ for the closure point of $\Gamma$ the \emph{set defined by $\Gamma$}.
\end{definition}

We formulated the definition of $\Gamma^\xi$ in set theory. However, we can formulate the same definition in a subsystem of second-order arithmetic by replacing ordinals with elements of a sufficiently long well-order.
We do not know in general that we have a sufficiently long well-order terminating the inductive definition. However, we can still state which set is defined by a given inclusive operator.

We want to make sure the set defined by an operator is $\ACA_0$-provably unique even when it may not $\ACA_0$-provably exist. If we work over $\ATR_0$, the trichotomy of well-orders guarantees the sequence $\lag \Gamma^\xi \mid \xi \in \field(\alpha)\rag$ given by an arithmetical inclusive operator $\Gamma$ only depends on the ordertype of $\alpha$. The next proposition says we can weaken $\ATR_0$ to $\ACA_0$ for an \added{inclusive arithmetical operator:}

\begin{proposition}[$\ACA_0$] \label{Proposition: Uniqueness of the Gamma defined set over ACA0}
    Let $\Gamma$ be an inclusive arithmetical operator. For a given well-order $\alpha$, let us say a set $X$ is a \emph{chain of $\Gamma$-inductive definition along $\alpha$} if for each $\xi\in \field(\alpha)$ and $n$,
    \begin{equation*}
        n\in (X)_\xi \iff \exists \eta<_\alpha \xi [n\in \Gamma((X)_\eta)],
    \end{equation*}
    where $(X)_\xi = \{m\mid \lag \xi,m\rag\in X\}$.
    A chain of $\Gamma$-inductive definition $X$ along $\alpha$ is \emph{exact} if $\alpha$ has the maximal element $\mu$ and its predecessor $\mu^-$, and 
    \begin{enumerate}
        \item $\eta<_\alpha\xi <_\alpha\mu$ implies $(X)_\eta \subsetneq (X)_\xi$.
        \item $(X)_\mu = (X)_{\mu^-}$.
    \end{enumerate}

    Suppose that $\alpha$ and $\beta$ are well-orders, $X$ an exact chain of $\Gamma$-inductive definition along $\alpha$, $Y$ an exact chain of $\Gamma$-inductive definition along $\beta$. Then there is an isomorphism $f\colon \alpha\to\beta$ such that $(X)_\xi = (Y)_{f(\xi)}$ for each $\xi\in \field(\alpha)$.
\end{proposition}
\begin{proof}
    Let $\mu_\alpha$, $\mu^-_\alpha$ be the maximal element and its predecessor of $\alpha$, and similarly for $\mu_\beta$ and $\mu^-_\beta$.
    We prove the following assertion by arithmetical transfinite induction (which is a theorem of $\ACA_0$; See \cite[Lemma V.2.1]{Simpson2009}):
    For each $i\in\field(\alpha)$, either $i=\mu_\alpha$ or there is a unique $j\in\field(\beta)$ such that $j\neq \mu_\beta$, $(X)_i = (Y)_j$.

    Suppose that $i\neq \mu_\alpha$ and the claim holds for $i'<_\alpha i$, so for each $i'<i$ we have a unique $j'\in \field(\beta)\setminus \{\mu_\beta\}$ such that $(X)_{i'} = (Y)_{j'}$.
    Such $j'$ is uniformly arithmetically definable from $i'$, $X$, $Y$, so we can unambiguously denote such $j'$ by $f(i')$ for some arithmetically definable function $f$.
    Let us define
    \begin{equation*}
        j = \min{}_\beta\{k\in\field(\beta)\mid \forall i'<_\alpha i [f(i') < k]\}.
    \end{equation*}
    We claim $(X)_i = (Y)_j$: First,
    \begin{equation*}
        n\in (X)_i \iff \exists i'<_\alpha i [n\in \Gamma((X)_{i'})] \implies \exists j'<_\beta j [n\in (Y)_{j'}]
    \end{equation*}
    so $(X)_i \subseteq (Y)_j$. Conversely, if $j'<_\beta j$, then there is $i'<_\alpha i$ such that $j'\le_\beta f(i')$, so we have the reversed inclusion. The construction of the isomorphism $f\colon\alpha\to\beta$ is straightforward.
\end{proof}
Hence the following $\Sigma^1_2$-formula witnesses the set $H$ defined by an arithmetical $\Gamma$ is a $\Sigma^1_2$-singleton set: ``There is a well-order $\alpha$ with the maximal element $\mu$ with its predecessor, and an exact chain of $\Gamma$-inductive definition, we have $H=(X)_\mu$.'' Moreover, $\ACA_0$ proves the uniqueness of the real satisfying the $\Sigma^1_2$-definition. (Note that $\ACA_0$ does not prove the existence of the real in general.)

\begin{remark} \label{Remark: Pi 1 1 singleton definition for an arithmetical ind def set}
    We \added{could} make the $\Sigma^1_2$-\added{singleton} definition into a $\Pi^1_1$-\added{singleton} definition by using linear orders with a $\Pi^0_m$-transfinite induction for a sufficiently large $m$ in place of well-orders.
    More precisely, if $\Gamma$ is $\Pi^0_m$-definable, the proof of \autoref{Proposition: Uniqueness of the Gamma defined set over ACA0} shows the following $\Pi^1_1$-formula defines the set $H$ defined by $\Gamma$:
    \begin{quote}
        Suppose that $\alpha$ is a linear order with the maximal element $\mu$ with its predecessor, such that $\Pi^0_{m+99}$-transfinite induction along $\alpha$ holds. For every exact chain of $\Gamma$-inductive definition, we have $H=(X)_\mu$.
    \end{quote}
    $\ACA_0$ proves there is at most one real satisfying the $\Pi^1_1$-formula.
\end{remark}

\subsection{Ordinal notation system}

In this section, we follow notions in \cite[\S8]{AczelRichter1974IDReflAdm}.
First, let us fix the following notations:
\begin{itemize}
    \item $\sfsup(a,b) = \lag 0,a,b\rag$.
    \item $a \ovee b = \lag 1,a,b\rag$.
\end{itemize}
Also, for a natural number $a$, $x$ and a real number $R$, $[a](\cdot,R)$ denotes the $a$th $R$-primitive recursive function. 

As an illustrative example for an ordinal notation system, the following ordinal notation resembles Kleene's $\calO$, where $|a| = \min\{\alpha\mid a\in M_{\alpha+1}\}$:
\begin{enumerate}
    \item $M_0=\varnothing$.
    \item $M_{\alpha+1} = M_\alpha \cup \{0\} \cup \{ \sfsup(a,b) \mid a\in M_\alpha ,\ \added{\forall^0 x}  [ [b](x,M_{|a|})\in M_\alpha]\}$.
    \item $M_\delta = \bigcup_{\alpha<\delta} M_\alpha$ if $\delta$ is limit.
\end{enumerate}
Here we understand $\sfsup(a,b)$ as the notation for the ordinal 
\begin{equation*}
    \max\big\{|a|+1, \sup_{x<\omega} \big(|[b](x,M_{|a|})|+1\big) \big\}.
\end{equation*}
Now let us use the notations for $X\subseteq\omega$:
\begin{itemize}
    \item $\calF(X) = \{x \mid \lag x,x\rag\in X\}$, and
    \item $X_{<x} = \{y\mid \lag y,x\rag\in X\land \lag x,y\rag\notin X\}$. 
    \item $M^*_\alpha = \{\lag x,y\rag \mid x,y\in M_\alpha\land |x|\le |y|\}$.
\end{itemize}
If we define 
\begin{equation*}
    \Theta(X) = \{0\}\cup \{\sfsup(a,b)\mid a\in \calF(X)\land \added{\forall^0 x} [ [b](x,X_{<a})\in M_\alpha]\},
\end{equation*}
then the previous ordinal notation system satisfies $M_{\alpha+1} = M_\alpha\cup\Theta(M_\alpha^*)$. Let us take it as a general definition for an ordinal notation system:
\begin{definition}
    Let $\Theta$ be an inclusive operator. $\frakM^\Theta = (M^\Theta,|\cdot|)$ is defined by
    \begin{enumerate}
        \item $M_0=\varnothing$.
        \item $M_{\alpha+1} = M_\alpha\cup\Theta(M_\alpha^*)$.
        \item $M_\delta = \bigcup_{\alpha<\delta} M_\alpha$ if $\delta$ is limit.
        \item $M^\Theta = \bigcup_{\alpha\in\Ord}M_\alpha$.
    \end{enumerate}
    $|M|$ is the least ordinal $\alpha$ such that $M_\alpha=M_{\alpha+1}$.
\end{definition}
The definition of $M_{\alpha+1}$ depends not only $M_\alpha$, but also on $M_\alpha^*$ for $a\in M_\alpha$, so the definition of $M_\alpha$ itself is not precisely the inductive definition.
We get around this issue by defining $M_\alpha^*$ inductively:

\begin{lemma}[{\cite[Lemma 8.13]{AczelRichter1974IDReflAdm}}]
    For an operator $\Theta$, if we define
    \begin{equation*}
        \Theta_\le(X) = \{\lag x,y\rag \mid x\in \calF(X)\land y\in \Theta(X)\setminus\calF(X)\} \cup \{\lag x,y\rag \mid x,y\in \Theta(X)\setminus\calF(X)\},
    \end{equation*}
    then $\Theta_\le^\alpha = M_\alpha^*$ for every ordinal $\alpha$.
\end{lemma}
Note that if $\Theta$ is $\Pi^0_m$ or $\Sigma^0_m$, then so is $\Theta_\le$ since $\calF$ is primitive recursive. 
The following type of ordinal notation systems are guaranteed to have a nice behavior, which the coding lemma will characterize:
\begin{definition}
    A notation system $\frakM=(M,|\cdot|)$ is \emph{Richterian}\footnote{\cite{AczelRichter1974IDReflAdm} used the word \emph{good}, but the word `good' is not good to use since it is used too much elsewhere.} if $\frakM=\frakM^\Theta$ where $\Theta(X) = \frakR(X) \cup \Phi(X)$, where $\frakR(X)$ is the \emph{Richter operator}
    \begin{equation*}
        \frakR(X) = \{0\}\cup \{\sfsup(a,b)\mid a\in \calF(X)\land \added{\forall^0 x} [[b](x,X_{<a})\in\calF(X)]\} \cup \{a\ovee b\mid a\in \calF(X)\lor b\in\calF(X)\},
    \end{equation*}
    and $\Phi(X)$ is always disjoint from $\{0\}\cup \{\sfsup(a,b),\ a\ovee b\mid a,b\in\omega\}$.

    For a Richterian notation system $\frakM$, an ordinal $\lambda\le |M|$ is \emph{$\frakM$-Richterian} if $\frakR(M^*_\lambda)\subseteq M_\lambda$. (In some sense, an ordinal is $\frakM$-Richterian if $M_\lambda$ is closed under the Richter operator.) $|M|$ is $\frakM$-Richterian, but in many cases, there are $\frakM$-Richterian ordinals below $|M|$.
\end{definition}

The role of $a\ovee b$ is rather technical, and the author does not know about its necessity. Richter also mentioned in \cite{Richter1971Mahlo} that he does not know if the $\ovee$ operator is necessary. Also, the hierarchy given by a Richerian operator increases \added{under} Turing reducibility in the following sense:
\begin{lemma}[{\cite[Lemma A.3]{AczelRichter1974IDReflAdm}}] \label{Lemma: Richterian Hierarchy Turing reducibility}
    If $|a|+1<\beta$, then $\mathsf{TJ}(M_{|a|})\le_\sfT M_\beta$ uniformly in $a$. Moreover, the reduction can be primitive recursive.
\end{lemma}
\begin{proof}
    Let $e$ be a recursive function such that
    \begin{equation*}
        [e(a,x)](t,M_{|a|})= 
        \begin{cases}
            a & \text{if }\lnot T^{M_{|a|}}(x,x,t), \\
            1 & \text{otherwise.}
        \end{cases}
    \end{equation*}
    (Note that $1\notin M$.) Then we can see that 
    \begin{equation*}
        x\notin \mathsf{TJ}(M_{|a|}) \iff \forall t \lnot  T^{M_{|a|}}(x,x,t) \iff 
        \forall t [e(a,x)](t,M_{|a|}) = a\in M_{|a|+1} \iff \sfsup(a,e(a,x))\in M_\alpha,
    \end{equation*}
    so $\bbN\setminus \mathsf{TJ}(M_{|a|}) \le_\sfT M_\alpha$ is witnessed by $x\mapsto \sfsup(a,e(a,x))$. Since Kleene's $T$-predicate is primitive recursive, we can choose $e$ to be primitive recursive.
\end{proof}

The next proposition says $M_{|x|}$ and $M^*_{|x|}$ are Turing equivalent:
\begin{proposition} \label{Proposition: M hierarchy and M star hierarchy are Turing equivalent}
    Let $\frakM=(M,|\cdot|)$ be a Richterian notation system.
    For each $x\in M$, $M_{|x|} \equiv_\sfT M^*_{|x|}$.
\end{proposition}
\begin{proof}
    $M_{|x|} \le_\sfT M^*_{|x|}$ is clear. (In fact, $M_{|x|} \le_1 M^*_{|x|}$.) For $M^*_{|x|}\le_\sfT M_x$, suppose that $\lag u,v\rag\in M^*_{|x|}$ so $u,v\in M_{|x|}$.
    Let $\id$ and $c_y$ be uniform primitive recursive codes for the identity function and the constant function of value $y$, respectively.
    First, let us observe that if $|u|<|x|$, then
    \begin{equation*}
        |u| > |v| \iff \sfsup(u,c_v)\in M_{|x|}.
    \end{equation*}
    Also, $|u|<|x|$ iff $\sfsup(u,\id)\in M_{|x|}$ since $|\sfsup(u,\id)|=|u|+1$.
    Hence for \emph{every} $u$ and $v$,
    \begin{equation*}
        \lag u,v\rag\in M^*_{|x|} \iff u,v\in M_{|x|} \land [\sup(u,\id)\in M_{|x|}\to \sup(u,c_v)\in M_{|x|}]. \qedhere 
    \end{equation*}
\end{proof}

The following \emph{coding lemma} (which should be called the coding theorem) first appeared in \cite{Richter1971Mahlo}, and its full \added{details} appeared in \cite[\S A]{AczelRichter1974IDReflAdm}. According to \cite{Richter1971Mahlo}, its initial idea goes back to Gandy's unpublished \added{proof} of $|\Pi^0_1| = \omega_1^\CK$.
\begin{theorem}[Richter's Coding Lemma, {\cite[Lemma 8.17]{AczelRichter1974IDReflAdm}}]
    Let $\frakM=(M,|\cdot|)$ be a Richterian notation system and let $T_\frakM = \{\lag x,\alpha\rag\mid \alpha\in\Ord\land x\in M_\alpha\}$. Then
    \begin{enumerate}
        \item Every $\frakM$-Richterian ordinal is admissible relative to $T_\frakM$.
        \item For every $\Sigma_1$-formula $\phi(v_0,\cdots,v_{n-1})$ in the language of primitive recursive set functions on ordinals with the extra predicate for $T_\frakM$, there is a primitive recursive function $h$ such that for every $\frakM$-Richterian ordinal $\lambda$,
        \begin{equation*}
            a_0,\cdots a_{m-1}\in M_\lambda\land \lambda\vDash \phi(|a_0|,\cdots,|a_{m-1}|) \iff h(a_0,\cdots, a_{m-1})\in M_\lambda.
        \end{equation*}
        \item If $\lambda$ is $\frakM$-Richterian, then for $X\subseteq \omega$, $X$ is $\lambda$-r.e. in $T_\frakM\restriction \lambda$ iff $X\le_\sfm M_\lambda$.
    \end{enumerate}
\end{theorem}
We will only use the first clause of Richter's Coding Lemma in this paper.

\subsection{Ordinal notation system for the least recursively inaccessible ordinal}
Let us consider the following Richterian notation system:

\begin{definition}
    Let us define $J$ as follows:
    \begin{equation*}
        n\in J(X) \iff n\in \frakR(X) \lor [\frakR(X)\subseteq \field(X) \land n\in \{\lag 2,x\rag \mid x\in \calF(X)\}].
    \end{equation*}
\end{definition}
Roughly, $x\mapsto \lag 2,x\rag$ maps $x$ to an ordinal notation denoting the next admissible ordinal greater than $|x|$. We can see that $J$ is an arithmetical Richterain operator, so the closure point of $\frakM^J=(M,|\cdot|)$ is admissible. 
From the informal description, we can suspect that $|\lag 2,x\rag|$ for $x\in M$ is admissible, and the following proposition confirms it is:
\begin{proposition} \label{Proposition: M-Richterian ordinal for J-operator}
    For $x\in M$, $|\lag 2,x\rag|$ is $\frakM$-Richterian. In fact, $|\lag 2,x\rag|$ is the least $\frakM$-Richterian ordinal above $|x|$.
\end{proposition}
\begin{proof}
    $\lambda=|\lag 2,x\rag|$ means $\lambda$ is the least ordinal such that $\lag 2,x\rag\in M_{\lambda+1} = M_\lambda \cup J(M^*_\lambda)$, which happens only when $\frakR(M^*_\lambda) \subseteq M_\lambda$. Hence $\lambda$ is $\frakM$-Richterian.
    If $\lambda'>|x|$ is another $\frakM$-Richterian ordinal, then $\lag 2,x\rag \in M_{\lambda'+1}$. It shows $|\lag 2,x\rag|\le \lambda'$.
\end{proof}

By Richter's coding lemma, $|M|$ is not only admissible but also limit admissible. Hence $|M|$ is recursively inaccessible.
Now let us claim that $|M|$ is the least recursively inaccessible ordinal. Observe that $J(X)$ works as follows: If $X$ is not $\frakR$-saturated, take $\frakR(X)$. If $X$ is $\frakR$-saturated, take $\{\lag 2,x\rag \mid x\in \calF(X)\}$.
Also, observe that $\frakR$ is $\Pi^0_1$, and $X\mapsto \{\lag 2,x\rag\mid x\in \calF(X)\}$ is primitive recursive (equivalently, $\Pi^0_0$.)
Let us generalize the definability complexity of $J$ (and that of $J_\le$) as follows:
\begin{definition}
    A \emph{$[\Pi^0_1,\Pi^0_0]$-operator} is an operator $\Gamma$ of the form
    \begin{equation*}
        n\in \Gamma(X)\iff n\in \Gamma_0(X) \lor [\Gamma_0(X)\subseteq X \implies n\in \Gamma_1(X)]
    \end{equation*}
    for a $\Pi^0_1$-operator $\Gamma_0$ and a primitive recursive operator $\Gamma_1$.
\end{definition}

Then we have the following:
\begin{proposition} \label{Proposition: The closure point of a Pi01-Pi00 operator}
    If $\Gamma$ is a $[\Pi^0_1,\Pi^0_0]$-operator. then $|\Gamma|$ is less than or equal to the least recursively inaccessible ordinal.
\end{proposition}

The proof of \autoref{Proposition: The closure point of a Pi01-Pi00 operator} depends on the proof of the following lemma appearing in \cite[Lemma 2.2]{Richter1971Mahlo}. We include its proof for completeness.

\begin{lemma} \label{Lemma: Admissible ordinal is a closure pt of a Pi01 operator}
    Let $\Gamma$ be a $\Pi^0_1$-operator and $\lambda$ an admissible ordinal. Then $\Gamma^\lambda = \Gamma^{\lambda+1}$.
\end{lemma}
\begin{proof}
    First, since `$n\notin \Gamma(X)$' is an r.e. predicate on $\bbN\times\calP(\bbN)$, we can find a recursive predicate $R$ on $\bbN^3$ such that 
    \begin{equation*}
        n\notin \Gamma(X) \iff \exists  x\in \Seq(X) \exists y\in \Seq(\bbN\setminus X) R(n,x,y).
    \end{equation*}
    Here $\Seq(X)$ is the set of natural numbers coding a finite sequence over $X$.
    In particular,
    \begin{equation*}
        n\in \Gamma(X) \iff \forall x\in \Seq(X) \forall y [R(n,x,y)\to y\in \bbN\setminus \Seq(\bbN\setminus X)].
    \end{equation*}
    Note that both $X\mapsto \Seq(X)$ and $X\mapsto \bbN\setminus \Seq(\bbN\setminus X)$ are monotone. Now let $\lag\Gamma^\xi\mid \xi\in\Ord\rag$ be a sequence defined by $\Gamma$. We claim that if $\lambda$ is admissible and $n\in \Gamma^{\lambda+1}=\Gamma(\Gamma^\lambda)$, then $n\in \Gamma^\lambda$.

    Suppose that $n\in \Gamma^{\lambda+1}=\Gamma(\Gamma^\lambda)$, which is equivalent to
    \begin{equation*}
        \forall x\in \Seq(\Gamma^\lambda) \forall y [R(n,x,y)\to y\in \bbN\setminus \Seq(\bbN\setminus \Gamma^\lambda)].
    \end{equation*}
    Now for an admissible set $M$ of height $\lambda$, we claim that there is an $\Sigma_1^M$-definable normal function $f\colon\lambda\to\lambda$ such that for $\alpha<\lambda$,
    \begin{equation*}
        \forall x\in \Seq(\Gamma^\alpha) \forall y [R(n,x,y)\implies y\in \bbN\setminus\Seq(\bbN\setminus \Gamma^{f(\alpha)})].
    \end{equation*}
    Note that $\lag \Gamma^\xi\mid \xi<\lambda\rag$ is $\Sigma_1$-definable over $M$. For $x\in \Seq(\Gamma^\lambda)$ and $y\in \bbN$, define
    \begin{equation*}
        g(x,y) = 
        \begin{cases}
            \min \{\alpha\mid y\in \bbN \setminus\Seq(\bbN\setminus \Gamma^\alpha)\} & \text{If }R(n,x,y),\\
            0 & \text{otherwise.}
        \end{cases}
    \end{equation*}
    Then $g\colon \Seq(\Gamma^\lambda)\times \bbN\to \lambda$ is also $\Sigma_1$-definable over $M$.
    For $\alpha<\lambda$, define $h(\alpha)=\sup\{g(x,y)\mid x\in \Seq(\Gamma^\alpha)\land y\in \bbN\}$, and define $f$ by $\Sigma_1$-recursion over $M$ so that $f(0)=h(0)$, $f(\alpha+1)=f(\alpha)+h(\alpha)$, and $f(\delta)=\sup_{\alpha<\delta}f(\alpha)$ if $\delta$ limit.
    Then $f$ is \added{the} desired function.

    Since $f$ is normal and $\Sigma_1$-definable over $M$, $f$ has a fixed point $\beta<\lambda$. Hence we have
    \begin{equation*}
        \forall x\in \Seq(\Gamma^\beta) \forall y [R(n,x,y)\to y\in \bbN\setminus \Seq(\bbN\setminus \Gamma^\beta)],
    \end{equation*}
    so $n\in \Gamma^{\beta+1}\subseteq \Gamma^\lambda$.
\end{proof}

\begin{proof}[Proof of \autoref{Proposition: The closure point of a Pi01-Pi00 operator}]
    Suppose that $\Gamma$ is a $[\Pi^0_1,\Pi^0_0]$-operator given by a $\Pi^0_1$-operator $\Gamma_0$ and a $\Pi^0_0$-operator $\Gamma_1$.
    Modifying the proof of \autoref{Lemma: Admissible ordinal is a closure pt of a Pi01 operator} gives the following: If $\alpha$ is admissible, then $\Gamma(\Gamma^\alpha) = \Gamma_1(\Gamma^\alpha)$.
    Now let $\lambda$ be a recursively inaccessible ordinal. Suppose that $n\in \Gamma^{\lambda+1} = \Gamma_1(\Gamma^\lambda)$. 
    Since $\Gamma_1$ is (primitive) recursive, answering the query $n\in \Gamma_1(\Gamma^\lambda)$ uses finitely many questions about the membership questions on $\Gamma^\lambda$. In addition, $\lambda$ is limit, so we have that
    \begin{equation*}
        n\in \Gamma_1(\Gamma^\lambda)\iff \exists \xi<\lambda \forall \eta<\lambda [\eta\ge \xi \to n\in \Gamma_1(\Gamma^\eta)].
    \end{equation*}
    Since admissible ordinals are cofinal in $\lambda$, we can find an admissible $\eta\ge\xi$ for $\xi$ given in the previous statement.
    For such $\eta$, $n\in \Gamma_1(\Gamma^\eta) = \Gamma^{\eta+1} \subseteq \Gamma^\lambda$.
\end{proof}

Hence we have
\begin{corollary}
    For $\frakM^J = (M,|\cdot|)$, $|M|$ is the least recursively inaccessible ordinal.
\end{corollary}
Note that the proof of \autoref{Proposition: The closure point of a Pi01-Pi00 operator} (In fact, the modification of \autoref{Lemma: Admissible ordinal is a closure pt of a Pi01 operator}) also shows every admissible ordinal is $\frakM^J$-Richterian. 
From the remaining part of this section, let us prove the following:

\begin{theorem} \label{Theorem: The next admissible level and hyperjump}
    For $\frakM^J = (M,|\cdot|)$, if $x\in M$ then $\calO^{M_{|x|}}\equiv_\sfT M_{|\lag 2,x\rag|}$ uniformly in $x$.
\end{theorem}
Here $\calO^X$ is the set of codes for $X$-primitive recursive functions coding a well-order; i.e., set of $e$ such that $Z_e^X:=\{\lag x,y\rag \mid [e](\lag x,y\rag,X)=0\}$ is a well-order. By Kleene normal form theorem, $\calO^X$ is Turing equivalent to a universal $\Pi^1_1[X]$-set. 

\begin{proof}
    Let us first prove $M_{|\lag 2,x\rag|} \le_\sfm \calO^{M_{|x|}}$.
    We have two cases: If $x$ is not of the form $\lag 2,y\rag$, then $M_{|x|}$ is not closed under the Richter operator.
    By \autoref{Proposition: M-Richterian ordinal for J-operator}, $|\lag 2,x\rag|$ is the least $\frakM$-Richterian ordinal greater than $|x|$. 
    This means that $M_{|\lag 2,x\rag|}$ is the least fixed point of $\frakR$ containing $M_{|x|}$, so $M_{|\lag 2,x\rag|}$ is $\Pi^1_1[M_{|x|}]$.
    If $x$ is of the form $\lag 2,y\rag$, $M_{|x|}$ is closed under the Richter operator so the previous formula gives $M_{|x|}$ and not $M_{|\lag 2,x\rag|}$.
    However, $M_{|\lag 2,x\rag|}$ is the least fixed point of $\frakR$ containing $M_{|x|}\cup \{\lag 2,y\rag\mid y\in M_{|x|}\}$.
    Hence we still have $M_{|\lag 2,x\rag|}$ is $\Pi^1_1[M_{|x|}]$.
    
    For uniformity, let us find a $\Pi^1_1$-formula $\phi(X,n)$ equivalent to the following:
    \begin{quote}
        $n$ is in the least fixed point of $\frakR$ containing $X$, \emph{or if} $\frakR(X)\subseteq X$, then $n$ is in the least fixed point of $\frakR$ containing $X\cup \{\lag 2,y\rag\mid y\in X\}$.
    \end{quote}
    
    Hence we can provide the reduction as follows: By Kleene normal form theorem, we can find a primitive recursive code $e$ such that $\phi(X,n)$ iff $[e](\cdot, n,X)$ codes a \emph{strict} well-order (so reflexivity fails). By smn-theorem, we can find $e'$ such that $[e](x,n,X)=\{[e'](n)\}(x,X)$ for every $n,x\in\bbN$ and a real $X$.
    Then we have
    \begin{equation*}
        n\in M_{|\lag 2,x\rag|} \iff [e'](n) \in \calO^{M_{|x|}}.
    \end{equation*}
    For $\calO^{M_{|x|}}\le_\sfT M_{|\lag 2,x\rag|}$, we follow the proof in \cite[\S I.4]{Sacks2017higher}. 
    First, let us consider the primitive recursive function $h$, taking an index $e$ for an $X$-primitive recursive well-order $Z^X_e$ and $n\in\bbN$ returning an $X$-primitive recursive index $h(e,n)$, so it satisfies
    \begin{equation*}
        [h(e,n)](\lag u,v\rag, X) = 
        \begin{cases}
            0 & \text{if }[e](z,X)=0\text{ for }z=\lag u,v\rag,\lag u,n\rag,\lag v,n\rag, \\
            1 & \text{otherwise}.
        \end{cases}
    \end{equation*}
    We can find a \emph{primitive recursive} $h$ satisfying the above equation by effectively finding a primitive recursive code for the function
    \begin{equation*}
        e,n\mapsto (\lag u,v\rag, X \mapsto \min(\{1\}\cup \{[e](z,X) \mid z=\lag u,v\rag,\lag u,n\rag,\lag v,n\rag\})).
    \end{equation*}
    Intuitively, $h(e,n)$ is an index for an initial segment of $Z^X_e$ below $n$, if $n$ is in the field of the well-order. (Otherwise, $h(e,n)$ codes the empty order.)
    By Kleene's recursion theorem, there is a partial recursive function $f$ such that for every $e$ on which $f(e)$ is defined,
    \begin{equation*}
        f(e) = \sfsup(x, \bblambda n.f(h(e,n))).
    \end{equation*}
    and the left-hand side is defined iff the right-hand side is defined. 
    Now we claim $f$ is the desired reduction; i.e., 
    \begin{equation*}
        e \in \calO^{M_{|x|}} \iff f(e)\in M_{|\lag 2,x\rag|}.
    \end{equation*}
    We prove the left-to-right implication by induction on the ordertype of $Z^{M_{|x|}}_e$:
    Note that $Z^{M_{|x|}}_{h(e,n)}$ is a proper initial segment of $Z^{M_{|x|}}_e$ for every $n$. Hence by the inductive hypothesis, $f(h(e,n))\in M_{|\lag 2,x\rag|}$ for every $n\in\bbN$.
    Since $M_{|\lag 2,x\rag|}$ is closed under the Richter operator $\frakR$ and $\sfsup(x,\id)\in M_{|\lag 2,x\rag|}$, we have $f(e) = \sfsup(x, \bblambda n.f(h(e,n))) \in M_{|\lag 2,x\rag|}$.
    For the right-to-left implication, we apply induction on $|f(e)|$: If $f(e)\in M_{|\lag 2,x\rag|}$, then $f(h(e,n))\in M_{|\lag 2,x\rag|}$ and $|f(h(e,n))| < |f(e)|$ for every $n$.
    By the inductive hypothesis, $f(h(e,n))\in M_{|\lag 2,x\rag|}$ implies $Z^{M_{|x|}}_{h(e,n)}$ is a well-order for every $n$, so $Z^{M_{|x|}}_e$ is also a well-order.
\end{proof}

Now let us \emph{define} the iterated hyperjump of the empty set with the help of the $J$-hierarchy:

\begin{definition}
    Let $\alpha$ be an ordinal less than the least recursively inaccessible ordinal.
    For $x\in M$ such that $|x|$ is the $\alpha$th admissible ordinal\footnote{More precisely, the $\alpha$th ordinal that is either admissible or a limit of admissible ordinals. For technical convenience, we take 0 as the zeroth admissible ordinal.}, we define $\HJ^\alpha(\emptyset)$ by $M^*_{|x|}$.
\end{definition}

Let us analyze the previous arguments to see some of them carries over $\ACA_0$ with a minimal assumption:
\begin{remark} \label{Remark: How to define iterated HJs}
    We have a $\Pi^1_1$-formula defining the singleton $\{\HJ^\alpha(\emptyset)\}$ such that $\ACA_0$ proves there is at most one real $X$ satisfying the $\Pi^1_1$-formula.
    To get the desired $\Pi^1_1$-formula, let us fix a notation $x\in M$ such that $|x|=\omega^\CK_\alpha$, and modify the definition of $J$ as follows:
    \begin{equation*}
        n \in J_x(X) \iff n\in \frakR(X) \lor [\frakR(X)\subseteq \field(X) \land n\in \{\lag 2,y\rag\mid y\in \calF(X)\}\cup \{k\in\bbN\mid x\in \field X\}]
    \end{equation*}
    Then $J$-hierarchy and $J_x$-hierarchy agree upon level $|x|$, but the $|x|+1$th level of the $J_x$-hierarchy becomes trivial. Hence $\HJ^\alpha(\emptyset)$ is equal to the last non-trivial level of the $J_x$-hierarchy. Hence with the help of the terminologies in \autoref{Proposition: Uniqueness of the Gamma defined set over ACA0}, the following formula defines a singleton for $\HJ^\alpha(\emptyset)$:
    \begin{quote}
        For every well-order $\gamma$ with a maximal element $\mu$ and its immediate predecessor $\mu^-$, and an exact chain of $J_x$-inductive definition $A$ along $\gamma$, $(A)_\mu=\bbN \neq (A)_{\mu^-}$ and $X=(A)_{\mu^-}$.
    \end{quote}
    The previous formula is $\Sigma^1_2$ due to the universal quantifier over a well-order. Employing the trick in \autoref{Remark: Pi 1 1 singleton definition for an arithmetical ind def set}, we obatin a $\Pi^1_1$-formula $\phi_x(X)$ defining a singleton for $\HJ^\alpha(\emptyset)$.
    By the proof of \autoref{Proposition: Uniqueness of the Gamma defined set over ACA0}, $\ACA_0$ proves there is at most one $X$ satisfying $\phi_x(X)$.
\end{remark}

\begin{remark}
    We can see that the proof of \autoref{Theorem: The next admissible level and hyperjump} works over $\ACA_0$ with the assumption $M_{|\lag 2,x\rag|}$ exists. We may use elements of well-orders in place of ordinals at the beginning of the proof of \autoref{Theorem: The next admissible level and hyperjump}, and Kleene normal form theorem is a theorem of $\ACA_0$. The last part of the proof can be seen the application of arithmetical transfinite induction, which is a theorem of $\ACA_0$.
    \autoref{Proposition: M hierarchy and M star hierarchy are Turing equivalent} is also $\ACA_0$-provable.

    Now for a set $X$, let $\HJ(X)$ be the universal $\Pi^1_1[X]$-set. $\ACA_0$ does not prove it exists, but we can still express the definition of $\HJ(X)$ in the language of second-order arithmetic. Kleene normal form theorem says $\ACA_0$ proves $\calO^X\equiv_\sfT \HJ(X)$ if both exist, and in fact, the existence of one of $\calO^X$ or $\HJ(X)$ implies the existence of the other.
    In sum, $\ACA_0$ proves the following: If $M_{|\lag 2,x \rag|}$ exists, then $\HJ(M_{|x|})\equiv_\sfT M_{|\lag 2,x\rag|}$.
\end{remark}

\subsection{Iterated hyperjumps and iterated Spector classes}

The following theorem shows the connection between iterated hyperjumps of the empty set and the iterated Spector class over $\bbN$:
\begin{theorem} \label{Theorem: Successor Spector class as a Pi 1 1 R pointclass}
    For every $\xi$ less than the least recursively inaccessible ordinal, $\SP^{\xi+1}_\bbN = \Pi^1_1[\HJ^\xi(\emptyset)]$.
\end{theorem}
\begin{proof}
    The case $\xi=0$ is clear. Now let us prove the equality by induction on $\xi$. First, combining \cite[Theorem 6B.5, 8A.1]{Moschovakis1974EIAS} and \cite[Corollary 9A.3]{Moschovakis1974EIAS} gives that for each real $X$, the set of $\Pi^1_1[X]$-relations (equivalently, $\Pi^1_1$-relations over the structure $(\bbN;X)$) is precisely the least Spector class over $\bbN$ containing $X$ and $\bbN\setminus X$. We will use this fact in the remaining part of the proof.
    Hence if $\frakM$ is a Spector class over $\bbN$ such that $X,\bbN\setminus X\in \frakM$ and $\phi(x,X)$ is a $\Pi^1_1[X]$-formula, then the set $\{x\in\bbN\mid \phi(x,X)\}$ is in $\frakM$.
    Moreover, it is a standard fact that the class of $\Pi^1_1[X]$-sets is always a Spector class for a real $X$.

    Let us consider the case $\xi=\eta+1$ first. By the inductive hypothesis, $\SP^\xi_\bbN = \Pi^1_1[\HJ^\eta(\emptyset)]$, which means that every set in $\SP^\xi_\bbN$ is Turing reducible into $\HJ(\HJ^\eta(\emptyset)) \equiv_\sfT \HJ^\xi(\emptyset)$.
    This shows
    \begin{equation*}
        \SP^{\xi}_\bbN \subsetneq \SP^{\xi}_\bbN\cup \{\HJ^\xi(\emptyset),\bbN\setminus \HJ^\xi(\emptyset)\}\subseteq \Pi^1_1[\HJ^\xi(\emptyset)],
    \end{equation*}
    so $\SP^{\xi+1}_\bbN \subseteq \Pi^1_1[\HJ^\xi(\emptyset)]$.
    Conversely, if $\frakM$ is a Spector class over $\SP^\xi_\bbN$ and $\HJ^\xi(\emptyset), \bbN\setminus \HJ^\xi(\emptyset)\in \frakM$, then $\Pi^1_1[\HJ^\xi(\emptyset)]\subseteq \frakM$ by the minimality of $\Pi^1_1[\HJ^\xi(\emptyset)]$. This shows $\Pi^1_1[\HJ^\xi(\emptyset)]\subseteq \SP^{\xi+1}_\bbN$.

    Now consider the case when $\xi$ is limit. By \autoref{Lemma: Richterian Hierarchy Turing reducibility}, every $\HJ^\eta(\emptyset)$ is Turing reducible to $\HJ^\xi(\emptyset)$ for $\eta<\xi$. Hence by the inductive hypothesis, every set in $\SP^\xi_\bbN = \bigcup_{\eta<\xi}\SP^\eta_\bbN$ is Turing reducible to $\HJ^\xi(\emptyset)$. This shows $\SP^{\xi+1}_\bbN \subseteq \Pi^1_1[\HJ^\xi(\emptyset)]$.
    
    Conversely, suppose that $\frakM\supseteq \SP^\xi_\bbN$ is a Spector class over $\bbN$.
    We work with an admissible companion $A_\frakM$ (cf. \cite[\S 9E]{Moschovakis1974EIAS}) instead of $\frakM$, so $A_\frakM$ is an admissible set such that
    \begin{enumerate}
        \item $A_\frakM$ is the least admissible set containing $\Delta$-sets in $\frakM$; i.e., a set $X\in \frakM$ such that $\bbN\setminus X\in \frakM$.
        \item There is a relation $R_\frakM \in A_\frakM$ such that a relation $P\subseteq\bbN^n$ is in $\frakM$ iff $P$ is $\Sigma_1(R_\frakM)$ on $A_\frakM$.
        \item Every $\Delta$-set in $\frakM$ is rud definable in $R_\frakM$ over $A_\frakM$.
    \end{enumerate}
    Since $\SP^\xi_\bbN\subseteq \frakM$, we have $\SP^\xi_\bbN \subseteq A_\frakM$.
    Now fix $a,b\in \bbN$ such that $|\sfsup(a,b)|=\omega^\CK_\xi$, so the sequence $\lag |[b](n,M_{|a|})| : n\in\omega\rag$ is a cofinal sequence below $\omega^\CK_\xi$.
    Then observe that for each $n<\omega$, $x\in M^*_{|[b](n,M_{|a|})|}$ if and only if
    \begin{align*}
        A_\frakM \vDash \exists F \exists\gamma\in\Ord & 
        \textstyle  \Big[ F\colon (\gamma+2)\to \calP(\bbN)\land \forall \xi<\gamma+2 \Big(F(\xi) = \bigcup_{\eta<\xi} J_{[b](n,M_{|a|})}(F(\eta))\Big)
        \\& \textstyle F(\gamma+1)=\bbN \neq F(\gamma) \land x\in F(\gamma) \Big].
    \end{align*}
    
    Now consider the set $H$ such that $x\in H$ if and only if
    \begin{align*}
        A_\frakM \vDash \exists n\in\omega \exists F\in Y \exists\gamma\in\Ord\cap Y & 
        \textstyle  \Big[ F\colon (\gamma+2)\to \calP(\bbN)\land \forall \xi<\gamma+2 \Big(F(\xi) = \bigcup_{\eta<\xi} J_{[b](n,M_{|a|})}(F(\eta))\Big)
        \\& \textstyle F(\gamma+1)=\bbN \neq F(\gamma) \land x\in F(\gamma) \Big].
    \end{align*}
    
    By $\Sigma_1$-collection, the above formula is equivalent to a $\Sigma_1$-formula with parameter $R_\frakM$ originated from $M_{|a|}$.
    Hence $H$ is $\Sigma_1(R_\frakM)$-defined over $A_\frakM$, so $H\in \frakM$.
    Moreover, $H=\bigcup_{n<\omega} M_{|[b](n,M_{|a|})|}^* = M^*_{|\sfsup(a,b)|} = \HJ^\xi(\emptyset)$.
    Similarly, $x\notin H$ if and only if
    \begin{align*}
        A_\frakM \vDash \forall n\in\omega \exists F\in Y \exists\gamma\in\Ord\cap Y & 
        \textstyle  \Big[ F\colon (\gamma+2)\to \calP(\bbN)\land \forall \xi<\gamma+2 \Big(F(\xi) = \bigcup_{\eta<\xi} J_{[b](n,M_{|a|})}(F(\eta))\Big)
        \\& \textstyle F(\gamma+1)=\bbN \neq F(\gamma) \land x\in \bbN \setminus F(\gamma) \Big].
    \end{align*}
    The above formula is also equivalent to a $\Sigma_1(R_\frakM)$-formula, so $\bbN\setminus H\in \frakM$.
    Hence $\frakM$ is a Spector class containing $\HJ^\xi(\emptyset)$ and its complement, so $\Pi^1_1[\HJ^\xi(\emptyset)]\subseteq \frakM$.
\end{proof}

\begin{theorem} 
    For $\xi$ less than the least recursively inaccessible ordinal, $\HJ^\xi(\emptyset)$ is a $\Sigma^1_2$-singleton real; Furthermore, 
    \begin{enumerate}
        \item We have a $\Sigma^1_2$-formula $\phi(X)$ with a unique free variable $X$ such that $\ACA_0$ proves there is at most one $X$ satisfying $\phi(X)$.
        \item Every set in $\SP^{\xi}_\bbN$ is Turing reducible into $\HJ^\xi(\emptyset)$. Hence every definable set over $\SP^{\xi}_\bbN$ is definable over $(\bbN;\HJ^\xi(\emptyset))$, and vice versa.
        \item If $\xi$ is a successor, then $ \HJ^\xi(\emptyset)\in \SP^\xi_\bbN$.
    \end{enumerate}
\end{theorem}
\begin{proof}
    The first clause of the statement follows from \autoref{Proposition: Uniqueness of the Gamma defined set over ACA0} and \autoref{Remark: How to define iterated HJs}. The proof of the second clause of the statement is implicit in the proof of \autoref{Theorem: Successor Spector class as a Pi 1 1 R pointclass}.

    Now let us prove the third clause. Suppose $\xi=\eta+1$, then $\HJ^\eta(\emptyset),\bbN\setminus \HJ^\eta(\emptyset)\in \SP^\xi_\bbN$ and $\HJ^\xi(\emptyset)$ is defined by a monotone inductive definition by the Richter operator starting from $\HJ^\eta(\emptyset)$. Hence $\HJ^\xi(\emptyset)$ is in $\SP^\xi_\bbN$ since $\SP^\xi_\bbN$ is closed under monotone inductive definition whose initial set and its complement is in $\SP^\xi_\bbN$.
\end{proof}

Now let us prove \autoref{Proposition: Connecting Pohlers' ordinal and Pi 1 1 R PTO}. 
Before restating the statement of the proposition, let us recall that $T$ is a first-order theory sound to $\SP^\xi_\bbN$. We understand a set in $\SP^\xi_\bbN$ as a predicate, so the language of $T$ is the language of first-order arithmetic plus new predicates corresponding to sets in $\SP^\xi_\bbN$.
The author guesses the next proposition holds for every $\xi$ less than the least recursively inaccessible ordinal, but we prove it only for successor ordinals.
If $\xi$ is a limit, then $\HJ^\xi(\emptyset)$ is not in $\SP^\xi_\bbN$, causing technical complications.
\begin{proposition} 
    Let $\xi$ be a successor ordinal less than the least recursively inaccessible ordinal.
    Suppose that $T$ is a strictly acceptable axiomatization of $\SP^{\xi}_\bbN$ containing $\Th(\SP^\xi_\bbN)$.
    Then we have 
    \begin{equation*}
        \delta^{\SP^{\xi}_\bbN}(T) = \lvert\ACA_0 + \text{`$\HJ^\xi(\emptyset)$ exists'} +  T\restriction \HJ^\xi(\emptyset)\rvert_{\Pi^1_1[\HJ^\xi(\emptyset)]},
    \end{equation*}
    where $T\restriction R$ is the theory obtained from $T$ by restricting its language to that of first-order arithmetic with an extra unary predicate for $R$.
\end{proposition}
\begin{proof}
    For a notational convenience, let $H=\HJ^\xi(\emptyset)$. 
    We first claim that $\delta^{\SP^\xi_\bbN}(T) = \delta^{(\bbN;H)}(T\restriction H)$: 
    $\delta^{\SP^\xi_\bbN}(T) \ge \delta^{(\bbN;H)}(T\restriction H)$ is clear.
    For $\delta^{\SP^\xi_\bbN}(T) \le \delta^{(\bbN;H)}(T\restriction H)$, suppose that $\prec$ is an arithmetical-in-$\SP^\xi_\bbN$ well-order such that $T\vdash \TI(\prec,X)$.
    Now let us find arithmetical-in-$H$ definition of $\prec$: Suppose that $\sigma(x,y)$ is a sentence defining $\prec$.
    For each $A\in \SP^\xi_\bbN$, let us fix a partial $H$-recursive function $f_A$ witnessing $A\le_\sfT H$, so $f_X(x)=1$ iff $x\in A$. (If $A=H$, take the characteristic function of $H$.)
    Then consider the formula $\sigma'(x,y)$ obtained from $\sigma(x,y)$ by replacing every occurrence of $n\in A$ for $A\in \SP^\xi_\bbN$ with $f_A(n)=1$. 
    Then $\sigma$ and $\sigma'$ define the same binary relation, i.e., $\SP^\xi_\bbN$ satisfies $\forall x,y [\sigma(x,y)\lr \sigma'(x,y)]$.
    Since $T$ satisfies every true first-order statement over $\SP^\xi_\bbN$, we have that $T\vdash \forall x,y [\sigma(x,y)\lr \sigma'(x,y)]$.
    If $\prec'$ is a binary relation defined by $\sigma'$, then $\prec$ and $\prec'$ are the same and $T\vdash \TI(\prec, X)\lr \TI(\prec',X)$.
    This completes the proof of the desired inequality.

    Now let us claim $\delta^{(\bbN;H)}(T\restriction H) = |\ACA_0 + \text{`$H$ exists'} +  T\restriction H|_{\Pi^1_1[H]}$.
    For a technical convenience, let us consider the extension $\ACA_0[H]$ of $\ACA_0$: The language of $\ACA_0[H]_T$ is the language of second-order arithmetic with a constant set symbol $H$, and it comprises the following axioms:
    \begin{enumerate}
        \item Axioms of $\ACA_0$. 
        \item For the $\Pi^1_1$-formula $\phi(X)$ defining a $\Pi^1_1$ singleton $\{\HJ^\xi(\emptyset)\}$ presented in \autoref{Remark: How to define iterated HJs}, we have $\phi(H)$.
        \item Every theorem in $T$ in the language of the first-order arithmetic (possibly with a free variable, allowing to express pseudo $\Pi^1_1$-sentences) with a unary predicate $H$ corresponding to $\HJ^\xi(\emptyset)$.
    \end{enumerate}
    We can see that $\ACA_0[H]_T$ is a conservative extension of $\ACA_0$ + `$H$ exists' + $T\restriction H$.
    By an argument for the conservativity of $\ACA_0$ over $\mathsf{PA}$ (cf. \cite[Lemma IX.1.3]{Simpson2009}), $\ACA_0[H]_T$ is a conservative extension of $T\restriction H$.
    Hence 
    \begin{equation} \label{Formula: delta N H equality}
        \delta^{(\bbN;H)}(T\restriction H) = \sup \{\operatorname{opt}(\prec) \mid \prec \text{ is arithmetical-in-$H$ and } \ACA_0[H]_T \vdash \forall^1 X\ \TI(\prec,X)\}.
    \end{equation}
    $\ACA_0$ proves arithmetical transfinite induction on $\prec$ is equivalent to the well-foundness of $\prec$, so \autoref{Proposition: Pi 1 1 R PTO and arithmetical in R WO length} implies the right-hand-side of \eqref{Formula: delta N H equality} is equal to $|\ACA_0[H]_T|_{\Pi^1_1[H]}$.
\end{proof}

\printbibliography

\end{document}